\documentclass[reqno]{amsart}
\usepackage{eurosym}
\usepackage{amsmath}
\usepackage{amsfonts}
\usepackage{amssymb}
\usepackage{amsfonts}
\usepackage{epsfig}
\usepackage{caption}
\usepackage{subcaption}
\usepackage{multirow}
\usepackage{enumitem}
\usepackage{stackengine}
\stackMath

\newtheorem{theorem}{Theorem}

\newtheorem{example}[theorem]{Example}

\newtheorem{remark}[theorem]{Remark}

\newtheorem{fact}[theorem]{Fact}

\begin{document}

\title[Long-time relative error analysis] {Long-time relative error analysis for linear ordinary differential equations with perturbed initial value}
\author[S. Maset] {S.\ Maset \\
Dipartimento di Matematica, Informatica e Geoscienze \\
Universit\`{a} di Trieste \\
maset@units.it}

\begin{abstract}
		{We investigate the propagation of initial value perturbations along the solution of a linear ordinary differential equation \( y'(t) = Ay(t) \). This propagation is analyzed using the relative error rather than the absolute error. Our focus is on the long-term behavior of this relative error, which differs significantly from that of the absolute error. The present paper is a practical sequel to the theoretical papers \cite{M1,M2} on the long-time behavior of the relative error: it includes applicative examples and important issues not addressed in \cite{M1,M2}. In addition, the present paper shows that understanding the long-term behavior provides insights into the growth of the relative error over all times, not just at large times. Therefore, it represents a crucial and fundamental aspect of the conditioning of linear ordinary differential equations, with applications in, for example, non-normal dynamics.}
\end{abstract}
\maketitle



\section{Introduction}

Consider a linear $n-$dimensional Ordinary Differential Equation (ODE) 
\begin{equation}
	\left\{ 
	\begin{array}{l}
		y^{\prime }\left( t\right) =Ay\left( t\right) ,\ t\in\mathbb{R}, \\ 
		y\left( 0\right) =y_{0},
	\end{array}
	\right.  \label{ode}
\end{equation}
where $A\in \mathbb{R}^{n\times n}$ and $y_0,y(t)\in\mathbb{R}^n$. In the present paper, we are interested in understanding how a perturbation of the initial value $y_{0}$ is
propagated along the solution $y(t)=\mathrm{e}^{tA}y_0$ of (\ref{ode}) over
a long time interval.

The next fact A) is well known.
\begin{itemize}
	\item [A)] \emph{For a generic perturbation of $y_0$,
		the perturbation of $y(t)$  asymptotically (as $t\rightarrow +\infty$)  vanishes exponentially if the rightmost eigenvalues of $A$ have negative real part and asymptotically diverges exponentially if such eigenvalues have positive real part}.
\end{itemize}
Thus, if we are concerned about how large the perturbation can become for large $t$,
we should be reassured by knowing that the rightmost eigenvalues of $A$ have negative real part and truly concerned by knowing that these eigenvalues have positive real part.

However, \emph{for a generic $y_0$, also $y(t)$  asymptotically vanishes exponentially if the rightmost eigenvalues of $A$ have negative real part and asymptotically diverges exponentially if such eigenvalues have positive real part}. Thus, to know that the perturbation of the solution asymptotically vanishes, or diverges, exponentially does not help us understand if it is really significant, when compared to the solution. 

Comparing the perturbation of the solution to the solution itself means considering the \emph{relative error} of the perturbed solution. The well-known fact A) is a description of   the long-time behavior of the \emph{absolute error} of the perturbed solution. It is important also to give a similar description of the long-time behavior of the relative error. This description is fact B) on page 17.

Understanding the long-time behavior of the relative error of the perturbed solution is fundamental in real-world systems described by mathematical models based on linear ODEs (\ref{ode}), which require simulation through the integration of these ODEs in the presence of uncertainty in the initial value. This is particularly true when the solution becomes small or large, compared to the initial value, a situation where looking at the absolute error of the perturbed solution can have little significance. In Subsection \ref{examples} below, we examine two such models and emphasize the importance of considering the relative error of the perturbed solution in cases of uncertainty in the initial value.

Put differently, we are interested in the \emph{relative conditioning} of the
problem 
\begin{equation}
	y_0\mapsto \mathrm{e}^{tA}y_0  \label{due}
\end{equation}
for large $t$, i.e. we are interested in studying how a relative error in the input $y_0$ is propagated to the output $\mathrm{e}^{tA}y_0$ for large $t$. 

The present paper is a sequel to the theoretical papers \cite{M1,M2} on the long-time relative conditioning of the problem (\ref{due}) and it contains applications to real-world systems,  experimental tests, and other practical issues related to the results in \cite{M1,M2}.  

\subsection{Literature}

In the literature, the relative conditioning of the matrix exponential function is a well-studied topic (see, e.g., \cite{Levis1969}, \cite{van2006}, \cite
{Kagstrom1977}, \cite{Van1977}, \cite{moler2003}, \cite{Mohy2008}, \cite
{Zhu2008}, \cite{Mohy2011}, \cite{Ed}, and \cite{ALMOHY2017}). Some of these papers are mainly focused on
computational aspects and algorithms, while others consider how the relative conditioning of the problem $A\mapsto \mathrm{e}^{tA}$
depends on $t$. However, for this problem, a general characterization of the long-time behavior of the relative conditioning is still not known.

In contrast, the relative conditioning of the action of the matrix exponential $\mathrm{e}
^{tA}$ on a vector with respect to perturbations of this vector, i.e., the relative 
conditioning of the problem (\ref{due}),  has
received little attention. The reason could be that (\ref{due}) is perceived as an easy linear problem, without all the complications involved in the non-linearity of the matrix exponential function problem. However, once one examines the dynamics of the relative conditioning, the issue is no longer straightforward (see \cite{M1,M2}).

A study of how the conditioning of the problem (\ref{due}) depends on $t$ was given
in the paper \cite{Maset2018}, but the analysis was confined to the case of a normal matrix $A$. The paper \cite{M1} extends this study to a general complex linear ODE and \cite{M2} delves into the results of \cite{M1} for the real case.

\subsection{Two models}\label{examples}

In many cases, a mathematical model based on a linear ODE (\ref{ode}) is constructed after a careful selection of the parameters appearing in the entries of the matrix \( A \) by the model creators. As a result, the matrix \( A \) can be regarded as fixed and reliable. The model is then used repeatedly with different initial conditions, often by users who can have little control over the accuracy of these initial conditions. In this context, where the initial conditions are inputs to a fixed model, the analysis of the relative error of the perturbed solution can be of interest.

We show the importance of understanding the behavior of the relative error of the perturbed solution through two  mathematical models based on a linear ODE (\ref{ode}).

\subsubsection{Gross Domestic Product and National Debt} \label{twom}

The papers \cite{Dmitriev2019,Dmitriev2020} present a mathematical model for the Gross Domestic Product (GDP) and National Debt (ND) of a given country. This model is based on the linear system of ODEs
\begin{equation}
\left\{ 
\begin{array}{l}
Q^{\prime }\left( t\right) =a_{11} Q\left( t\right)
+a_{12} B\left( t\right) \\
B^{\prime }\left( t\right)=a_{21} Q\left( t\right)
+a_{22}B\left( t\right),
\end{array}
\right.  \label{ODEdebt}
\end{equation}
where $Q(t)$ is the GDP and $B(t)$ is the ND at the time $t$. The unit for time is  $1\ \mathrm{yr}$ and the unit for GDP and ND is the initial GDP, i.e., $Q(0)=1$. A possible instance for the coefficients (given in \cite{Dmitriev2019}) is
\begin{equation}
a_{11}=0.08,\ a_{12}=-0.07,\ a_{21}=0.03,\ a_{22}=-0.02. \label{instanceGDPND}
\end{equation}
For this instance, the  matrix $A$ of the ODE has the positive eigenvalues $0.05$ and $0.01$. 

Suppose we are interested in simulating the growth of GDP and ND over a period of $50\ \mathrm{yr}$ by integrating (\ref{ODEdebt}) for a given initial ND (remember that the initial GDP is set at $1$). We assume, as in \cite{Dmitriev2019}, that the initial ND is $0.60$. Due to uncertainty in the available data, this initial ND may not be the actual value. Suppose the actual initial ND is $0.61$. We set $B(0)$ to this actual initial ND of $0.61$, while our assumed initial value of $0.60$ for the simulation is the perturbed value $\widetilde{B}(0)$. (Conversely, we could set $B(0) = 0.60$ and $\widetilde{B}(0) = 0.61$. See Remark \ref{remark}).

Since the eigenvalues of $A$ are positive, the perturbation of the initial ND grows exponentially in the solution.
After $50\ \mathrm{yr}$, for the actual initial value $(Q(0),B(0))=(1,0.61)$, we have
$$
(Q(50),B(50))=(8.84,4.09)
$$
and, for the perturbed initial value $(\widetilde{Q}(0),\widetilde{B}(0))=(1,0.60)$, we have
$$
(\widetilde{Q}(50),\widetilde{B}(50))=(9.02,4.15). 
$$
Here, we denote the perturbed solution of (\ref{ODEdebt}) by $(\widetilde{Q}(t),\widetilde{B}(t))$. 

Now, the question is:
\begin{itemize}
\item [] after $50\ \mathrm{yr}$, is the effect of the initial uncertainty significant?
\end{itemize}
 There are two ways to answer this.
\begin{itemize}
\item [1)] The perturbation of the initial ND is $1$ \emph{second decimal figure} (from $61$ to $60$ second decimal figures) and the perturbations  of GDP and ND  after $50\ \mathrm{yr}$ are, respectively, $18$ \emph{second decimal figures} and $6$ \emph{second decimal figures} (from $884$ to $902$  for GPD and from $409$ to $415$ for ND), respectively. This demonstrates a significant increase of the perturbation, as an effect of exponential growth. In Figure \ref{QB}, we see in the plane $(Q,B)$ the points $(Q(0),B(0))$ and $(\widetilde{Q}(0),\widetilde{B}(0))$ on the left and the points $(Q(50),B(50))$ and $(\widetilde{Q}(50),\widetilde{B}(50))$ on the right. On the left and on the right, we use the same scale with both axes of length 1 GDP unit.     
\begin{figure}[tbp]
\includegraphics[width=0.8\textwidth]{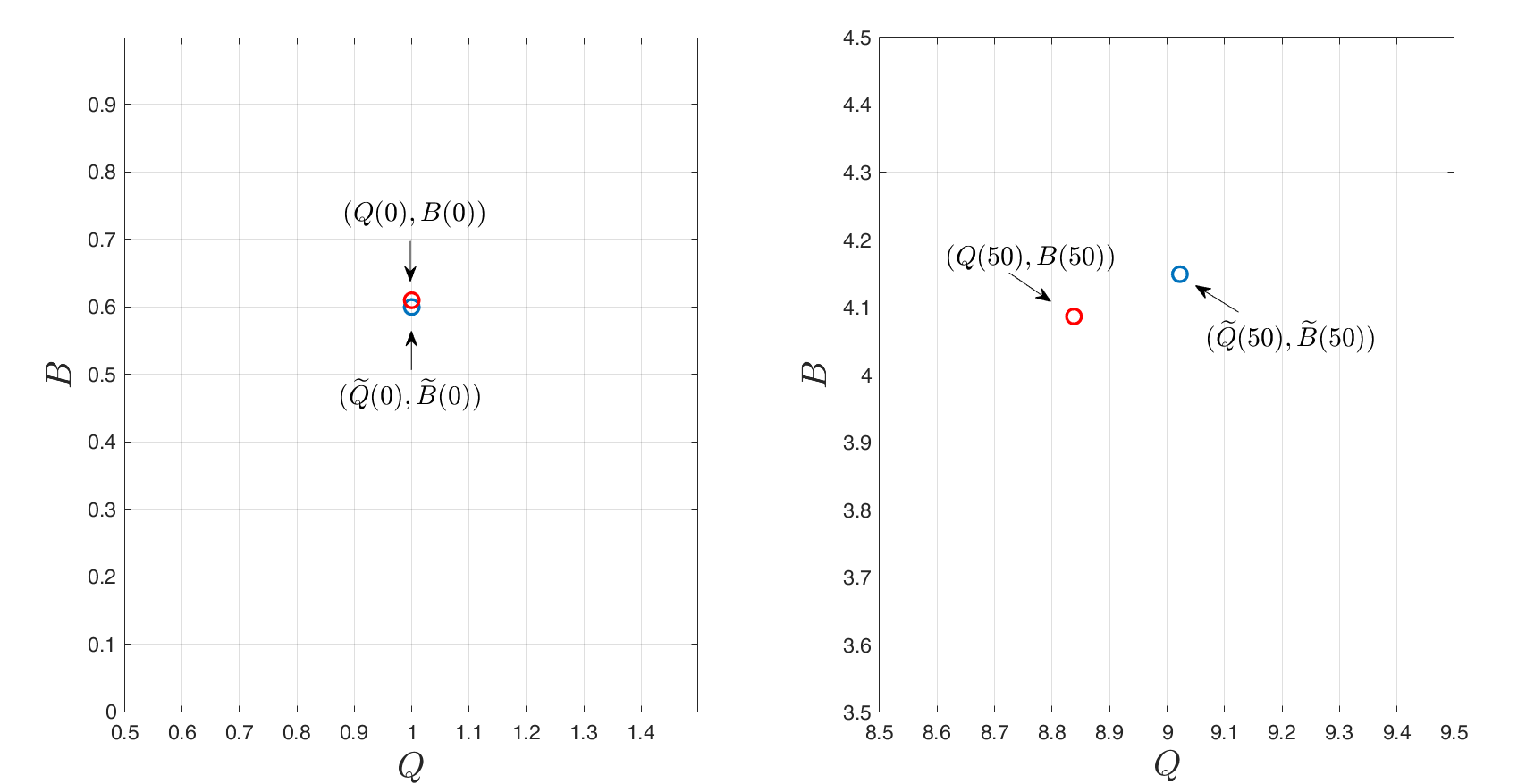}
\caption{Left: points $(Q(0),B(0))$ and $(\widetilde{Q}(0),\widetilde{B}(0))$. Right: points $(Q(50),B(50))$ and $(\widetilde{Q}(50),\widetilde{B}(50))$.}
\label{QB} 
\end{figure}
The growth of the perturbation is evident. In  Figure \ref{Figura4}, we see the absolute error of the perturbed solution
\begin{equation}
e(t)=\left\Vert (\widetilde{Q}\left(t\right),\widetilde{B}(t))-\left(Q\left(t\right),B(t)\right)\right\Vert_2 \label{abs},
\end{equation}
i.e., the Euclidean distance between the points $(\widetilde{Q}\left(t\right),\widetilde{B}(t))$ and $\left(Q\left(t\right),B(t)\right)$, for $t\in[0,50]$.
\begin{figure}[tbp]
\includegraphics[width=0.8\textwidth]{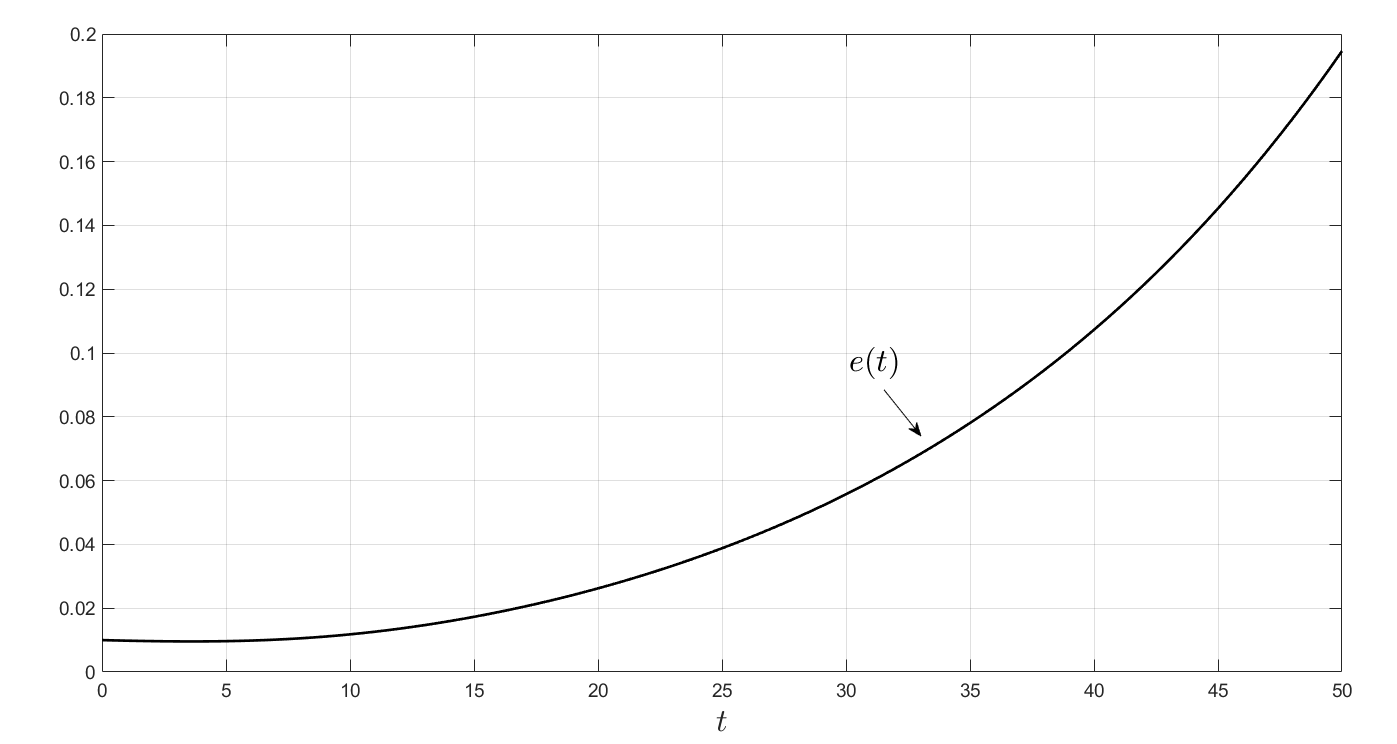}
\caption{Absolute error $e(t)$ for the GDP-ND model.}
\label{Figura4} 
\end{figure} 
After $50\ \mathrm{yr}$, the absolute error is about $20$ times the initial absolute error.

\item [2)] The perturbation of the initial ND is $1$ \emph{second significant figure} (from $61$ to $60$ second significant figures) and the perturbation  of GDP and ND after $50\ \mathrm{yr}$ is $1.8$ \emph{second significant figures} and $0.6$ \emph{second significant figures} (from $88.4$ to $90.2$ and from $40.9$ to $41.5$), respectively. This shows that the perturbed GDP and ND are still close to the unperturbed values after $50\ \mathrm{yr}$, by taking into account the growth in magnitude of GDP and ND. In Figure \ref{QB2}, similarly to  Figure \ref{QB}, we see the points $(Q(0),B(0))$ and $(\widetilde{Q}(0),\widetilde{B}(0))$ on the left and the points $(Q(50),B(50))$ and $(\widetilde{Q}(50),\widetilde{B}(50))$ on the right, but now on the left  we use a scale with axes of length $\left\Vert (Q(0),B(0))\right\Vert_2$, and on the right we use a scale with axes of length $\left\Vert (Q(50),B(50))\right\Vert_2$. 
\begin{figure}[tbp] 
\includegraphics[width=0.8\textwidth]{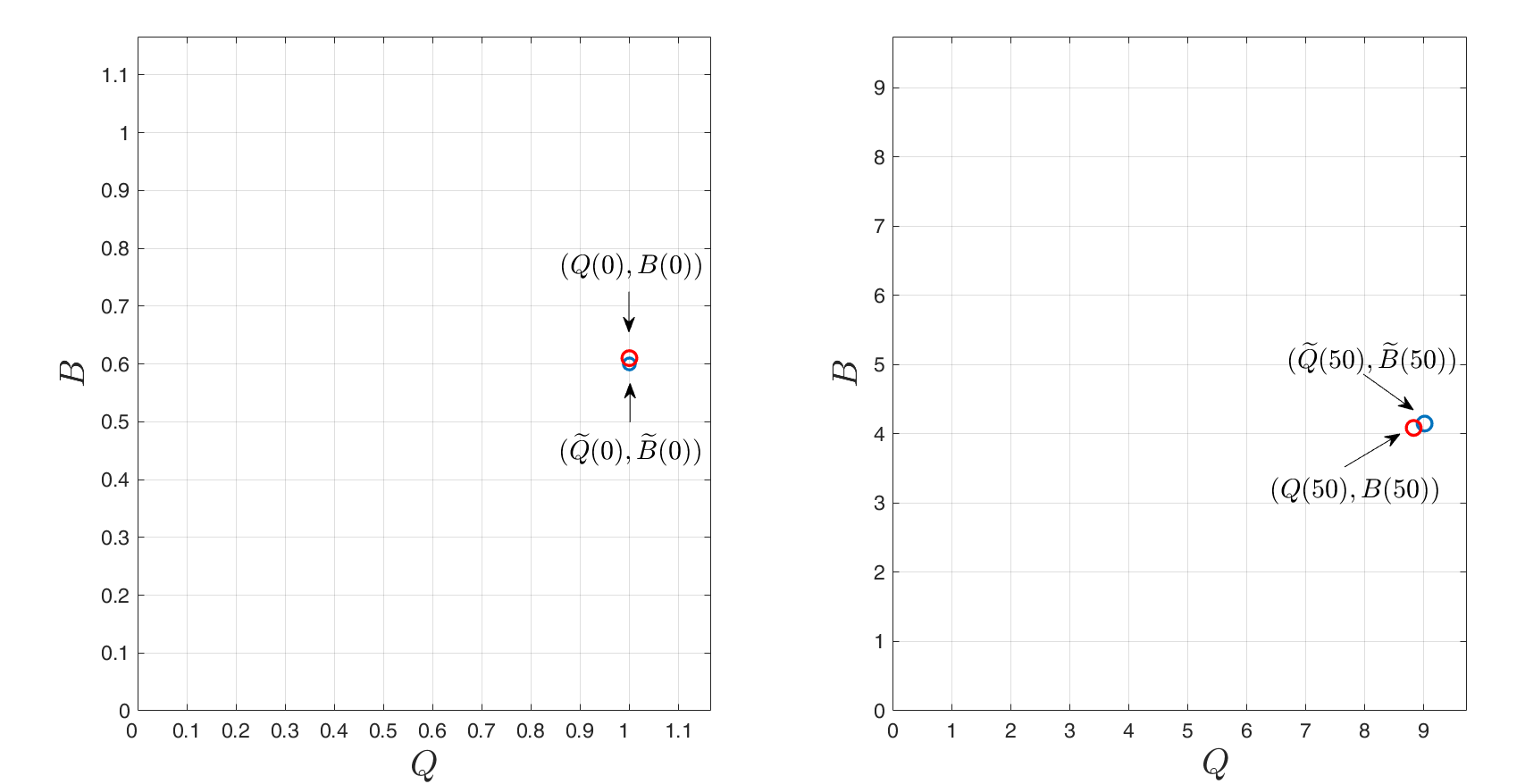}
\caption{Left: points $(Q(0),B(0))$ and $(\widetilde{Q}(0),\widetilde{B}(0))$. Right: points $(Q(50),B(50))$ and $(\widetilde{Q}(50),\widetilde{B}(50))$.}
\label{QB2} 
\end{figure}
The closeness of $(Q(50),B(50))$ and $(\widetilde{Q}(50),\widetilde{B}(50))$ is comparable to the closeness of $(Q(0),B(0))$ and $(\widetilde{Q}(0),\widetilde{B}(0))$. This closeness is also confirmed by Figure \ref{Figura3}, where it appears that GDP and perturbed GDP, as well as ND and perturbed ND, plotted as functions of time over $50$ $\mathrm{yr}$ are barely distinguishable.
\begin{figure}[tbp]
\includegraphics[width=0.8\textwidth]{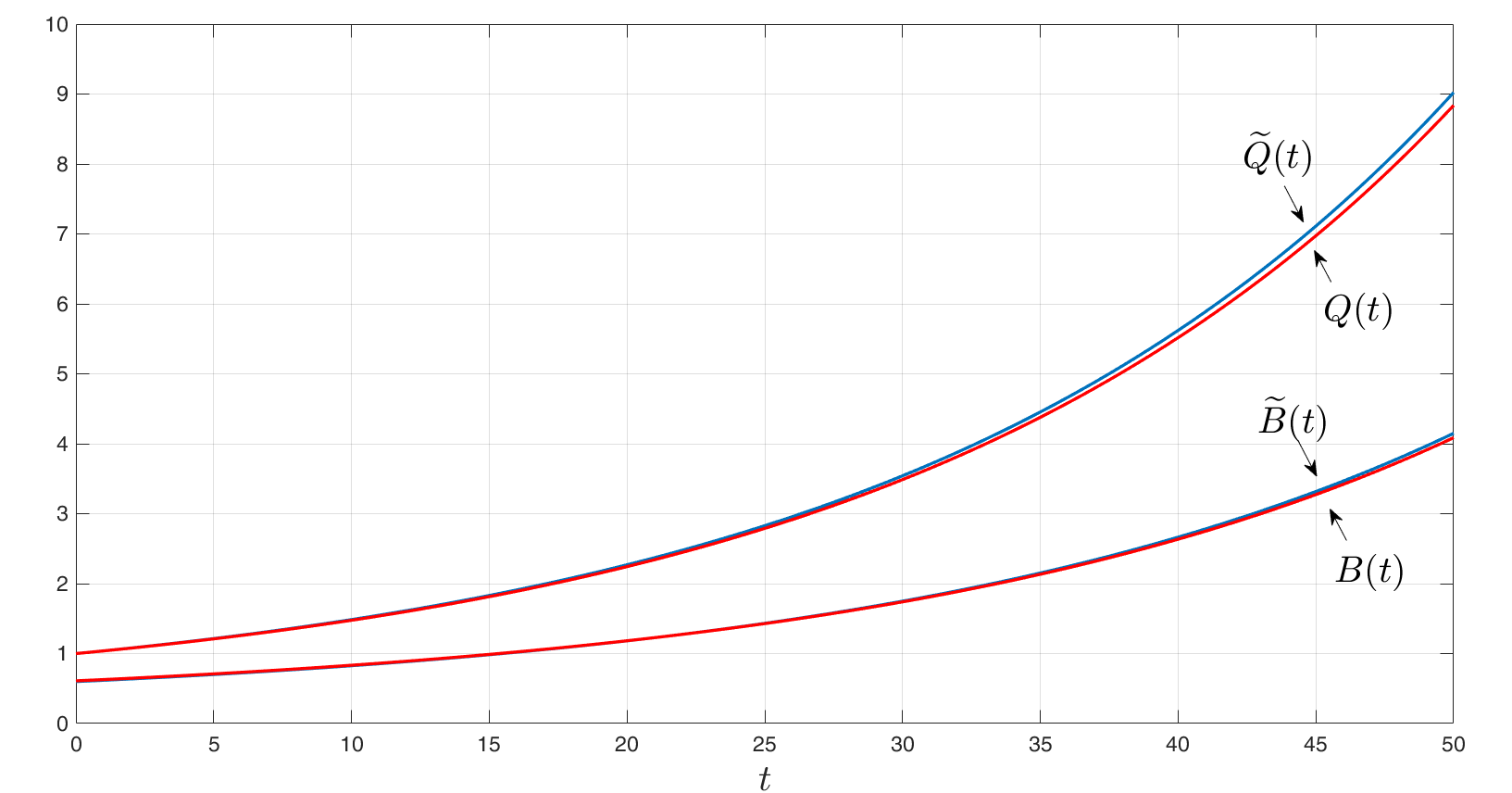}
\caption{GDP $Q(t)$, perturbed GDP $\widetilde{Q}(t)$, ND $B(t)$ and perturbed ND $\widetilde{B}(t)$.}
\label{Figura3} 
\end{figure}
In Figure  $\ref{Figura5}$, we see the relative error of the perturbed solution
\begin{equation}
\delta(t)=\frac{\left\Vert (\widetilde{Q}\left(t\right),\widetilde{B}(t))-\left(Q\left(t\right),B(t)\right)\right\Vert_2}{\left\Vert\left(Q\left(t\right),B(t)\right)\right\Vert_2} \label{rel}
\end{equation}
for $t\in[0,50]$. It is the distance between the points $(\widetilde{Q}\left(t\right),\widetilde{B}(t))$ and $\left(Q\left(t\right),B(t)\right)$ when, as in Figure \ref{QB2}, we use the length $\left\Vert (Q(t),B(t))\right\Vert_2$ of the axes as unit length.
\begin{figure}[tbp]
\includegraphics[width=0.8\textwidth]{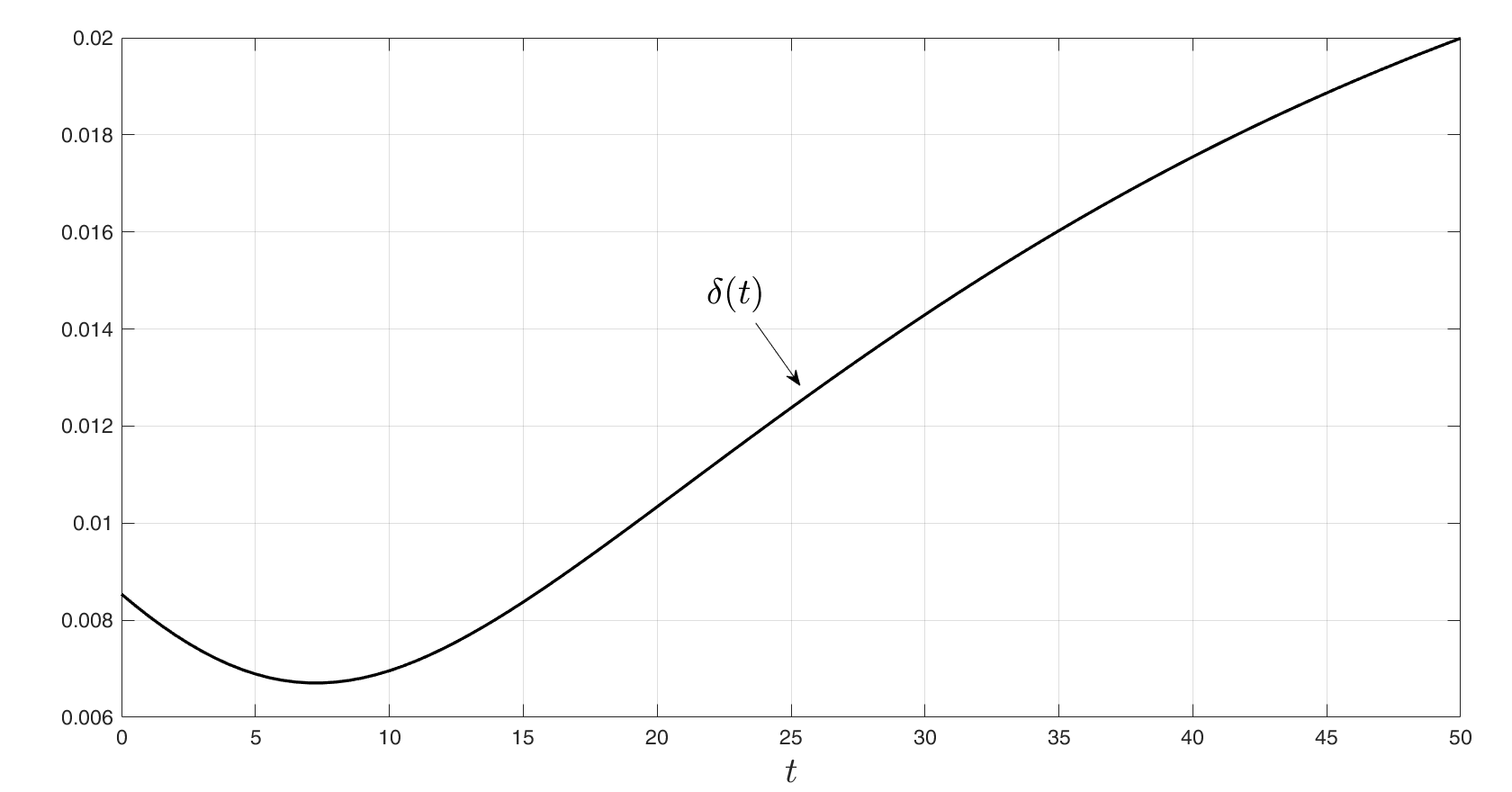}
\caption{Relative error $\delta(t)$ for the GDP-ND model.}
\label{Figura5} 
\end{figure}
Unlike the exponentially growing absolute error (\ref{abs}), the relative error (\ref{rel}) grows much less. After $50\ \mathrm{yr}$, it is about $2$ times the initial relative error.
\end{itemize}

If one believes that answer 2) is more appropriate (meaning that after 50 years, the effect of the initial perturbation is not significant because the perturbation of the solution is not much larger than its initial value, when compared to the solution), it becomes crucial to understand the behavior of the relative error (\ref{rel}) rather than that of the absolute error (\ref{abs}).

The relative error (\ref{rel}) refers to the relative error of the solution vector $(Q(t),B(t))$, i.e., it is a normwise relative error. One might find it more interesting to consider the relative errors of the two components $Q(t)$ and $B(t)$.
In the specific integration over $50$ $\mathrm{yr}$  we are considering, the relative errors of the components do not have order of magnitude larger than that of the relative error (\ref{rel}), since the components of the solution are not small compared to the Euclidean norm of the solution: we see in Figure \ref{QB2} that the point $(Q(50),B(50))$ is not close to the axes. It is noteworthy that in Figure \ref{Figura3}, $Q(t)$ and $\widetilde{Q}(t)$, as well as $B(t)$ and $\widetilde{B}(t)$, are nearly indistinguishable due to  small componentwise relative errors.

About componentwise relative errors, see Subsection \ref{nwcw} below.

This mathematical model of GDP and ND will be revisited in Section \ref{Examples} in light of the results presented in this paper.

\subsubsection{Building heating}\label{buildingheating}

The example in this subsection is a toy mathematical model of building heating (see the book \cite{Gust2022}).
Similar more complex non-toy models are used in literature for modeling real building heating (e.g., see \cite{MH1995}, \cite{Buz1998} and \cite{Tol2023}).

Consider a building constituted by
basement, main floor and attic. Let $x_1(t)$, $x_2(t)$ and $x_3(t)$ be the temperatures of basement, main floor
and attic, respectively, at the time $t$. By using \emph{Newton's law of cooling}
\begin{eqnarray*}
&&\text{rate of change of internal temperature } \\
&&\propto \text{external temperature $-$ internal temperature,}
\end{eqnarray*}
we derive the linear system of ODEs
\begin{equation}
\left\{ 
\begin{array}{l}
x_1^{\prime }\left( t\right) =k_{g1}\left( x_{g}-x_1\left( t\right) \right)
+k_{12}\left( x_2\left( t\right) -x_1\left( t\right) \right)+f_1\\
x_2^{\prime }\left( t\right) =k_{a2}\left( x_{a}-x_2\left( t\right) \right)
+k_{12}\left( x_1\left( t\right) -x_2\left( t\right) \right) +k_{23}\left(
x_3\left( t\right) -x_2\left( t\right) \right) +f_{2} \\
x_3^{\prime }\left( t\right) =k_{a3}\left( x_{a}-x_3\left( t\right) \right)
+k_{23}\left( x_2\left( t\right) -x_3\left( t\right) \right) +f_{3},
\end{array}
\right.  \label{building-heating}
\end{equation}
where
\begin{itemize}
\item $x_{g}$ is temperature of the ground and $x_{a}$ is the outdoor
temperature of the air;
\item $k_{g1},k_{a2},k_{a3},k_{12},k_{23}$ are positive proportionality constants in the Newton's law depending on the thermal insulation of the floors;
\item $f_{1}$, $f_{2}$ and $f_{3}$ are constant forcing terms due to heaters in the basement, main floor and
attic, respectively.
\end{itemize}

The ODE (\ref{building-heating}) has the form
\begin{equation}
x^\prime(t)=Ax(t)+b, \label{ODEex}
\end{equation}
where
$$
A=
\left[
\begin{array}{ccc}
-k_{g1}-k_{12} & k_{12} & 0  \\
k_{12} &  -k_{a2}-k_{12}-k_{23} & k_{23}  \\
0 & k_{23} &  -k_{a3}-k_{23}
\end{array}
\right]\text{\ \ and\ \ } 
b=\left[
\begin{array}{ccc}
k_{g1}x_g+f_1\\
k_{a2}x_a+f_2\\
k_{a3}x_a+f_3
\end{array}
\right].
$$
Since the matrix $A$ is symmetric and strictly diagonally dominant with negative diagonal entries, it has real negative eigenvalues. Consequently, the ODE (\ref{ODEex}) has a globally asymptotically stable equilibrium point
\begin{equation*}
x_{eq}=-A^{-1}b.\label{eq}
\end{equation*}

The transient 
$$
y(t)=x(t)-x_{eq}
$$
satisfies the ODE
\begin{equation}
y^\prime(t)=Ay(t). \label{ODEhouse}
\end{equation}
Suppose we want to simulate the transient by integrating (\ref{ODEhouse}), and there is uncertainty in the initial temperatures $y_0$ at time $t=0$.  Of course, $y_0$ and $y(t)$ are temperatures with respect to the equilibrium temperatures.

Therefore, we have an initial value $y_0$, representing the actual initial temperatures, and a perturbed initial value $\widetilde{y}_0$, representing the initial temperatures available to us for simulation. (Alternatively, we could assume that $y_0$ represents the available temperatures and $\widetilde{y}_0$ represents the actual temperatures. See Remark \ref{remark}). 

\begin{figure}[tbp]
	\includegraphics[width=0.8\textwidth]{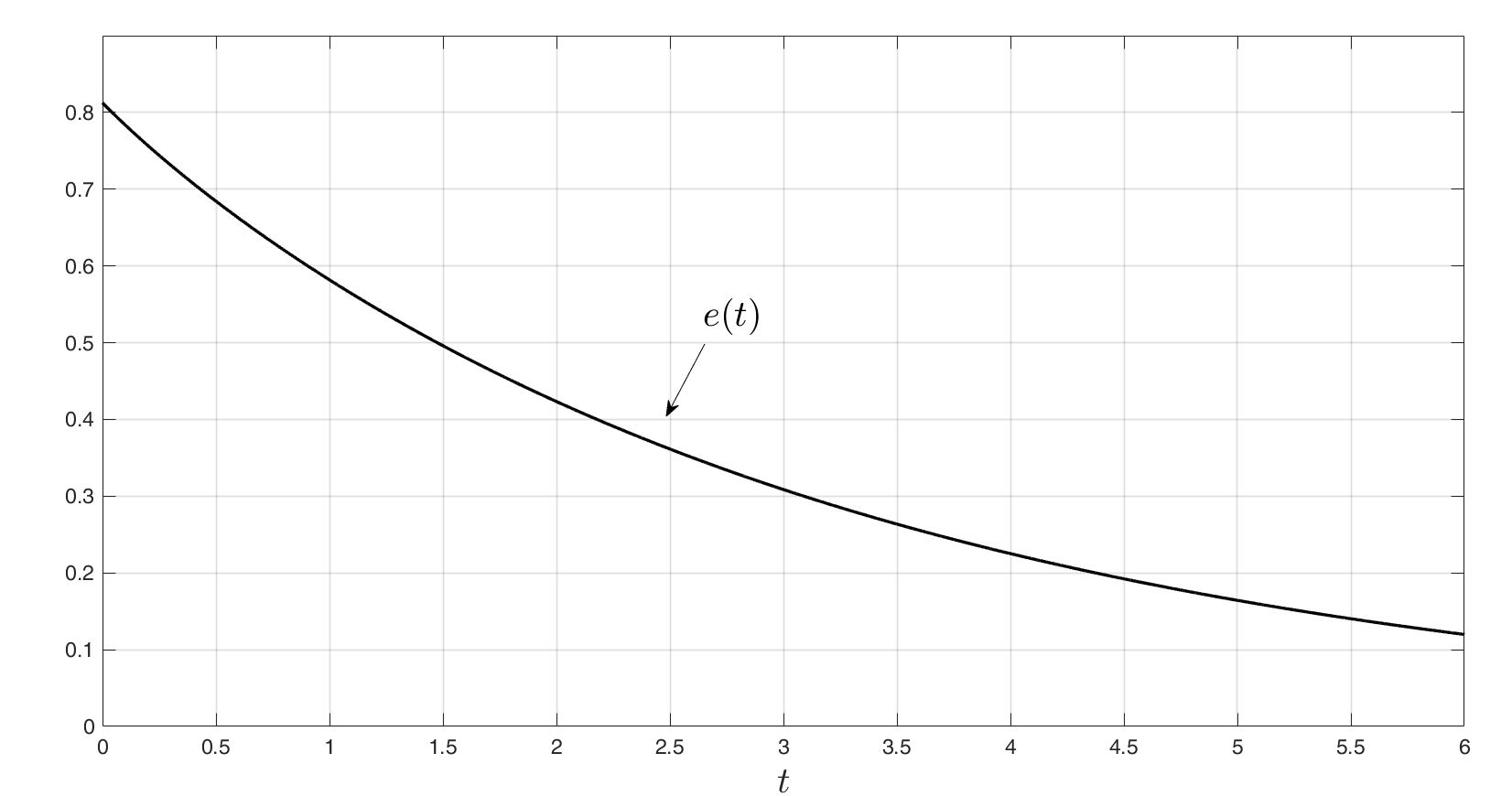}
	\caption{Absolute error $e(t)$ for the building heating model.}
	\label{absolutebh} 
\end{figure}

\begin{figure}[tbp]
	\includegraphics[width=0.8\textwidth]{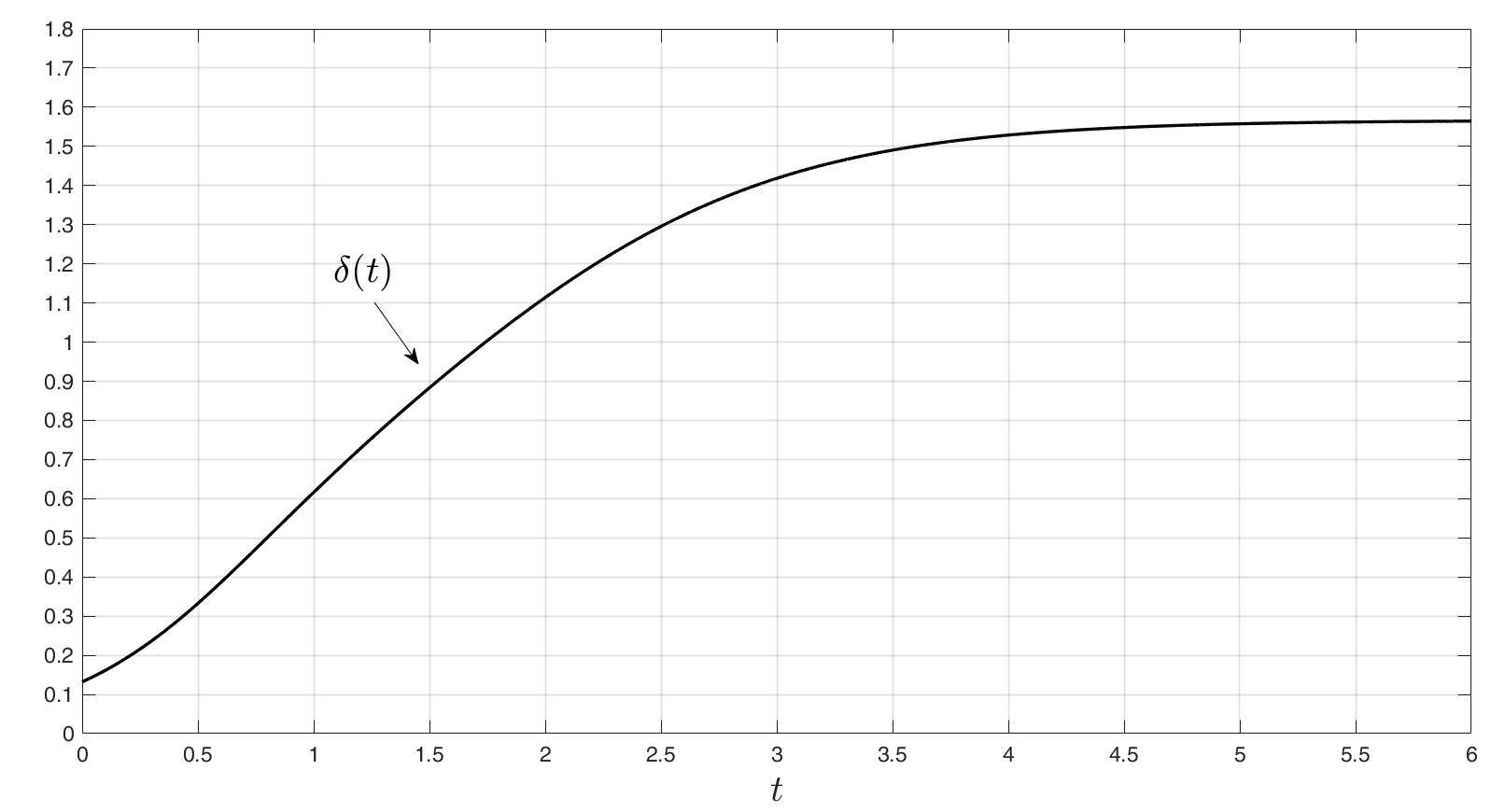}
	\caption{Relative error $\delta(t)$ for the building heating model.}
	\label{relativebh} 
\end{figure}

For illustrating our considerations, in (\ref{building-heating}) set
\begin{equation}
	k_{g1}=0.5\ \frac{\frac{^{\circ}\mathrm{C}}{\mathrm{h}}}{^{\circ}\mathrm{C}},\ k_{o2}=0.25\ \frac{\frac{^{\circ}\mathrm{C}}{\mathrm{h}}}{^{\circ}\mathrm{C}},
	k_{o3}=0.25\ \frac{\frac{^{\circ}\mathrm{C}}{\mathrm{h}}}{^{\circ}\mathrm{C}},\ k_{12}=0.5\ \frac{\frac{^{\circ}\mathrm{C}}{\mathrm{h}}}{^{\circ}\mathrm{C}},\ k_{23}=1\ \frac{\frac{^{\circ}\mathrm{C}}{\mathrm{h}}}{^{\circ}\mathrm{C}} \label{instancehb}
\end{equation}
(temperatures are measured in Celsius degrees and the time in hours). For this particular instance, the eigenvalues of $A$ are $-0.31519\ \mathrm{h}^{-1}$, $-1.0560\ \mathrm{h}^{-1}$ and $-2.6288\ \mathrm{h}^{-1}$. 
Suppose
\begin{equation}
	y_0=\left(3.5\ ^{\circ}\mathrm{C}, -4.4\ ^{\circ}\mathrm{C}, 2.5\ ^{\circ}\mathrm{C}\right) \label{y0pert}
\end{equation}
and   
\begin{equation}
	\widetilde{y}_0=\left(4\ ^{\circ}\mathrm{C}, -4\ ^{\circ}\mathrm{C},3\ ^{\circ}\mathrm{C}\right).\label{y0}
\end{equation}
In Figures \ref{absolutebh} and  \ref{relativebh} we show the absolute error
$$
e(t)=\left\Vert \widetilde{y}(t) -y(t)\right\Vert_2
$$
and the relative error
$$
\delta(t)=\frac{\left\Vert \widetilde{y}(t) -y(t)\right\Vert_2}{\left\Vert y(t)\right\Vert_2},
$$
where $\widetilde{y}$ is the perturbed (simulated) solution of (\ref{ODEhouse}), for $t\in \left[0,6\ \mathrm{h}\right]$.

Here, the scenario is completely reversed compared to the GDP-ND model. Due to the negative eigenvalues, the absolute error decreases exponentially and the initial perturbation is reduced by about seven times in $6$ hours, reassuring us about the effect of the uncertainty in the initial temperatures, when the transient is simulated by integrating (\ref{ODEhouse}). However, this is true only when comparing the perturbation to the initial value $y_0$, i.e., by looking at the absolute error. When comparing the perturbation to the solution $y(t)$, i.e., by looking at the relative error, the situation is much more concerning, as the initial perturbation is magnified by about twelve times. This could significantly impact our understanding of the transient and the decisions we make based on the simulation.

For example, suppose we want to determine when the Euclidean norm of the solution drops below $0.5\ ^{\circ}\mathrm{C}$. When this occurs, we could decide that the transient phase has ended. The simulation with initial value $\widetilde{y}_0$ shows that this occurs at time $\widetilde{t}^{\ast} = 3.0876\ \mathrm{h}$. However, at this time $\widetilde{t}^{\ast}$, the relative error of the perturbed solution is greater than 1 (see Figure \ref{relativebh}), meaning that the norm of the perturbation is larger than the norm of the actual solution. This implies that the norm $\Vert y(\widetilde{t}^\ast)\Vert_2$ of the actual solution  could be less than half the norm $\Vert \widetilde{y}(\widetilde{t}^\ast)\Vert_2=0.5^\circ\mathrm{C}$ of the perturbed solution. Indeed, we have
$\Vert y(\widetilde{t}^\ast)\Vert_2=0.2093\ ^\circ\mathrm{C}$. The norm of the actual solution dropped below $0.5\ ^\circ\mathrm{C}$ at \(t^\ast = 1.6362\ \mathrm{h}\), in roughly half the time. Therefore, the simulation with the initial temperatures $\widetilde{y}_0$ available to us provides a highly inaccurate estimate of the end transient time. Observe that the relative error $\delta(0)=0.1320$ of $\widetilde{y}_0$ with respect to $y_0$ is amplified into a relative error of $0.8871$ of $\widetilde{t}^\ast$ with respect to $t^\ast$.

Note that, to understand how different \(\widetilde{t}^\ast\) and \({t}^\ast\) can be, it is natural to consider the relative error $\delta(t)$. In fact, we have
$$
\frac{\Vert \widetilde{y}(t) \Vert_2}{\Vert y(t) \Vert_2}=1+\xi(t),
$$
where $\xi(t)$ satisfies $\vert \xi(t)\vert \leq \delta(t)$. Consequently, $t^\ast$ and $\widetilde{t}^\ast$ are such that
\begin{equation}
\Vert y(t^\ast) \Vert_2=0.5\ ^\circ\mathrm{C}\text{\ and\ }\Vert y(\widetilde{t}^\ast) \Vert_2=\frac{0.5\ ^\circ\mathrm{C}}{1+\xi(\widetilde{t}^\ast)}, \label{twotimes}
\end{equation}
where $\left\vert \xi(\widetilde{t}^\ast)\right\vert\leq \delta(\widetilde{t}^\ast)$.

Similarly to the previous GDP-ND model, we will revisit this building heating model in Section \ref{Examples}.

\subsection{The asymptotic behavior of the relative error}

The present paper deals with the \emph{asymptotic behavior}, i.e., the long-time behavior, of the relative error of the perturbed solution.

One could observe that, in case of a solution that decays to zero or diverges, this analysis might have limited relevance since asymptotically the solution is zero or infinite and then a long-time simulation of the solution is not very interesting. However, in such a case we might be interested in simulating the solution only up to a certain size threshold, beyond which it has become too small or too large to be of further interest. Therefore, if the relative error becomes close to, or of the same order of magnitude as, its asymptotic behavior before the solution has reached this threshold, then the analysis of the asymptotic behavior has interest.

In the GDP-ND and building heating models, although the solution diverges to infinity or converges to zero, we are still interested in understanding its evolution before it becomes too large or too small. The relative error after 50 years in the GDP-ND model and the relative error at the end transient time $\widetilde{t}^\ast$ in the building heating model have the same order of magnitude as their asymptotic values: see Figure \ref{Figura55} and Figure \ref{relativebh}.

In the  present paper, in addition to studying the asymptotic behavior of the relative error, we also investigate how rapidly this asymptotic behavior is attained. As expected, the non-normality of the matrix $A$ adversely affects the rapid attainment of the asymptotic behavior. However, as it will be shown, a high non-normality of $A$ does not necessarily imply a late onset of the asymptotic behavior, which could lead to a loss of interest in such behavior.

In conclusion, we can say that although the asymptotic behavior of the relative error does not fully describe the propagation to the solution of the perturbation of the initial value, since it can miss a possible initial growth  of the relative error to values much larger than the asymptotic behavior, it nonetheless constitutes an important piece in the qualitative study of the relative conditioning of the problem  (\ref{due}).

However, Section 7 shows experimentally that it is very rare to have such a large initial growth of the relative error. Therefore, \emph{the asymptotic behavior might be not only an important piece of the study of conditioning, but the crucial and fundamental piece}.

\subsection{Plan of the paper}
The asymptotic analysis of the relative error of the perturbed solution given in the papers \cite{M1,M2} is long and the details quite technical. This is due to the fact that \cite{M1} addressed the general case involving an arbitrary matrix $A$ in (\ref{ode}), a choice that necessitates dealing with the Jordan Canonical Form of $A$ and generalized eigenvectors. Moreover, the presence of complex eigenvalues further complicates matters in \cite{M2}. Indeed, by considering only diagonalizable matrices with real eigenvalues, the analysis would be considerably shorter. Finally, the definition of asymptotic form used in \cite{M1} also implies some effort in proving asymptotic results.

Due to their lengths, the papers \cite{M1,M2} only include theoretical results. All the non-theoretical practical issues are moved to the present paper, which is organized as follows.

Section 2 introduces two condition numbers presented in \cite{M1} for the problem (\ref{due}). Section 3 introduces notations and definitions, already used in \cite{M1} and \cite{M2}, for the description and the analysis of the asymptotic behavior of these two condition numbers. Section 4 recalls the results of \cite{M2} regarding the asymptotic condition numbers, i.e., the asymptotic forms of the condition numbers. Section 5 focuses on the more important of the two condition numbers and addresses fundamental questions such as asymptotic well-conditioning, the onset of asymptotic behavior, and the effect of the non-normality of the matrix $A$. In  Section 6, the GDP-ND and building heating models are revisited and three additional examples illustrate the contents of the previous sections. Section 7 shows that the asymptotic behavior of the relative error can also provide insight into the non-asymptotic behavior in most cases. Section 8 illustrates how the results of the present paper can be applied in non-normal dynamics. Section 9 presents the conclusions and can even be read now to gain a better idea of the subject of this paper.

The paper presents several numerical experiments in which a large number of random instances of the ODE (\ref{ode}) are generated by sampling the entries of \( A \) and the components of \( y_0 \) from the standard normal distribution. This means that all entries of \( A \) and all components of \( y_0 \) in every instance are independently sampled. It is worth noting that the entries of \( A \) are drawn from the standard normal distribution rather than from the normal distribution with mean zero and standard deviation \( \frac{1}{\sqrt{n}} \), which is often employed. We choose not to use the \( \frac{1}{\sqrt{n}} \) scaling because we aim to include a broad range of values for the norm of the matrix \( A \), thereby generating potentially stiffer or more extreme instances of the ODE. Nevertheless, performing the experiments with the scaled version yields essentially the same results and, consequently, leads to the same conclusions.

All numerical experiments are carried out in MATLAB, with the matrix exponential values $\mathrm{e}^{tA}$ computed using the \texttt{expm} function.

\section{The condition numbers}

\begin{figure}[tbp]
\centering
\begin{subfigure}[b]{1\textwidth}
         \centering
         \includegraphics[width=0.4\textwidth]{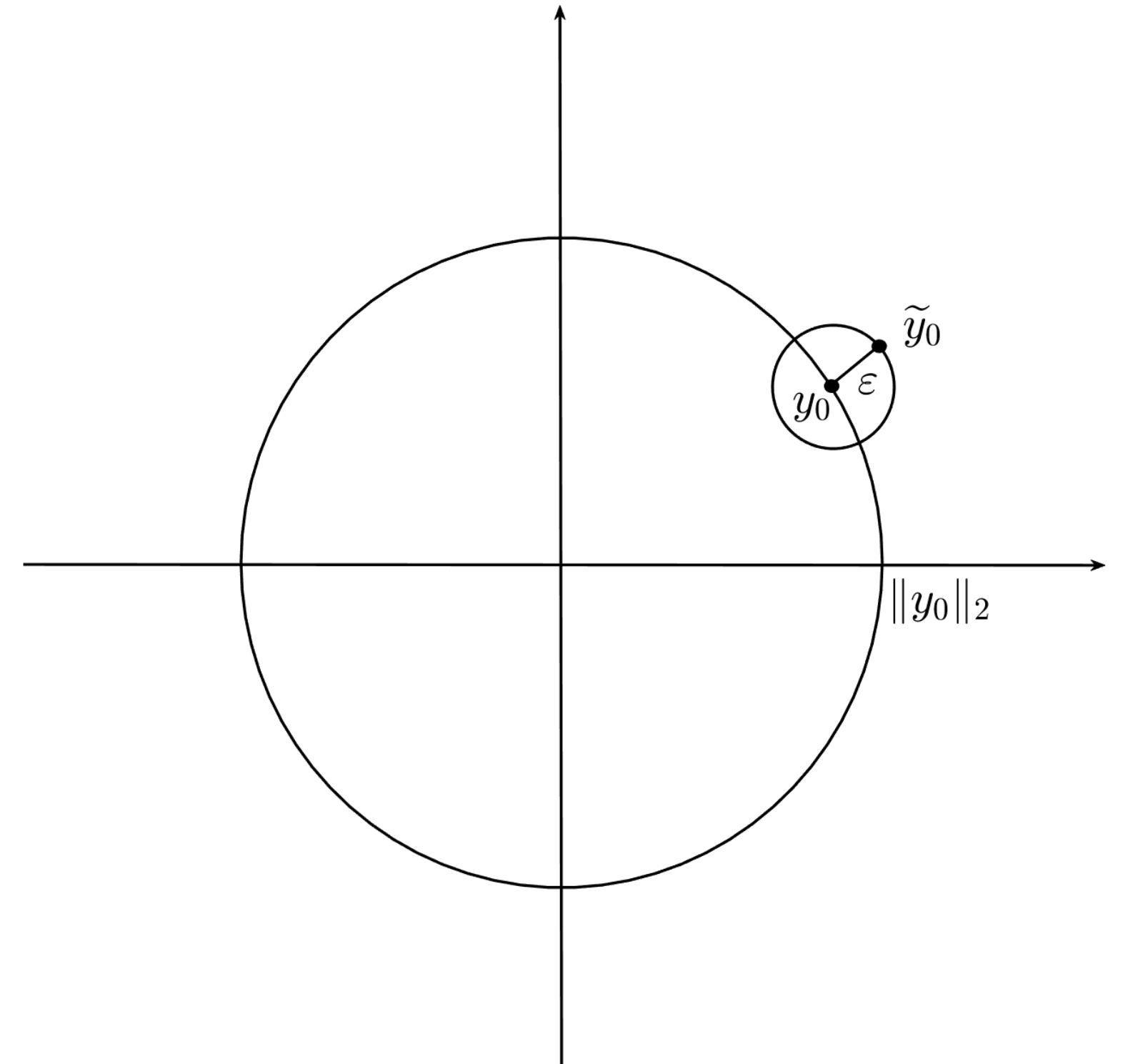}
     \end{subfigure}
\hfill  
\quad \begin{subfigure}[b]{1\textwidth}
         \centering
         \includegraphics[width=0.4\textwidth]{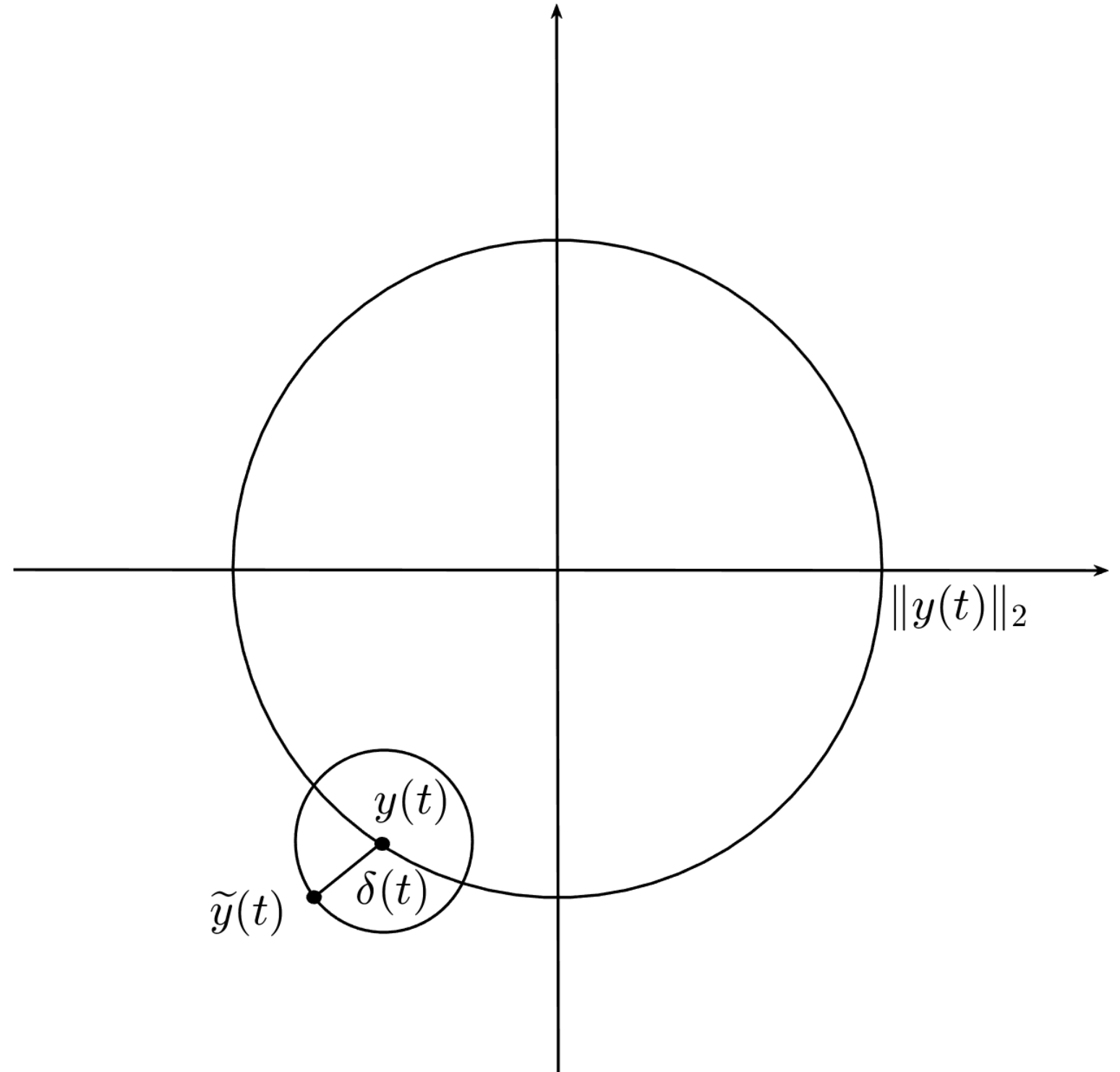}
     \end{subfigure}
\hfill
\caption{The relative error $\delta(t)$ when $\Vert y(t)\Vert_2$ is the scale unit.}
\label{scala}
\end{figure}

Assume that the initial value $y_{0}\neq 0$ of (\ref{ode}) is perturbed to $
\widetilde{y}_{0}$ and then the solution $y$ is perturbed to $\widetilde{y}$
. Fixed  an arbitrary
vector norm $\left\Vert \ \cdot \ \right\Vert $,  
we are interested in relating the normwise relative error 
\begin{equation}
\varepsilon :=\frac{\left\Vert \widetilde{y}_{0}-y_{0}\right\Vert }{
\left\Vert y_{0}\right\Vert }  \label{varepsilon}
\end{equation}
of the perturbed initial value $\widetilde{y}_0$ to the normwise relative error 
\begin{equation}
\delta \left( t\right) :=\frac{\left\Vert \widetilde{y}\left( t\right)
-y\left( t\right) \right\Vert }{\left\Vert y\left( t\right) \right\Vert }  \label{deltat}
\end{equation}
of the perturbed solution $\widetilde{y}(t)$. Observe that $\delta(0)=\varepsilon$. In the case of the Euclidean norm as vector norm, $\delta(t)$ is the distance between $\widetilde{y}(t)$ and $y(t)$ when $\Vert y(t)\Vert_2$ is used as scale unit (see Figure \ref{scala}).

\begin{remark}\label{averaged}
	Observe that the normwise relative error (\ref{deltat}) is invariant under multiplication of the vector norm by a constant. Therefore, these relative errors remain the same for both the \( p \)-norm and the mean \( p \)-norm,
	$$
	\Vert x \Vert_{p,\text{mean}} = \left( \sum\limits_{i=1}^n \frac{1}{n} \vert x_i \vert^p \right)^{\frac{1}{p}}, \quad x \in \mathbb{C}^n.
	$$
	Using such mean norms avoids the expected growth of the $p$-norms with the dimension \( n \).
\end{remark}

By writing the perturbed initial value as 
\begin{equation*}
\widetilde{y}_0=y_0+\varepsilon \Vert y_0\Vert \widehat{z}_0,
\end{equation*}
where $\widehat{z}_0\in\mathbb{C}^n$ is a unit vector (i.e., $\Vert \widehat{z
}_0 \Vert=1$) indicating the \emph{direction of perturbation}, we obtain 
\begin{equation}
\delta \left( t\right) =K\left( t,y_{0},\widehat{z}_{0}\right) \cdot  \label{magnification}
\varepsilon ,
\end{equation}
where 
\begin{equation*}
K\left( t,y_{0},\widehat{z}_{0}\right) :=\frac{\left\Vert \mathrm{e}^{tA}
\widehat{z}_{0}\right\Vert }{\left\Vert \mathrm{e}^{tA}\widehat{y}
_{0}\right\Vert } 
\end{equation*}
with $\widehat{y}_{0}:=\frac{y_{0}}{\left\Vert y_{0}\right\Vert }$ the \emph{normalized initial value}. We define $K(t,y_0,\widehat{z}_0)$ as the \emph{
directional pointwise condition number} of the problem (\ref{due}).

We also introduce 
\begin{equation*}
K\left( t,y_{0}\right) :=\max\limits_{\substack{ \widehat{z}_0\in \mathbb{C
}^{n}  \\ \Vert \widehat{z}_0\Vert =1}}K\left( t,y_{0},\widehat{z}
_0\right) =\frac{\left\Vert \mathrm{e}^{tA}\right\Vert }{\left\Vert \mathrm{e
}^{tA}\widehat{y}_{0}\right\Vert }, 
\end{equation*}
where $\left\Vert \mathrm{e}^{tA}\right\Vert $ is the matrix norm of $
\mathrm{e}^{tA}$ induced by the vector norm $\left\Vert \ \cdot \
\right\Vert $. We define $K(t,y_0) $ as the \emph{pointwise condition number} of the
problem (\ref{due}) (see \cite{Burgisser2013} for the definition of
condition number of a general problem). It is the worst condition number $
K(t,y_0,\widehat{z}_0)$ as the direction of perturbation $\widehat{
z}_0$ varies.

In general, we know nothing about the direction $\widehat{z}_{0}$ of the perturbation of $y_0$. Therefore, the condition number $K\left(t,y_{0}\right)$ is more useful than the condition number $K\left(t,y_{0},\widehat{z}_{0}\right)$. However, considering $K\left(t,y_{0},\widehat{z}_{0}\right)$ allows us to understand what happens for a general direction of perturbation, not necessarily the worst one. Moreover, in the discussion of the GDP-ND model in Section 6, the relevant condition number is $K\left(t,y_{0},\widehat{z}_{0}\right)$, since here we have information about the direction of perturbation.

\subsection{Normwise and componentwise relative errors} \label{nwcw}

We have introduced the \emph{normwise} relative errors (\ref{varepsilon}) and (\ref{deltat}),  but it may also be worthwhile to consider the \emph{componentwise} relative errors
\begin{equation*}
\varepsilon_l =\frac{\left\vert \widetilde{y}_{0l}-y_{0l}\right\vert }{
\left\vert y_{0l}\right\vert },\ l\in\{1,\ldots,n\}, \label{cre0}
\end{equation*}
of the perturbed initial value and 
\begin{equation*}
\delta_l \left( t\right) =\frac{\left\vert \widetilde{y}_l\left( t\right)
-y_l\left( t\right) \right\vert }{\left\vert y_l\left( t\right) \right\vert },\ l\in\{1,\ldots,n\}, \label{cre}
\end{equation*}
of the perturbed solution, where $y_{0l}$, $\widetilde{y}_{0l}$, $y_{l}(t)$ and $\widetilde{y}_{l}(t)$, $l\in\{1,\ldots,n\}$, are the components of $y_{0}$, $\widetilde{y}_{0}$, $y(t)$ and $\widetilde{y}(t)$, respectively. Componentwise relative errors were considered in \cite{FarooqMaset2021-2}, for $A$ diagonalizable. 

We can derive information regarding the componentwise relative errors from the normwise relative errors. In fact, if the vector norm $\left\Vert \ \cdot \ \right\Vert$ is a $p$-norm, then :
\begin{itemize}
\item [1)]
$
\varepsilon\leq \max\limits_{l\in\{1,\ldots,n\}}\varepsilon_l;
$
\item [2)]
$
 \delta(t)\leq \max\limits_{l\in\{1,\ldots,n\}}\delta_l \left( t\right);
$
\item [3)] 
$
\delta_l \left( t\right)\leq \frac{\Vert y(t)\Vert}{\vert y_l(t)\vert}\delta(t),\ l\in\{1,\ldots,n\}.
$
\end{itemize}

In particular, point 3) is useful for estimating the order of magnitude of $\delta_l(t)$, once the order of magnitude of $\delta(t)$ is known (recall the observation regarding relative errors of the two components $Q(t)$ and $B(t)$ of the GDP-ND model of Subsection \ref{twom}). In light of this, an important question is to understand how frequently components with a large ratio $\frac{\Vert y(t)\Vert}{\vert y_l(t)\vert}$ appear. To this aim, we consider the $\infty$-norm as vector norm. For other $p$-norms, observe that 
\begin{equation}
	\frac{\Vert y(t)\Vert}{\vert y_l(t)\vert}\leq n^{\frac{1}{p}}\frac{\Vert y(t)\Vert_\infty}{\vert y_l(t)\vert}. \label{n1/p}
\end{equation}

\begin{figure}[tbp] 
\includegraphics[width=1\textwidth]{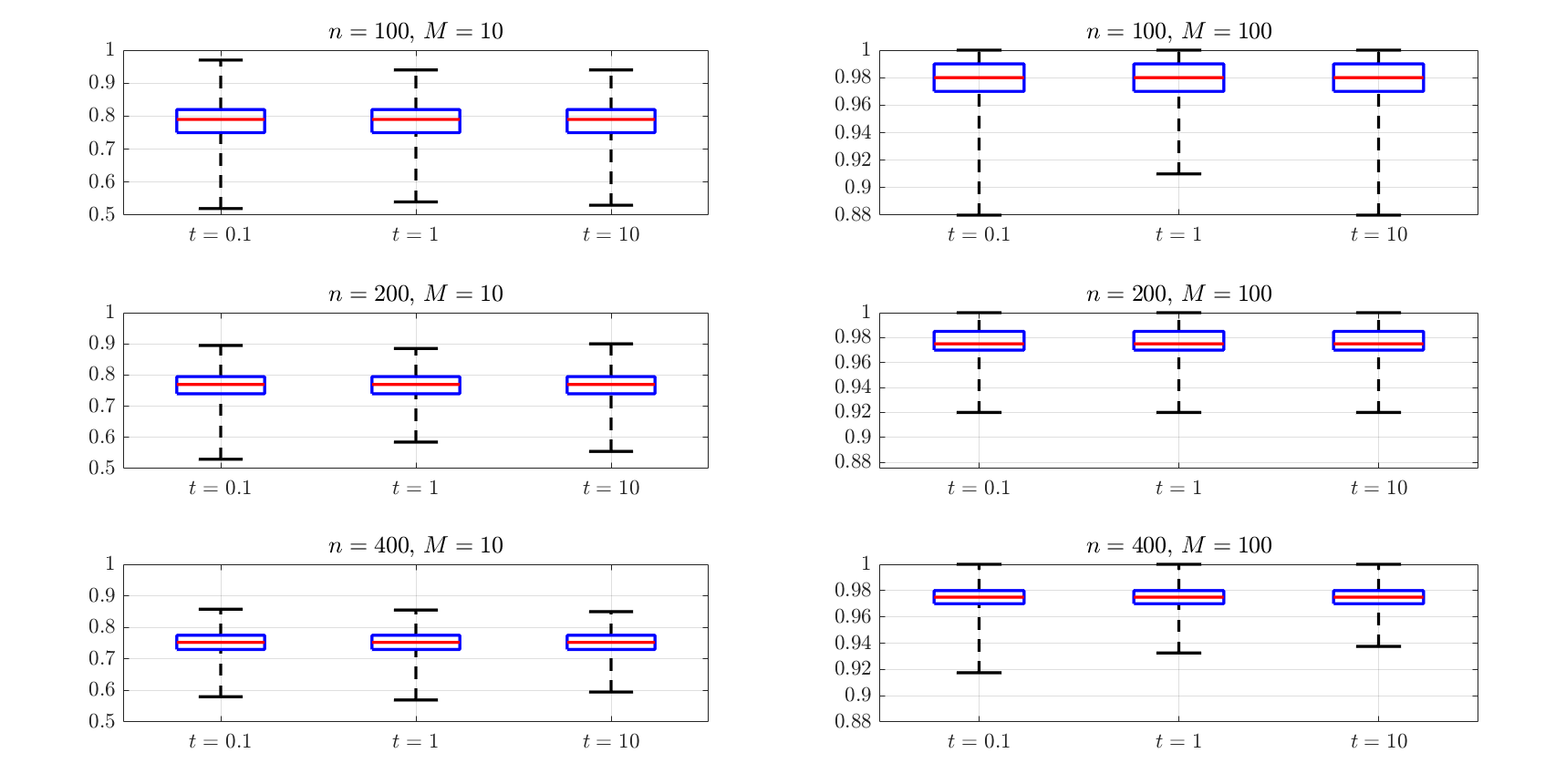}
\caption{Box plots of $r(t,M)$ in (\ref{fraction}) for $10\ 000$ random instances of (\ref{ode}) with $n=100,200,400$.}
\label{percentages} 
\end{figure}

In Figure \ref{percentages}, we see box plots of
\begin{equation}
r(t,M)=\frac{\text{number of components $y_l(t)$, $l\in\{1,\ldots,n\}$,  such that $\frac{\Vert y(t)\Vert_\infty}{\vert y_l(t)\vert}\leq M$}}{n}\label{fraction}
\end{equation}
for $10\ 000$ instances of (\ref{ode}), where $A$ of order $n=100,200,400$ and $y_0$ have entries sampled from the standard normal distribution. We consider $t=0.1,1,10$ and $M=10,100$.

We have evidence that the number of components with $\frac{\Vert y(t)\Vert_\infty}{\vert y_l(t)\vert}\leq M$ is a substantial (for $M=10$) or high (for $M=100$) percentage of the total number of components: in all instances the percentage is larger than $50\%$ for $M=10$ and $88\%$ for $M=100$; and, in three quarters of the instances the percentage is larger than $73\%$ for $M=10$ and $97\%$ for $M=100$.

Consequently, by point 3) above, we have evidence that \emph{a substantial percentage of the componentwise relative errors of the solution have an order of magnitude not larger than that of the normwise relative error, and a high percentage have an order of magnitude not larger than one plus that of the normwise relative error, when the $\infty$-norm is used as vector norm.}  For other $p$-norms, the factor $n^{\frac{1}{p}}$ in (\ref{n1/p}) must be taken into account.

In conclusion, analyzing the normwise relative error can provide valuable insights into the componentwise relative errors.

\section{Notations and definitions} \label{sectionAB}
We consider the spectrum of the matrix $A$ as partitioned into the sets $\Lambda_j$, $j\in\{1,\ldots,q\}$, where $\Lambda_j$ contains all the eigenvalues with the same real part $r_j$ and $r_1>r_2>\cdots>r_q$ holds. Observe that $\Lambda_1$ is the set of the rightmost eigenvalues of $A$.

Let $j\in\{1,\ldots,q\}$. We say that $\Lambda_j$ is \emph{simple real} if it consists of a real simple eigenvalue, and \emph{simple complex} if it consists of a single pair of complex conjugate simple eigenvalues. For $\Lambda_j$ simple real, we denote by $\lambda_j$ the real eigenvalue in $\Lambda_j$. For $\Lambda_j$ simple complex, we denote by $\lambda_j$ and $\overline{\lambda_j}$, with $\lambda_j$ having positive imaginary part $\omega_j$, the pair of complex conjugate simple eigenvalues in $\Lambda_j$.

\begin{remark} \label{remark1}
It is a generic case for the real matrix $A$ to have, for any $j\in\{1,\ldots,q\}$, $\Lambda_j$ simple real or simple complex.
\end{remark}

\subsection{$\boldsymbol{w^{(j)}}$, $\boldsymbol{v^{(j)}}$, $\boldsymbol{\widehat{\Theta }_j\left( t,u\right)}$, $\boldsymbol{\widehat{\Theta}_j\left(t\right)}$ and $\boldsymbol{f_j}$}

Let $j\in\{1,\ldots,q\}$.

For $\Lambda_j$ simple real, let $w^{(j)}$  (a real row vector)  and $v^{(j)}$ (a real column vector)  be left and right, respectively, 
eigenvectors corresponding to $\lambda_j$ such that $w^{(j)}v^{(j)}=1$ and let  $\widehat{w}^{(j)}=\frac{w^{(j)}}{\Vert w^{(j)} \Vert }$ and $\widehat{v}^{(j)}=\frac{v^{(j)}}{\Vert v^{(j)} \Vert }$ be their normalizations. Here,
	$$\Vert w^{\left( j\right)} \Vert=\max\limits_{\substack{ u\in \mathbb{R}^{n}  \\
			\Vert u\Vert=1}}\vert w^{\left( j\right) }u\vert
	$$
is the real induced norm of the real matrix $w^{(j)}$.

For $\Lambda_j$ simple complex, let $w^{(j)}$  (a complex row vector)  and $v^{(j)}$ (a complex column vector)  be left and right, respectively, eigenvectors corresponding to $\lambda_j$ such that $w^{(j)}v^{(j)}=1$ and let  $\widehat{w}^{(j)}=\frac{w^{(j)}}{\Vert w^{(j)} \Vert }$ and $\widehat{v}^{(j)}=\frac{v^{(j)}}{\Vert v^{(j)} \Vert }$ be their normalizations. Here, $$\Vert w^{\left( j\right) }\Vert=\max\limits_{\substack{ u\in \mathbb{C}^{n}  \\
			\Vert u\Vert=1}}\vert w^{\left( j\right) }u\vert$$ is the complex induced norm of the complex matrix $w^{(j)}$.
Moreover, given the polar forms
	\begin{eqnarray*}
		&&\widehat{v}_k^{\left( j\right) }=\left\vert \widehat{v}_{k}^{\left(
			j\right) }\right\vert \mathrm{e}^{\mathrm{i}\alpha _{jk}},\ k\in\{1,\ldots ,n\},\\
		&&\widehat{w}_l^{\left( j\right) }=\left\vert \widehat{w}_{l}^{\left(
			j\right) }\right\vert \mathrm{e}^{\mathrm{i}\beta _{jl}},\ l\in\{1,\ldots ,n\},
			\end{eqnarray*}
	of the components of the complex vectors $\widehat{v}^{\left( j\right) }$ and $\widehat{w}^{\left( j\right) }$ and the polar form
	\begin{eqnarray*}
		&&\widehat{w}^{\left( j\right) }u=\left\vert\widehat{w}^{\left( j\right)
		}u\right\vert\mathrm{e}^{\mathrm{i}\gamma _j\left( u\right) }
	\end{eqnarray*}
	of the complex scalar $\widehat{w}^{\left( j\right) }u$, where $u\in\mathbb{R}^n$, we introduce the vector
	\begin{equation*}
	\widehat{\Theta }_j\left( t,u\right) :=\left( \left\vert \widehat{v}
	_{k}^{\left( j\right) }\right\vert \cos \left( \omega _{j}t+\alpha
	_{jk}+\gamma_j\left( u\right) \right) \right) _{k=1,\ldots
		,n}\in \mathbb{R}^{n}
\end{equation*}
and the matrix
\begin{equation*}
	\widehat{\Theta }_j\left( t\right) =\left[ \left\vert \widehat{v}
	_{k}^{\left( j\right) }\right\vert \left\vert \widehat{w}_{l}^{\left(
		j\right) }\right\vert \cos \left( \omega _{j}t+\alpha _{jk}+\beta _{jl}\right) \right] _{k,l=1,\ldots ,n}\in \mathbb{
		R}^{n\times n}.
\end{equation*}

Finally, for both cases $\Lambda_j$ simple real and $\Lambda_j$ simple complex, we introduce the quantity
\begin{equation}
	f_j:=\Vert w^{(j)} \Vert\cdot \Vert v^{(j)} \Vert\in [1,+\infty).  \label{fj}
\end{equation}

\subsection{Asymptotic forms and approximation with a given precision} \label{withprecision}

For the description of the asymptotic behavior of the condition numbers, we use the following notion of asymptotic form. Let $a(t)$ and $b(t)$ be real functions of $
t\in \mathbb{R}$. We say that $b(t)$ is an \emph{asymptotic form} of $a(t)$ and write
\begin{equation*}
	a\left( t\right) \sim b\left( t\right),\ t\rightarrow +\infty,  \label{sim}
\end{equation*}
if
\begin{equation*}
	\lim\limits_{t\rightarrow +\infty }\frac{a(t)}{b(t)}=1.
\end{equation*}
In other words, $b(t)$ is an asymptotic form of $a(t)$ if
\begin{equation*}
	\lim\limits_{t\rightarrow +\infty }\chi(t)=0,
\end{equation*}
where $\chi(t)$ is the relative error of $a(t)$ with respect to $b(t)$.

The closeness of $a(t)$ to its asymptotic form $b(t)$ at a finite time $t$ is measured by using the following notion of \emph{approximation with a given precision}. For $\epsilon\geq 0$, we write
$$
a(t)\approx b(t)\text{\ with precision\ }\epsilon
$$
if $\vert \chi(t)\vert \leq\epsilon$.

\section{The asymptotic condition numbers}

Next two results (see \cite{M2}) describe the asymptotic forms of the condition numbers $K(t,y_0,\widehat{z}_0)$ and $K(t,y_0)$ in the generic case for $A$ where the set $\Lambda_1$  is simple real or simple complex (recall Remark \ref{remark1}), and the generic case for $y_0$ and $\widehat{z}_0$ given by $w^{(1)}y_0\neq 0$ and  $w^{(1)}\widehat{z}_0\neq 0$.
\begin{theorem}\label{Threal}
Assume $\Lambda_1$ simple real. For $y_0$ and $\widehat{z}_0$ such that $w^{(1)}y_0\neq 0$ and  $w^{(1)}\widehat{z}_0\neq 0$, we have
	\begin{equation*}
			K\left(t,y_{0},\widehat{z}_{0}\right)\sim  K_\infty\left(t,y_{0},\widehat{z}_{0}\right)=K_\infty\left(y_{0},\widehat{z}_{0}\right):=\frac{\left\vert
			\widehat{w}^{(1)}\widehat{z}_{0}\right\vert }{\left\vert \widehat{w}^{(1)}\widehat{y}
			_{0}\right\vert },\ t\rightarrow +\infty,
\end{equation*}
and
	\begin{equation*}
K(t,y_0)\sim  K_\infty\left(t,y_{0}\right)=K_\infty\left(y_{0}\right):=\frac{1}{\left\vert \widehat{w}^{(1)}\widehat{y}
		_{0}\right\vert  },\ t\rightarrow +\infty.
\end{equation*}
\end{theorem}

\begin{theorem} \label{Thcomplex}
	Assume $\Lambda_1$ simple complex. For $y_0$ and $\widehat{z}_0$ such that $w^{(1)}y_0\neq 0$ and  $w^{(1)}\widehat{z}_0\neq 0$, we have
	\begin{eqnarray*}
		K\left( t,y_{0},\widehat{z}_{0}\right)\sim K_\infty\left( t,y_{0},\widehat{z}_{0}\right) = \mathrm{OSF}(y_0,\widehat{z}_0)\cdot\mathrm{OT}(t,y_0,\widehat{z}_0),\ t\rightarrow +\infty,
	\end{eqnarray*}
	and
	\begin{eqnarray*}
	K(t,y_0)\sim K_\infty\left( t,y_{0}\right) = \mathrm{OSF}(y_0)\cdot\mathrm{OT}(t,y_0),\ t\rightarrow +\infty,
	\end{eqnarray*}
	where
	\begin{equation*}
	\mathrm{OSF}(y_0,\widehat{z}_0):=\frac{\left\vert
		\widehat{w}^{(1)}\widehat{z}_{0}\right\vert }{\left\vert \widehat{w}^{(1)}\widehat{y}
		_{0}\right\vert }\text{\ \ and\ \ }\mathrm{OSF}(y_0):=\frac{1}{\left\vert \widehat{w}^{(1)}\widehat{y}
		_{0}\right\vert }  
	\end{equation*}
	and
	\begin{equation}
	\mathrm{OT}(t,y_0,\widehat{z}_0):= \frac{\left\Vert \widehat{\Theta}_1
		\left( t,\widehat{z}_{0}\right) \right\Vert }{\left\Vert \widehat{\Theta}_1
		\left( t,\widehat{y}_{0}\right) \right\Vert }\text{\ \ and\ \ }\mathrm{OT}(t,y_0):=\frac{\left\Vert \widehat{\Theta}_1
		\left( t\right) \right\Vert }{\left\Vert \widehat{\Theta}_1
		\left( t,\widehat{y}_{0}\right) \right\Vert }.  \label{OTOT}
	\end{equation}
\end{theorem}
The constants $\mathrm{OSF}(y_0,\widehat{z}_0)$ and $\mathrm{OSF}(y_0)$ are called \emph{oscillation scale factors} and $\mathrm{OT}(t,y_0,\widehat{z}_0)$ and $\mathrm{OT}(t,y_0)$, which are periodic functions of $t$ of period $\frac{\pi}{\omega_1}$, are called \emph{oscillating terms}. Note that the oscillation scale factors depend on the moduli, and the oscillating terms on the angles, in the polar forms of the complex numbers $\widehat{w}^{(1)}\widehat{y}_{0}$ and $\widehat{w}^{(1)}\widehat{z}_{0}$.

	\subsection{The case $\Lambda_1$ simple complex and the Euclidean norm}\label{Euclidean}
	In this subsection, in the case of $\Lambda_1$ simple complex, we analyze the oscillating terms when the vector norm is the Euclidean norm.

	We introduce
	$$
	V_1:=\left\vert \left(\widehat{v}^{(1)}\right)^T\widehat{v}^{(1)}\right\vert
	\text{\ \ and\ \ }W_1:=\left\vert \widehat{w}^{(1)}\left(\widehat{w}^{(1)}\right)^T\right\vert,
	$$
	moduli of the complex numbers $\left(\widehat{v}^{(1)}\right)^T\widehat{v}^{(1)}$ and $\widehat{w}^{(1)}\left(\widehat{w}^{(1)}\right)^T$. Here $^T$, unlike $^H$, denotes pure transposition without conjugation.
	We have $V_1,W_1\in[0,1)$.
	
\begin{figure}[tbp]
				\centering
				\par
				\includegraphics[width=1.1\textwidth]{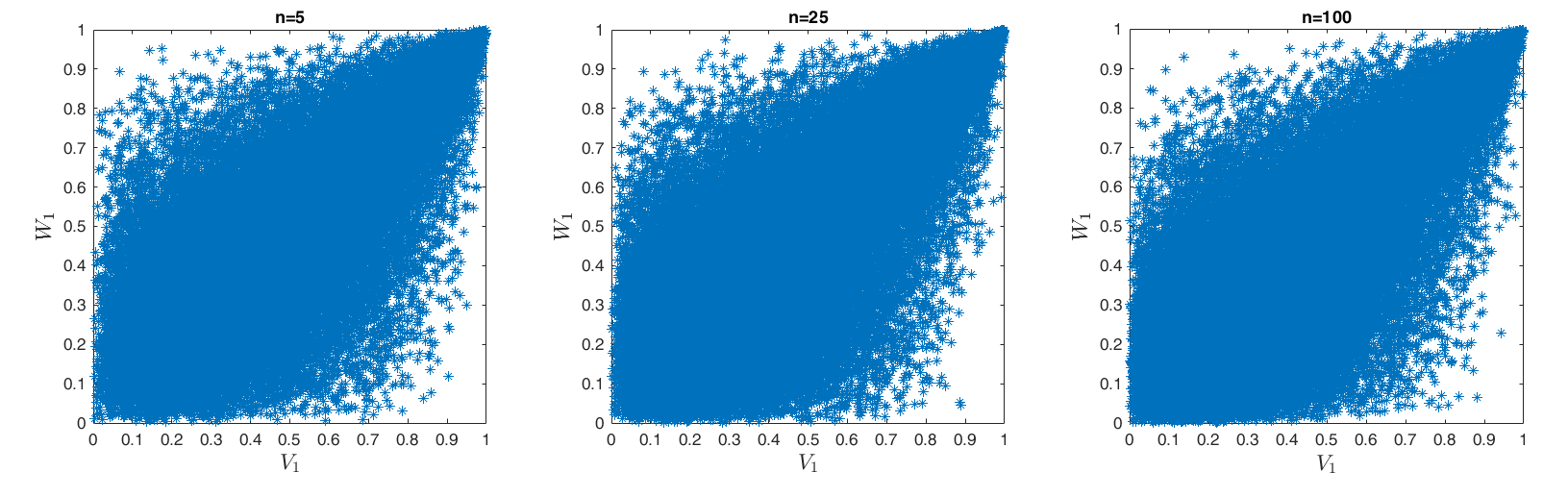}
				\caption{Scatter diagrams  $(V_1,W_1)$ for $50\ 000$ random matrices $A$ of order $n=5,25,100$.}
				\label{FigureVW}
	\end{figure}

\begin{figure}[tbp]
		\centering
		\par
		\includegraphics[width=1.1\textwidth]{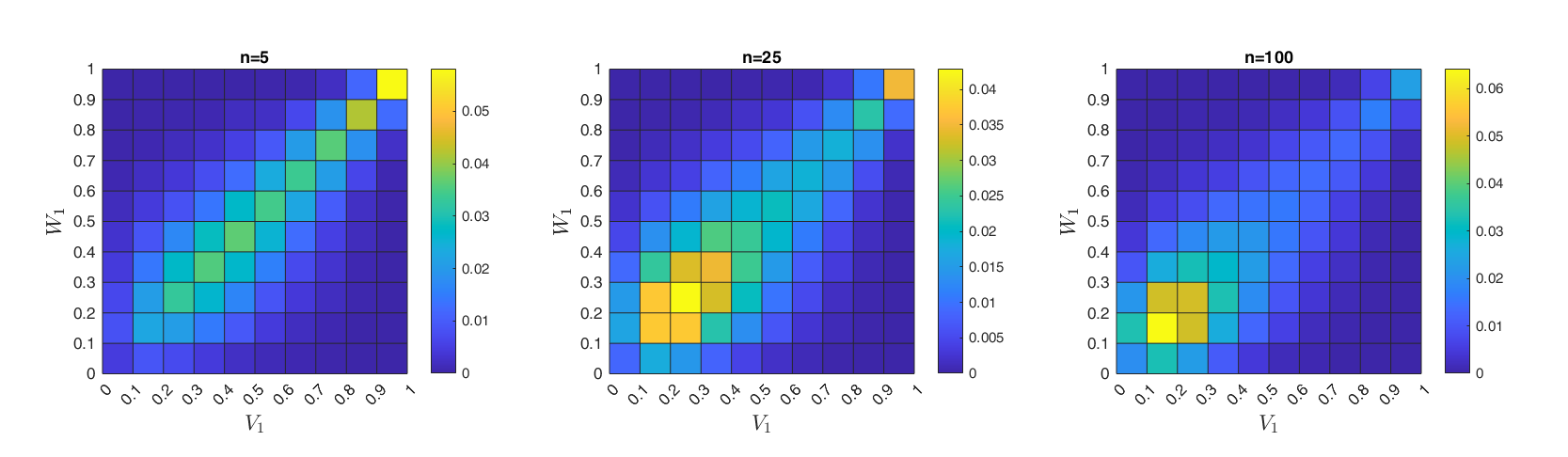}
		\caption{Bidimensional histograms with heatmaps of $(V_1,W_1)$ for the $50\ 000$ random matrices $A$ of order $n=5,25,100$ in Figure \ref{FigureVW}. Relative frequencies are reported.}
		\label{FigureVW1}
	\end{figure}

	In Figure \ref{FigureVW}, we see scatter diagrams of $(V_1,W_1)$ for $50\ 000$ random matrices $A$ of order $n=5,25,100$ such that $\Lambda_1$ is simple complex. The elements of $A$ are sampled from the standard normal distribution.  Observe that the distribution of the pairs in the square $[0,1)^2$ is not uniform, since the pairs tend to accumulate around the diagonal $V_1=W_1$. There are no instances close to the corners $(0,1)$ and $(1,0)$. In Figure \ref{FigureVW1}, we see bidimensional histograms with heatmaps of $(V_1,W_1)$. The distribution of $(V_1,W_1)$ appears symmetric.

	The next theorem (see \cite{M2}) specifies the ranges within which the oscillating terms $\mathrm{OT}\left(t,y_0,\widehat{z}_0\right)$ and $\mathrm{OT}\left(t,y_0\right)$ in (\ref{OTOT}) can vary. 
	\begin{theorem} \label{theorem9z}
		Suppose that the vector norm is the Euclidean norm. For $y_0$ and $\widehat{z}_0$ such that $w^{(1)}y_0\neq 0$ and  $w^{(1)}\widehat{z}_0\neq 0$, we have
		\begin{equation*}
		\min\limits_{\substack{y_0,\widehat{z}_0\in\mathbb{R}^n\\ \left\Vert \widehat{z}_0\right\Vert_2=1\\ w^{(1)}y_0\neq 0,w^{(1)}\widehat{z}_0\neq 0}}\min\limits_{t\in\mathbb{R}} \mathrm{OT}\left(t,y_0,\widehat{z}_0\right)=\sqrt{\frac{1-V_1}{1+V_1}}
		\end{equation*}
		and
		\begin{equation*}
			\max\limits_{\substack{y_0,\widehat{z}_0\in\mathbb{R}^n\\ \left\Vert \widehat{z}_0\right\Vert_2=1\\ w^{(1)}y_0\neq 0,w^{(1)}\widehat{z}_0\neq 0}}\max\limits_{t\in\mathbb{R}} \mathrm{OT}\left(t,y_0,\widehat{z}_0\right)=\sqrt{\frac{1+V_1}{1-V_1}}
		\end{equation*}
	Moreover, we have
	\begin{equation*}
		\min\limits_{\substack{y_0\in\mathbb{R}^n\\ w^{(1)}y_0\neq 0}}\min\limits_{t\in\mathbb{R}}\mathrm{OT}\left(t,y_0\right)=a^{\min}_{V_1W_1}\text{\ \ and\ \ }
		\max\limits_{\substack{y_0\in\mathbb{R}^n\\ w^{(1)}y_0\neq 0}}\max\limits_{t\in\mathbb{R}}\mathrm{OT}\left(t,y_0\right)=a^{\max}_{V_1W_1},
	\end{equation*}
	where
	\begin{equation*}
		a^{\min}_{V_1,W_1}=\left\{
		\begin{array}{l}
			\sqrt{\frac{(1-V_1)(1+W_1)}{2(1+V_1)}}\text{\ if\ }V_1\leq W_1\\
			\\
			\sqrt{\frac{1-W_1}{2}}\text{\ if\ }V_1\geq W_1
		\end{array}
		\right.\text{\ \ and\ \ }
		a^{\max}_{V_1W_1}=\sqrt{\frac{(1+V_1)(1+W_1)}{2(1-V_1)}} 
	\end{equation*}
	\end{theorem}

	In light of Theorems \ref{Threal}, \ref{Thcomplex} and \ref{theorem9z} and Remark \ref{remark1}, we can state a fact B) regarding the relative error of the perturbed solution, analogous to fact A) concerning the absolute error of the perturbed solution and presented at the beginning of the paper.
	\begin{itemize}
		\item [B)] \emph{For a generic matrix $A$, for a generic initial value $y_0$ and for a generic perturbation of $y_0$, the relative error of the perturbed solution asymptotically (as $t\rightarrow +\infty$) neither diverges nor decays to zero. Instead, it converges to a non-zero constant, in case of a real rightmost eigenvalue, and to a periodic oscillating function bounded and uniformly away from zero, in case of a complex conjugate pair of rightmost eigenvalues.}
	\end{itemize}
	
Indeed, Fact B) holds for any matrix $A$, not only for matrices with $\lambda_1$ simple real or simple complex (see \cite{M1}).
	
\subsubsection{$V_1$ not close to $1$} \label{V1c1}

By Theorem \ref{theorem9z} and the fact that $V_1,W_1\in [0,1)$, we have
\begin{equation}
	\frac{1}{k(V_1)}\leq \mathrm{OT}(t,y_0,\widehat{z}_0)\leq k(V_1)\text{\ \ and\ \ }\frac{1}{k(V_1)}\leq \mathrm{OT}(t,y_0)\leq k(V_1),\ t\in\mathbb{R}, \label{kV1}
\end{equation}
where 
$$
k(V_1)=\sqrt{\frac{2}{1-V_1}}.
$$
Therefore, we can conclude that if $V_1$ is not close to 1, then, for any $y_0$ and $\widehat{z}_0$, $\mathrm{OT}(t, y_0, \widehat{z}_0)$ and $\mathrm{OT}(t, y_0)$ do not assume large or small values as $t$ varies, i.e., the values of $\mathrm{OT}(t, y_0, \widehat{z}_0)$ and $\mathrm{OT}(t, y_0)$ as $t$ varies have the order of magnitude $1$. Consequently, \emph{the values of $K_\infty(t, y_0, \widehat{z}_0)$ and $K_\infty(t, y_0)$ as $t$ varies have the same order of magnitude as $\mathrm{OSF}(y_0, \widehat{z}_0)$ and $\mathrm{OSF}(y_0)$}, respectively.

\subsubsection{Is $V_1$ close to $1$?}\label{V1close1}

In light of the previous Subsection \ref{V1c1},  it is of interest to understand how frequently $V_1$ is close to 1. 
	
The following table shows, for the $50\ 000$ random matrices of Figures \ref{FigureVW} and \ref{FigureVW1}, the percentages of cases with $V_1$ greater than $0.9$, $0.99$, $0.999$ and $0.9999$.
\begin{equation}
	\begin{tabular}{|l|l|l|l|l|}
		\hline
		& $V_1>0.9$ & $V_1>0.99$ & $V_1>0.999$ & $V_1>0.9999$ \\
		\hline
		$n=5$ &  7\%   & 0.7\% &  0.06\% &  0.004\% \\
		\hline
		$n=25$ & 5\%   & 0.5\% &  0.05\% &  0.004\%  \\
		\hline
		$n=100$  &  3\%   & 0.3\% & 0.02\% & 0\%   \\
		\hline	  
	\end{tabular}
	\label{table}
\end{equation}

We have strong evidence that $V_1$ is rarely close to 1.

Regarding the possible order of magnitude of the oscillating terms when $V_1$ is close to $1$, note that the values of $k(V_1)$ (recall (\ref{kV1})) for $V_1 = 0.9, 0.99, 0.999, 0.9999$ in table (\ref{table}) are
$$
\begin{tabular}{|l|l|l|l|l|}
	\hline
	& $V_1=0.9$ & $V_1=0.99$ & $V_1=0.999$ & $V_1=0.9999$ \\
	\hline
	$k(V_1)$ &  $4.4721$   & $14.1421$ &  $44.7214$ &  $141.4214$ \\
	\hline  
\end{tabular}
$$

\section{Two important issues}

This section discusses two important issues concerning the most important asymptotic condition number $K_\infty(t,y_{0})$: the \emph{asymptotic well-conditioning} of problem (\ref{due}) and the \emph{onset of the asymptotic behaviour}. Moreover, the effect of the non-normality of $A$ on these two important issues is considered.

\subsection{The asymptotic well-conditioning}
It is of interest to know for which initial values $y_0$ the problem (\ref{due}) is \emph{asymptotically well-conditioned} at $y_0$, i.e., $K_\infty(t,y_0)$ does not assume large values as $t$ varies.

We assume $\Lambda_1$ simple real or simple complex, a generic case for the matrix $A$ (recall Remark \ref{remark1}). Moreover, we suppose $w^{(1)}y_0\neq 0$, a generic case for the initial value $y_0$.

\subsubsection{The case $\Lambda_1$ simple real} 

For $\Lambda_1$ simple real and $w^{(1)}y_0\neq 0$, we have
\begin{equation*}
	K_\infty(y_{0})=\frac{1}{\left\vert \widehat{w}^{(1)}\widehat{y}
		_{0}\right\vert}
\end{equation*}
(see Theorem \ref{Threal}). Thus, we have the following fact (we use the term fact, rather than theorem, because the asymptotic well-conditioning has not a rigorous definition, since it is based on a vague term like "large values").
\begin{fact}\label{fact1}
	Assume $\Lambda_1$ simple real. Let $y_0$ be such that $w^{(1)}y_0\neq 0$. The problem (\ref{due}) is asymptotically well-conditioned at $y_0$ if and only if 
	$$
	\frac{1}{\left\vert \widehat{w}^{(1)}\widehat{y}_{0}\right\vert}=\frac{\Vert w^{(1)}\Vert\ \Vert y_0 \Vert}{\vert w^{(1)}y_0\vert}
	$$
	is not large.
\end{fact}

\subsubsection{The case $\Lambda_1$ simple complex} \label{lambdacomplex}

Suppose that $\Lambda_1$ is simple complex and the vector norm is the Euclidean norm.

We have
\begin{equation*}
	K_\infty(t,y_{0})=\mathrm{OSF}(y_0)\cdot \mathrm{OT}(t,y_0),
\end{equation*}
where the oscillation scale factor $\mathrm{OSF}(y_0)$ is given by 
$$
\mathrm{OSF}(y_0)= \frac{1}{\left\vert \widehat{w}^{(1)}\widehat{y}
	_{0}\right\vert}.
$$
If $V_1$ is not close to $1$, then the values of $K_\infty(t,y_0)$ as $t$ varies have the same order of magnitude as $\mathrm{OSF}(y_0)$, since the values of $\mathrm{OT}(t,y_0)$ as $t$ varies have the order of magnitude $1$ (see Subsection \ref{Euclidean}). Thus, we have the following fact
\begin{fact}
	Assume $\Lambda_1$ simple complex and the Euclidean norm as vector norm. Let $y_0$ be such that $w^{(1)}y_0\neq 0$. If $V_1$ is not close to $1$, then the problem (\ref{due}) is asymptotically well-conditioned at $y_0$ if and only
	$$
	\mathrm{OSF}(y_0)=\frac{1}{\left\vert \widehat{w}^{(1)}\widehat{y}_{0}\right\vert}=\frac{\Vert w^{(1)}\Vert\ \Vert y_0 \Vert}{\vert w^{(1)}y_0\vert}
	$$
	is not large.
\end{fact}

If $V_1$ is close to $1$, then $\mathrm{OSF}(y_0)$ does not determine the asymptotic well-conditioning of the problem (\ref{due}), unlike when $V_1$ is not close to $1$, since it is no longer guaranteed that the values of $\mathrm{OT}(t, y_0)$ as $t$ varies maintain the order of magnitude of $1$.

Indeed, if $V_1$ is
close to $1$, $W_1$ is not close to $1$ and $y_0$ is the second right singular vector of the matrix
$$
R_1=\left[
\begin{array}{c}
	\mathrm{Re}\left(\widehat{w}^{(1)}\right)\\
	\mathrm{Im}\left(\widehat{w}^{(1)}\right)
\end{array}
\right],
$$
then $\mathrm{OSF}(y_0)$ is not large, but $K_\infty(t,y_0)$ assumes large values as $t$ varies, i.e., the problem (\ref{due}) is not asymptotically well-conditioned. In this situation, we have an oscillating term $\mathrm{OT}(t,y_0)$ assuming large values as $t$ varies. See \cite{M2}.

However, by recalling Subsection \ref{V1close1}, observe that $V_1$ is rarely close to $1$.

\subsection{The onset of the asymptotic behavior} It is of interest to know the  \emph{onset of the asymptotic behavior}, i.e., when $K(t,y_0)$ begins to be close to $
	K_\infty(t,y_{0})$.

We assume $\Lambda_j$, $j\in\{1,\ldots,q\}$, simple real or simple complex, a generic case for $A$ (recall Remark \ref{remark1}). Moreover, we suppose $w^{(1)}y_0\neq 0$, a generic case for $y_0$.

We recall from \cite{M2} the following result.  See Subsection \ref{withprecision} for the notation $\approx$ with precision.
\begin{theorem}\label{onset}
	Assume that, for any $j\in\{1,\ldots,q\}$, $\Lambda_j$ is simple single real or simple  single complex. For $y_0$ such that $w^{(1)}y_0\neq 0$, we have
	\begin{equation}
		K\left(t,y_{0}\right) \approx K_\infty\left(t,y_{0}\right)\text{\ with
					precision\ }\frac{\epsilon(t)+\epsilon(t,\widehat{y}_0)}{
					1-\epsilon(t,\widehat{y}_0)}  \label{appK}
		\end{equation}
whenever $\epsilon(t,\widehat{y}_0)<1$, where
	\begin{equation}
			\epsilon(t,\widehat{y}_0):=\sum\limits_{j=2}^{q}\mathrm{e}^{\left( r_{j}-r_{1
						}\right) t}\frac{f_j}{f_1}\cdot\frac{\left\vert \widehat{w}^{(j)}\widehat{y}_0\right\vert
				}{\left\vert \widehat{w}^{(1)}\widehat{y}_0\right\vert}G_j\left(t,\widehat{y}_0\right) \label{appK1}
		\end{equation}
and
	\begin{equation}
			\epsilon(t):=\sum\limits_{j=2}^{q}\mathrm{e}^{\left( r_{j}-r_{1
						}\right) t}\frac{f_j}{f_1}G_j(t), \label{appK2}
		\end{equation}
with $G_j(t,\widehat{y}_0)$ and $G_j(t)$, $j\in\{2,\ldots,q\}$, given as follows.
\begin{itemize}
	\item [1)]  	If both $\Lambda_j$ and $\Lambda_1$ are simple real, then
	\begin{equation*}
			G_j(t,u)=1\text{\ \ and\ \ }G_j(t)=1.
		\end{equation*}
	\item [2)]  If $\Lambda_j$ is simple complex and $\Lambda_1$ is simple real, then
	\begin{equation*}
			G_j(t,u)=2\left\Vert \widehat{\Theta}
			_j(t,u)\right\Vert\text{\ \ and\ \ }G_j(t)=2\left\Vert \widehat{\Theta}_j(t)\right\Vert.
		\end{equation*}
\item [3)]  If $\Lambda_j$ is simple real and $
			\Lambda_1$ is simple complex, then
			\begin{equation*}
						G_j(t,u)=\frac{1}{2\left\Vert 
									\widehat{\Theta}_1(t,u)\right\Vert}\text{\ \ and\ \ }G_j(t)=\frac{1}{2\left\Vert \widehat{\Theta}_1(t)\right\Vert}.
					\end{equation*}
			\item [4)]  If both $\Lambda_j$ and $\Lambda_1$ are simple complex, then
			\begin{equation*}
						G_j(t,u)=\frac{\left\Vert \widehat{
												\Theta}_j(t,u)\right\Vert}{\left\Vert \widehat{\Theta}_1(t,u)\right\Vert}
						\text{\ \ and\ \ }G_j(t)=\frac{
									\left\Vert \widehat{\Theta}_j(t)\right\Vert}{\left\Vert \widehat{\Theta}
									_1(t)\right\Vert}.
					\end{equation*}
\end{itemize}
\end{theorem}

In the following, $\widehat{t}$ is the time used as time unit. If $r_1\neq 0$, one possible choice for the time unit is the \emph{characteristic time} $\frac{1}{\vert r_1\vert}$.

\subsubsection{The case $\Lambda_1$ simple real}

The following result for the case $\Lambda_1$ simple real follows by Theorem \ref{onset}. 

\begin{theorem} \label{onset1}
Assume $\Lambda_1$ simple real and, for $j\in\{2,\ldots,q\}$, $\Lambda_j$ simple real or simple complex. Assume that a $p$-norm is used as vector norm. Let $y_0$ be such that $w^{(1)}y_0\neq 0$. For any $\epsilon>0$, we have
$$
K(t,y_0)\approx K_\infty(y_0)\text{\ with precision $\epsilon$}
$$
if 
\begin{eqnarray}
	&&\frac{t}{\widehat{t}}\geq \max\limits_{j\in\{2,\ldots,q\}}\frac{1}{\left(r_{1}-r_{j}\right) \widehat{t}}\left(\log\frac{2+\epsilon}{\epsilon}+\log(q-1)+\log\frac{f_j}{f_1}\right. \notag\\
	&&\quad\quad\quad\quad\quad\quad\quad\quad\quad\quad\quad  \left.
	+\max\left\{0,\log\frac{\left\vert \widehat{w}^{(j)}\widehat{y}
		_0\right\vert}{\left\vert \widehat{w}^{(1)}\widehat{y}
		_0\right\vert}\right\}+\log 2\right). \notag \\
	\label{realepsilon<}
\end{eqnarray}
\end{theorem}

\begin{proof}
Since the vector norm is a $p$-norm, we have (see \cite{M2})
	\begin{equation*}
			G_j(t,u)\leq 2\text{\ \ and\ \ }G_j(t)\leq 2
		\end{equation*}
in (\ref{appK1}) and (\ref{appK2}).  
As a consequence, we obtain
\begin{eqnarray*}
	\max\{\epsilon(t,\widehat{y}_0),\epsilon(t)\}&\leq& \sum\limits_{j=2}^{q}\mathrm{e}^{\left(r_{j}-r_{1 }\right) t}\frac{f_j}{f_1} \max\left\{1,
	\frac{\left\vert \widehat{w}^{(j)}\widehat{y}_0\right\vert}{\left\vert \widehat{w}^{(1)}\widehat{y}_0\right\vert}\right\}2\\
	&\leq&(q-1)\max\limits_{j\in\{2,\ldots,q\}}\mathrm{e}^{\left(r_{j}-r_{1 }\right) t}\frac{f_j}{f_1} \max\left\{1,
	\frac{\left\vert \widehat{w}^{(j)}\widehat{y}_0\right\vert}{\left\vert \widehat{w}^{(1)}\widehat{y}_0\right\vert}\right\}2.
	\label{sumssums}
\end{eqnarray*}
Therefore, for $\epsilon_0>0$, we have
$
 \max\{\epsilon(t,\widehat{y}_0),\epsilon(t)\}\leq \epsilon_0
$
if
\begin{equation*}
	\frac{t}{\widehat{t}}\geq \max\limits_{j\in\{2,\ldots,q\}}\frac{1}{\left(r_{1}-r_{j}\right) \widehat{t}}\left(\log\frac{q-1}{\epsilon_0}+\log\frac{f_j}{f_1}
	+\max\left\{0,\log\frac{\left\vert \widehat{w}^{(j)}\widehat{y}
		_0\right\vert}{\left\vert \widehat{w}^{(1)}\widehat{y}
		_0\right\vert}\right\}+\log 2\right).
\end{equation*}
Now, given $\epsilon>0$, consider $\epsilon_0$ such that
$$
\frac{2\epsilon_0}{1-\epsilon_0}=\epsilon,\text{\ i.e.\ },\ \epsilon_0=\frac{\epsilon}{2+\epsilon}, 
$$
and use (\ref{appK}).
\end{proof}
Therefore, the smaller  the right-hand side in (\ref{realepsilon<}), the earlier the onset of the asymptotic behavior.

\subsubsection{The case $\Lambda_1$ simple complex}

The following result for the case $\Lambda_1$ simple complex follows by Theorem \ref{onset}.

\begin{theorem} \label{onset2}
	Assume $\Lambda_1$ simple complex and, for $j\in\{2,\ldots,q\}$, $\Lambda_j$ simple real or simple complex. Assume that the Euclidean norm is used as vector norm. Let $y_0$ be such that $w^{(1)}y_0\neq 0$. For any $\epsilon>0$, we have
	$$
	K(t,y_0)\approx K_\infty(t,y_0)\text{\ with precision $\epsilon$}
	$$
	if 
	\begin{eqnarray}
		&&\frac{t}{\widehat{t}}\geq \max\limits_{j\in\{2,\ldots,q\}}\frac{1}{\left(r_{1}-r_{j}\right) \widehat{t}}\left(\log\frac{2+\epsilon}{\epsilon}+\log(q-1)+\log\frac{f_j}{f_1}\right. \notag\\
		&&\quad\quad\quad\quad\quad\quad\quad\quad\quad\quad\quad  \left.
		+\max\left\{0,\log\frac{\left\vert {\widehat w}^{(j)}\widehat{y}
			_0\right\vert}{\left\vert {\widehat w}^{(1)}\widehat{y}
			_0\right\vert}\right\}+\frac{1}{2}\log\frac{1}{1-V_1}+\log 2\right). \notag \\
		\label{complexepsilon<}
	\end{eqnarray}
\end{theorem}

\begin{proof}
Since the vector norm is the Euclidean norm, we have (see \cite{M2})
\begin{equation*}
	G_j(t,u)\leq \sqrt{\frac{1}{2(1-V_1)}}\text{\ \ and\ \ }G_j(t)\leq \sqrt{\frac{1}{a_{V_1W_1}}}
\end{equation*}
if $\Lambda_j$ is simple real, and
\begin{equation*}
	G_j(t,u)\leq \sqrt{\frac{1+V_j}{1-V_1}}\text{\ \ and\ \ }G_j(t)\leq \sqrt{\frac{(1+V_j)(1+W_j)}{a_{V_1W_1}}}
\end{equation*}
if $\Lambda_j$ is simple complex, where
$$
a_{V_1W_1}=\left\{
	\begin{array}{l}
				(1-V_1)(1+W_1)\text{\ if\ }V_1\leq W_1\\
				\\
				(1+V_1)(1-W_1)\text{\ if\ }V_1\geq W_1,
			\end{array}
	\right.
$$
$V_j=\left\vert \left(\widehat{v}^{(j)}\right)^T\widehat{v}^{(j)}\right\vert\in [0,1)$ and $W_j=\left\vert \widehat{w}^{(j)}\left(\widehat{w}^{(j)}\right)^T\right\vert\in [0,1)$. Hence, we have\begin{equation*}
	G_j(t,u)\leq \frac{2}{\sqrt{1-V_1}}\text{\ \ and\ \ }G_j(t)\leq \frac{2}{\sqrt{1-V_1}}
\end{equation*}
in (\ref{appK1}) and (\ref{appK2}), for both cases $\Lambda_j$ simple real and $\Lambda_j$ simple complex.

As a consequence, we obtain
\begin{eqnarray*}
	\max\{\epsilon(t,\widehat{y}_0),\epsilon(t)\}\leq \sum\limits_{j=2}^{q}\mathrm{e}^{\left(r_{j}-r_{1 }\right) t}\frac{f_j}{f_1} \max\left\{1,
	\frac{\left\vert \widehat{w}^{(j)}\widehat{y}_0\right\vert}{\left\vert \widehat{w}^{(1)}\widehat{y}_0\right\vert}\right\}\frac{2}{\sqrt{1-V_1}}.
\end{eqnarray*}
Now, the proof proceeds as in Theorem \ref{onset1}.
\end{proof}
	Therefore, the smaller  the right-hand side in (\ref{complexepsilon<}), the earlier the onset of the asymptotic behavior.
By comparing (\ref{realepsilon<}) and (\ref{complexepsilon<}), we observe the presence of the term $\frac{1}{2}\log\frac{1}{1-V_1}$ in (\ref{complexepsilon<}). If $V_1$ is close to $1$,  we could observe a delayed onset of the asymptotic behavior, compared to the case where $V_1$ is not close to $1$ or to the case where $\Lambda_1$ is simple real.

\subsection{Non-normal matrices}

Suppose that the vector norm is the Euclidean norm. For a normal matrix $A$, we have
$f_j=1,\ j\in\{1,\ldots,q\}$, and, in the case of  $\Lambda_1$ simple complex, $V_1=0$. On the other hand, a non-normal matrix $A$ can exhibit $f_j$ values that are arbitrarily large, as well as $V_1$ values that are arbitrarily close to $1$.

It is important to understand the impact of the non-normality of the matrix $A$ on the two aforementioned issues of asymptotic well-conditioning and onset of the asymptotic behavior.

\subsubsection{Asymptotic well-conditioning and non-normality}

	Regarding the asymptotic well conditioning, we can observe that \emph{the non-normality of $A$ has an impact only if $\Lambda_1$ is simple complex}. In fact,  if $\Lambda_1$ is simple real, or $\Lambda_1$ is simple complex and $V_1$ is not close to $1$, then the problem (\ref{due}) is asymptotically well-conditioned if and only if $\vert \widehat{w}^{(1)}\widehat{y}_0\vert^{-1}$ is not large, with the quantity $\vert \widehat{w}^{(1)}\widehat{y}_0\vert^{-1}$ unrelated to the non-normality. On the other hand, if $\Lambda_1$ is simple complex and $V_1$ is close to $1$, a case possible only for $A$ non-normal, the quantity $\vert \widehat{w}^{(1)}\widehat{y}_0\vert^{-1}=\mathrm{OSF}(y_0)$ does not determine the asymptotic well-conditioning of the problem (\ref{due}), as we have seen in the previous Subsection \ref{lambdacomplex}. In other words, for a non-normal matrix $A$, we could have an oscillating term $\mathrm{OT}(t,y_0)$ assuming large values as $t$ varies along with a nonlarge oscillation scale factor $\mathrm{OSF}(y_0)$.

\subsubsection{Onset of the asymptotic behavior and non-normality}
Regarding the onset of the asymptotic behavior, we can observe that if some ratio $\frac{f_j}{f_1}$, $j\in\{2,\ldots,q\}$, is large, or $\Lambda_1$ is simple complex and $V_1$ is close to 1, cases possible only for $A$ non-normal, \emph{the onset of the asymptotic behavior could be delayed}: recall (\ref{realepsilon<}) and (\ref{complexepsilon<}).

\begin{figure}[tbp]
	\includegraphics[width=0.8\textwidth]{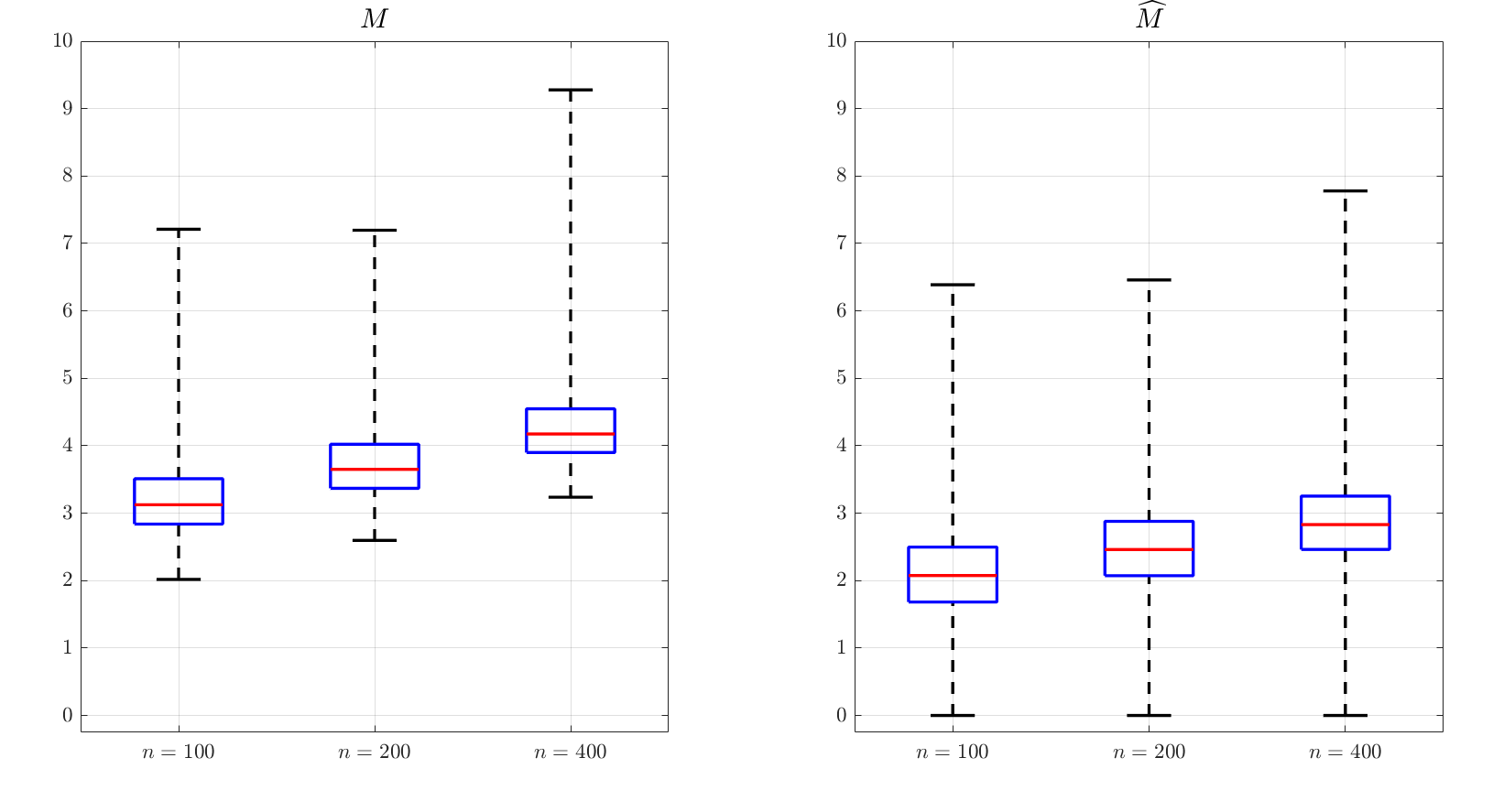}
	\caption{Box plots of $M$ and $\widehat{M}$ in (\ref{Maxima}) for $10\ 000$ random instances of $A$ with $n=100,200,400$.}
	\label{logfjf1104} 
\end{figure}

It's worth noting that the ratios $\frac{f_j}{f_1}$, $j\in\{2,\ldots,q\}$, can be much smaller than the large values of $f_j$ characterizing a high non-normality of the matrix $A$. In Figure \ref{logfjf1104}, for the Euclidean norm as vector norm, we see box-plots of
	\begin{equation}
	M=\max\limits_{j\in\{2,\ldots,n\}}\log f_j \text{\ \ and\ \ }\widehat{M}=\max\limits_{j\in\{2,\ldots,n\}}\log \frac{f_j}{f_1}, \label{Maxima}
	\end{equation}
for $10\ 000$ matrices $A$ of order $n=100,200,400$, whose elements are sampled from the standard normal distribution. It appears that the values $\widehat{M}$ are smaller than the values $M$. 

Moreover, the highly non-normal example in Subsection \ref{hnn} also highlights that the ratios $\frac{f_j}{f_1}$ can be significantly smaller than the values of $f_j$.

\subsection{Normal matrices}
Once again, assume that the vector norm is the Euclidean norm and consider $A$ normal.
 
In the case $\Lambda_1$ simple complex, we have $V_1=W_1=0$. Consequently, we have the constant over time values
$$
\mathrm{OT}(t,y_0)=\sqrt{\frac{1}{2}}\text{\ \ and\  }K_{\infty}(t,y_0)=\frac{1}{\vert \widehat{w}^{(1)} \widehat{y}_0 \vert} \sqrt{\frac{1}{2}}. 
$$
The asymptotic condition number $K_{\infty}(t,y_0)$ is constant in time, as in the case $\Lambda_1$ simple real.

Moreover, since $\frac{f_j}{f_1} = 1$ for $j \in \{2, \ldots, q\}$, and $V_1 = 0$ in the case $\Lambda_1$ simple complex, an earlier onset of the asymptotic behavior for $K(t,y_0)$ is expected compared to non-normal matrices: see (\ref{realepsilon<}) and (\ref{complexepsilon<}). 

Finally, the constant in time asymptotic condition number $K_{\infty}(t,y_0)$ is not only the limit as $t\rightarrow +\infty$ of $K(t,y_0)$, but it is also the supremum of $K(t,y_0)$ for $t\geq 0$ (see \cite{Maset2018}).

\section{Examples}\label{Examples}

To illustrate the contents of the present paper, we consider five examples of ODEs (\ref{ode}), where the first two are the models of Section 1 revisited.

\subsection{Gross Domestic Product and National Debt}

The first example is the ODE (\ref{ODEdebt}) of the GDP-ND model of Subsection \ref{twom}. Assume, as in the instance (\ref{instanceGDPND}), that there are two real eigenvalues. Moreover, since we set $Q(0)=1$, there is uncertainty only on $B(0)$ and therefore the direction of perturbation is $\widehat{z}_0=(0,\pm 1)$ (in any vector $p$-norm). In the Euclidean norm, we have
\begin{equation*}
	K_{\infty}(y_{0},\widehat{z}_{0})=\frac{\left\vert \widehat{w}^{(1)}
		\widehat{z}_{0}\right\vert}{\left\vert \widehat{w}^{(1)}\widehat{y}
		_{0}\right\vert}=\sqrt{1+B(0)^2}\cdot \frac{1}{\left\vert\frac{w^{(1)}_1}{w^{(1)}_2}+B(0)\right\vert}. 
\end{equation*}
In the instance (\ref{instanceGDPND}), we have $\frac{w^{(1)}_1}{w^{(1)}_2}=-1$ and then
\begin{equation}
	K_{\infty}(y_{0},\widehat{z}_{0})=\sqrt{1+B(0)^2}\cdot \frac{1}{\vert -1 +B(0)\vert}. \label{expression}
\end{equation}
For the initial ND $B(0)=0.61$ considered in Subsection \ref{twom}, we have $K_{\infty}(y_{0},\widehat{z}_{0})=3.0035$.

\begin{figure}[tbp]
	\includegraphics[width=0.8\textwidth]{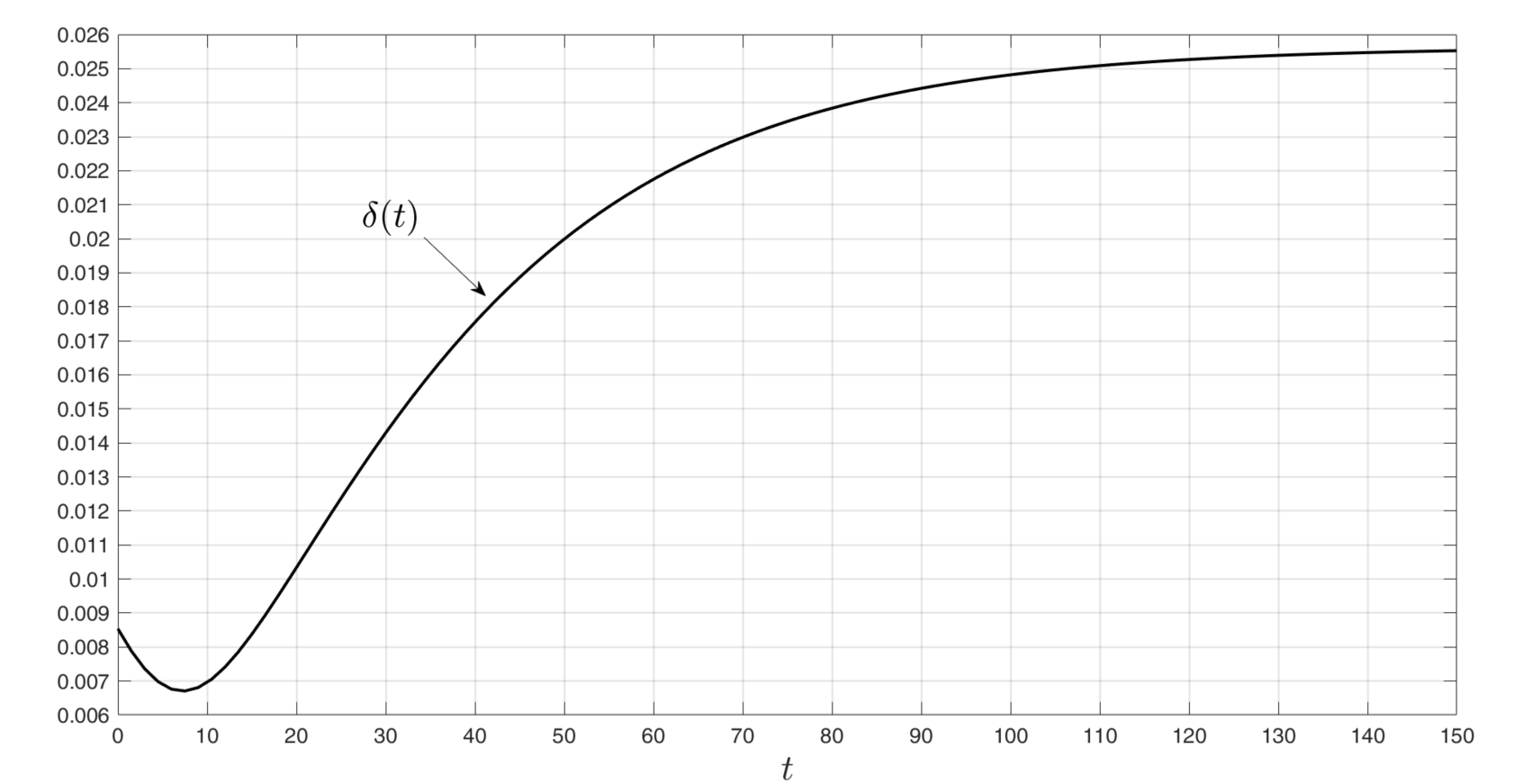}
	\caption{Relative error $\delta(t)$ for the GDP-ND model.}
	\label{Figura55} 
\end{figure}

In Figure \ref{Figura55}, which extends Figure \ref{Figura5} to the interval $[0,150]$ (the characteristic time is $20$), we can observe how the relative error $\delta(t)$ asymptotically approaches the value $0.025641$. By recalling (\ref{magnification}), we have
$$
\delta(t)=K(t,y_{0},\widehat{z}_{0})\cdot \delta(0)\text{\ \ and\ \ }\lim\limits_{t\rightarrow +\infty}\delta(t)=K_{\infty}(y_{0},\widehat{z}_{0})\cdot \delta(0)=3.0035\cdot \delta(0).
$$

The expression for $K_{\infty}(y_{0},\widehat{z}_{0})$ in (\ref{expression}) has a singularity at $B(0)=1$. In this case, where $w^{(1)}y_0=0$ holds, i.e., the initial value $(Q(0),B(0))=(1,1)$ has zero component along the eigenvector corresponding to the rightmost eigenvalue, the asymptotic condition number  $K_{\infty}(t,y_{0},\widehat{z}_{0})$ grows exponentially in time (see \cite{M1}).

\begin{figure}[tbp] 
	\includegraphics[width=0.8\textwidth]{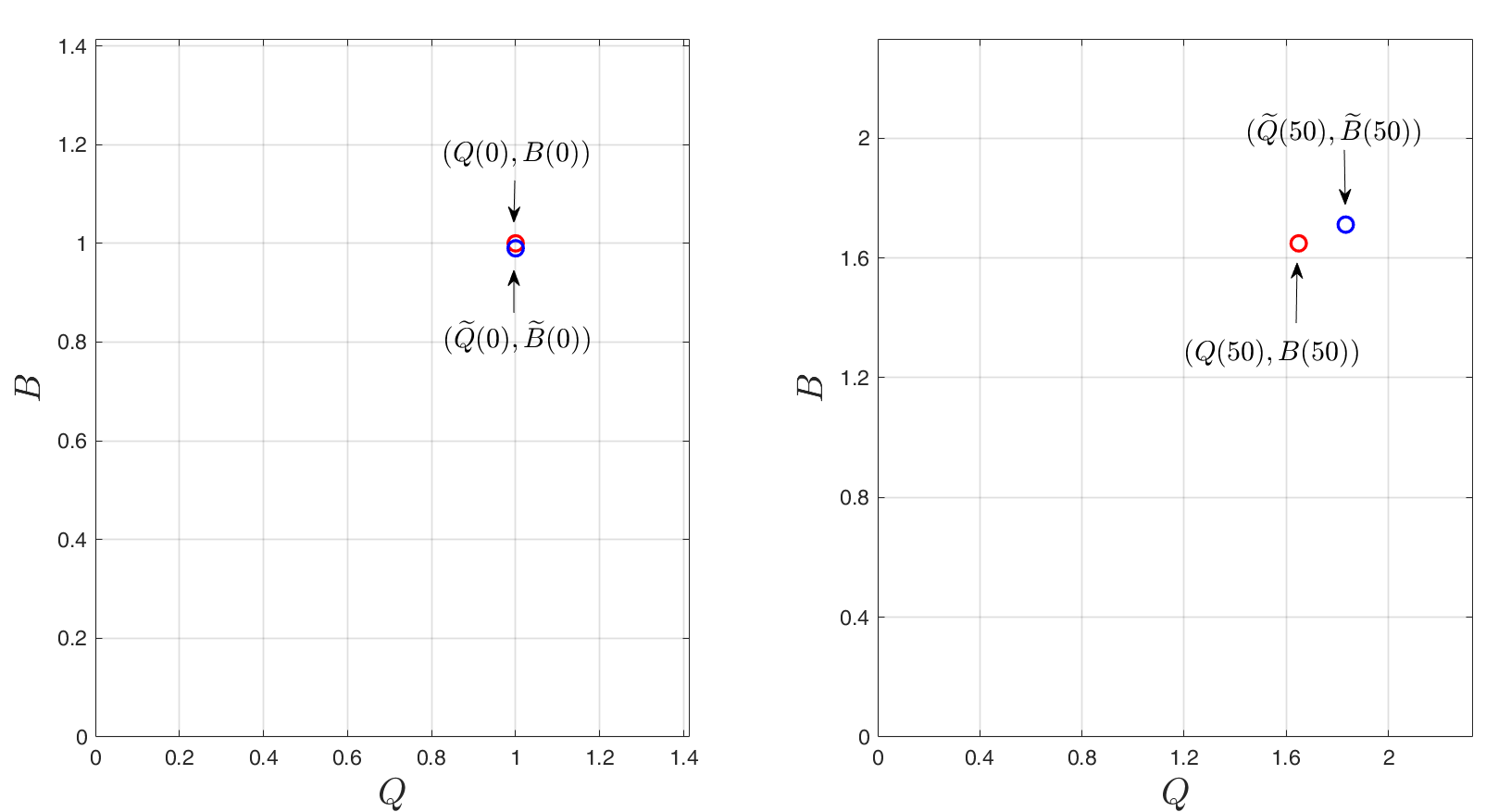}
	\caption{Singular case $w^{(1)}y_0=0$ for the GDP-ND model. Left: points $(Q(0),B(0))$ and $(\widetilde{Q}(0),\widetilde{B}(0))$. Right: points $(Q(50),B(50))$ and $(\widetilde{Q}(50),\widetilde{B}(50))$.}
	\label{QB11} 
\end{figure}

\begin{figure}[tbp]
	\includegraphics[width=0.8\textwidth]{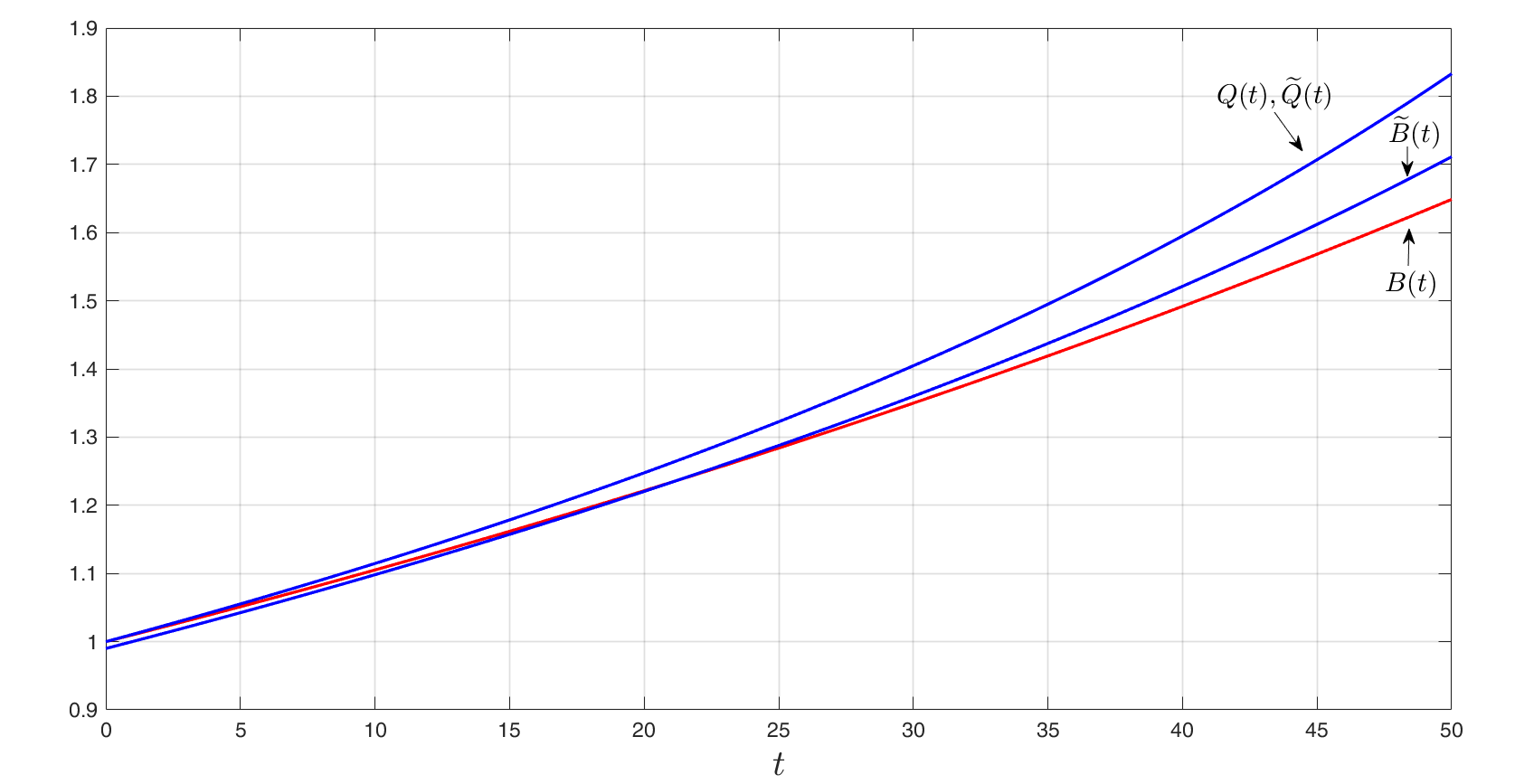}
	\caption{Singular case $w^{(1)}y_0=0$ for the GDP-ND model. GDP $Q(t)$, perturbed GDP $\widetilde{Q}(t)$, ND $B(t)$ and perturbed ND $\widetilde{B}(t)$.}
	\label{Figura311} 
\end{figure}

\begin{figure}[tbp]
	\includegraphics[width=0.8\textwidth]{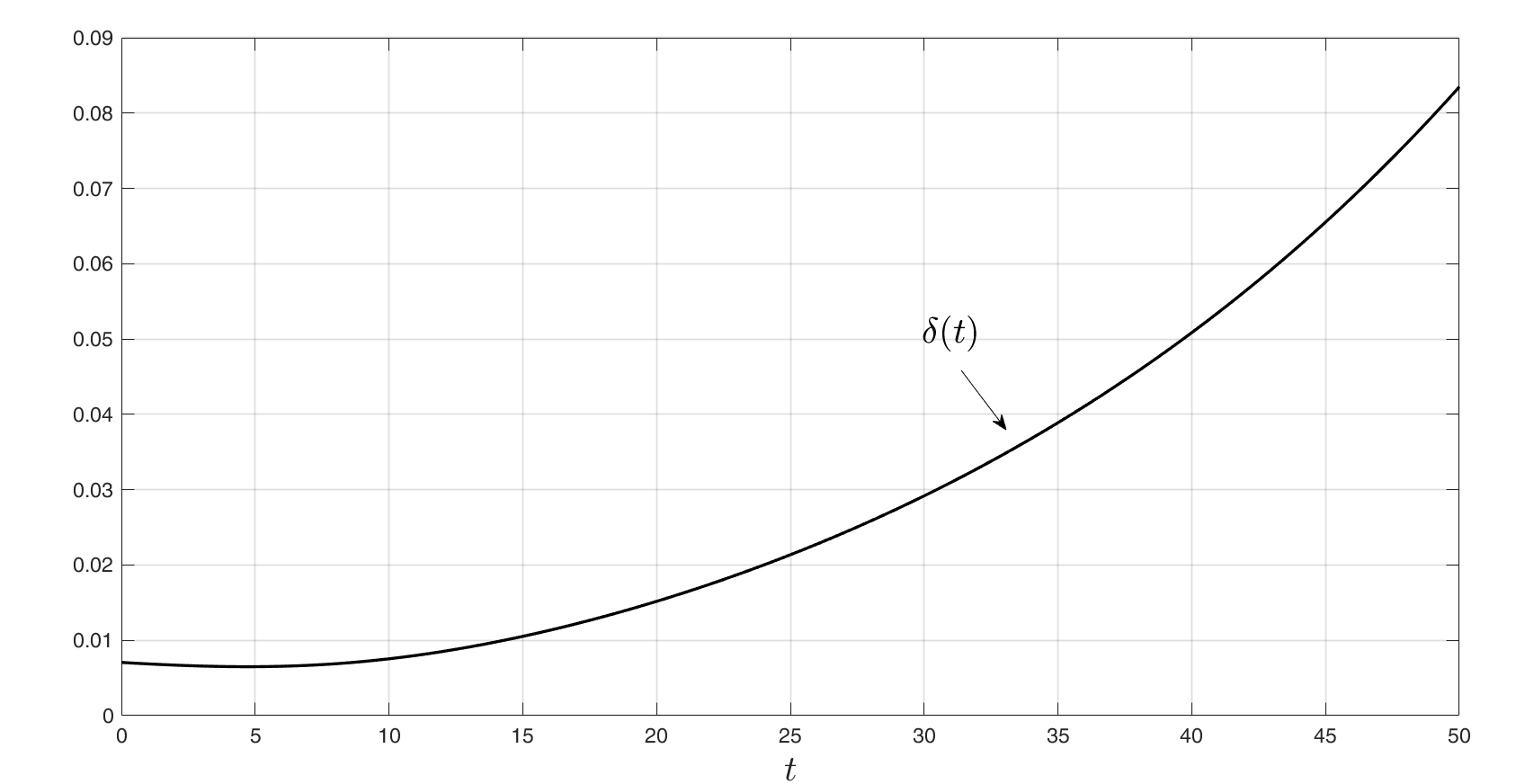}
	\caption{Singular case $w^{(1)}y_0=0$ for the GDP-ND model. Relative error $\delta(t)$.}
	\label{Figura511} 
\end{figure}

For $(Q(0),B(0))=(1,1)$ and $(\widetilde{Q}(0),\widetilde{B}(0))=(1,0.99)$, Figures \ref{QB11}, \ref{Figura311}, and \ref{Figura511} depict the same information as Figures \ref{QB2}, \ref{Figura3}, and \ref{Figura5} do for $(Q(0),B(0))=(1,0.61)$ and $(\widetilde{Q}(0),\widetilde{B}(0))=(1,0.60)$. The differences between the two cases are quite evident. After $50\ \mathrm{yr}$, the relative error  is about $12$ times the initial relative error.

\subsection{Building heating}

The second example is the ODE (\ref{ODEhouse}) of the building heating model of Subsection \ref{buildingheating}. Assume, as in the instance (\ref{instancehb}), that there are three (negative) eigenvalues. We have
\begin{equation}
K_{\infty}(y_{0},\widehat{z}_0)=\frac{\left\vert \widehat{w}^{(1)}\widehat{z}_0\right\vert}{\left\vert \widehat{w}^{(1)}\widehat{y}_0\right\vert}\text{\ and\ }K_{\infty}(y_{0})= \frac{1}{\left\vert \widehat{w}^{(1)}\widehat{y}_0\right\vert}.
\label{expressionbh}
\end{equation}
For the particular instance (\ref{instancehb}) and the Euclidean norm as the vector norm, we have the normalized left eigenvector
\begin{equation}
\widehat{w}^{(1)}=
\left[
\begin{array}{c}
	-0.4462\   -0.6111\   -0.6538
\end{array}
\right]. \label{wsamesign}
\end{equation}
corresponding to the rightmost eigenvalue. For the initial value and perturbed initial value in (\ref{y0pert}) and (\ref{y0}), we have $K_{\infty}(y_{0},\widehat{z}_0)$\ $= 11.8648$.

In Figure \ref{relativebh}, we observe the relative error $\delta(t) = K(t, y_0, \widehat{z}_0)\cdot \delta(0)$ over approximately two characteristic times (the characteristic time is $3.17 \ \mathrm{h}$). The relative error asymptotically approaches   
$$
\lim\limits_{t\rightarrow +\infty}\delta(t)=K_\infty(y_0,\widehat{z}_0)\cdot \delta(0)=11.8648\cdot\delta(0)
$$
as $t$ increases. Moreover, we have $K_{\infty}(y_{0}) = 12.1330$. This indicates that the perturbation in (\ref{y0}) is nearly the worst-case scenario.

The expression for $K_\infty(y_0)$ in (\ref{expressionbh}) has a singularity when  $\widehat{w}^{(1)}y_0=0$. For the instance  (\ref{instancehb}), this occurs when (\ref{y0pert}) is replaced with
\begin{equation}
	y_0=\left(3.5\ ^{\circ}\mathrm{C}, -5.2298\ ^{\circ}\mathrm{C}, 2.5\ ^{\circ}\mathrm{C}\right). \label{y0sing}
\end{equation}
In this case, the asymptotic condition number  $K_{\infty}(t,y_{0})$ grows exponentially in time (see \cite{M1}).

\begin{figure}[tbp]
	\includegraphics[width=0.8\textwidth]{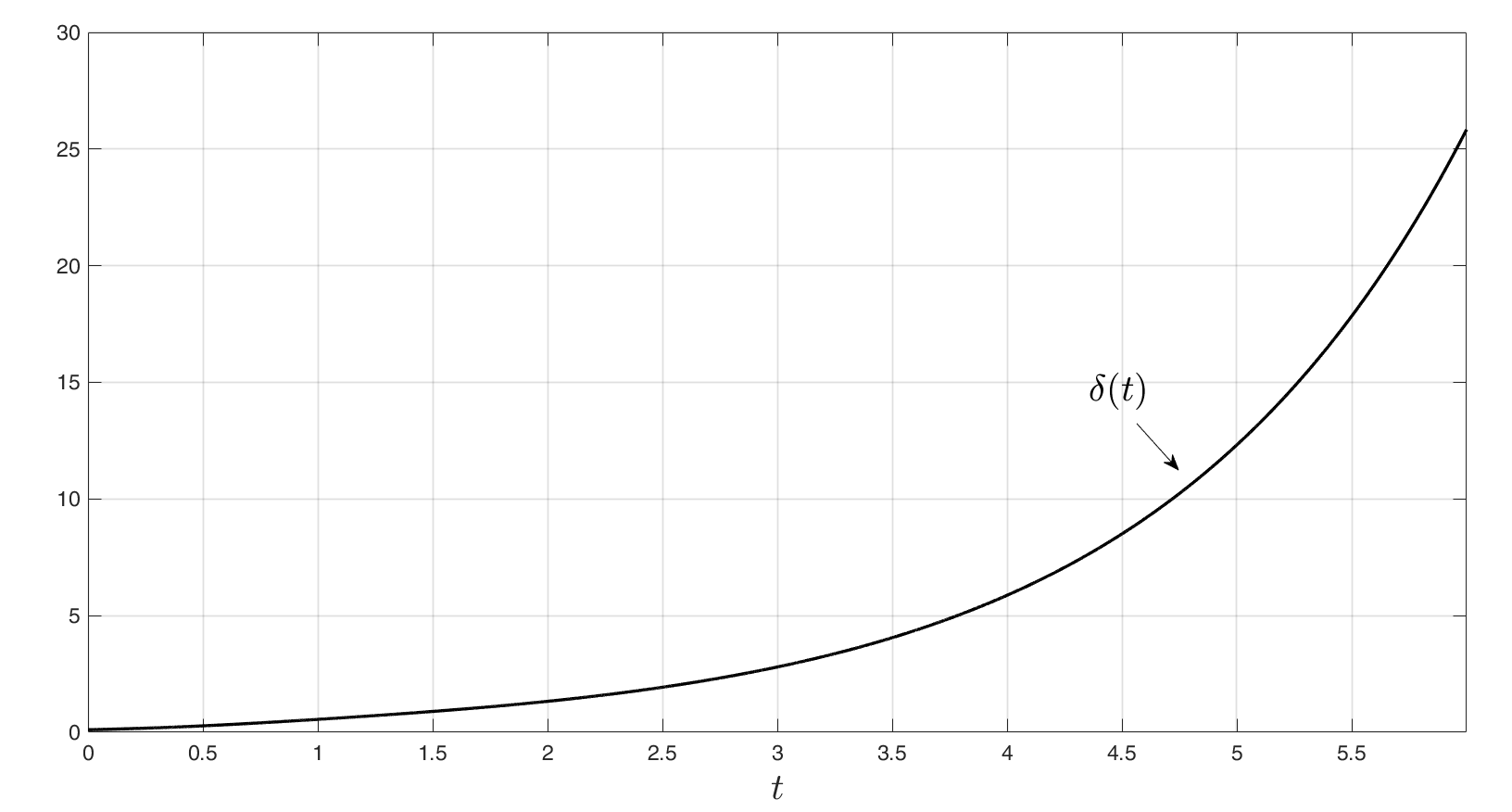}
	\caption{Singular case $w^{(1)}y_0=0$ for the heating building model. Relative error $\delta(t)$.}
	\label{Figuresingularityhb}
\end{figure}

In Figure \ref{Figuresingularityhb}, for $y_0$ given in (\ref{y0sing}) and
\begin{equation*}
	\widetilde{y}_0=\left(4\ ^{\circ}\mathrm{C}, -5\ ^{\circ}\mathrm{C}, 3\ ^{\circ}\mathrm{C}\right), \label{y0tildesing}
\end{equation*}
we observe the relative error $\delta(t)$ for $t \in [0, 6\ \mathrm{h}]$. The difference compared to Figure \ref{relativebh} is quite evident. After $6\ \mathrm{h}$, the relative error is 235 times the initial relative error.
\begin{remark}
In the instance (\ref{instancehb}), all the components of $\widehat{w}^{(1)}$ in (\ref{wsamesign}) have the same sign. This fact permits the following observation. If all the components of $y_0$ have the same sign, then 
$$
K_{\infty}(y_{0})= \frac{\Vert y_0\Vert_2}{\left\vert \widehat{w}^{(1)}_1\right\vert\vert y_{01}\vert+\left\vert \widehat{w}^{(1)}_2\right\vert\vert y_{02}\vert+\left\vert \widehat{w}^{(1)}_3\right\vert\vert y_{03}\vert}
$$
and then
$$
K_{\infty}(y_{0})\leq \frac{1}{\min\left\{\left\vert \widehat{w}^{(1)}_1\right\vert,\left\vert \widehat{w}^{(1)}_2\right\vert,\left\vert \widehat{w}^{(1)}_3\right\vert	\right\}}=\frac{1}{0.4462}=2.2411.
$$
This shows that if the temperatures of the three floors are all larger or all smaller than the equilibrium temperatures, then the relative error of the perturbed solution is not much larger than the initial relative error. This is not the case for the initial values (\ref{y0pert}) and (\ref{y0sing}), where the relative error grows significantly compared to its initial value.
\end{remark}

\begin{remark}\label{remark}

In the building heating model, as well as in the GDP-ND model, we have considered $y_0$ as the actual initial value and $\widetilde{y}_0$ as the initial value available to us for the simulation.  Therefore, we have considered $\delta(t)$ as the relative error of the simulated solution with respect to the actual solution. From this perspective, $y_0$ is the unknown actual initial value, and only an approximation is known, which is regarded as a perturbed value $\widetilde{y}_0$ of $y_0$.

However, to obtain a computable asymptotic condition number $K_\infty(y_0)$ for the heating building model, or a computable asymptotic condition number $K_\infty(y_0,\widehat{z}_0)$ for the GDP-ND model, one may consider $y_0$ as the available initial value and $\widetilde{y}_0$ as the actual initial value. Therefore, now we are considering $\delta(t)$ as the relative error of the actual solution with respect to the simulated solution. From this new perspective, $y_0$ is the available initial value, and the unknown actual initial value  is regarded as a perturbed value $\widetilde{y}_0$ of $y_0$.

In the case of the  building heating model, we have $K_\infty(y_0)=4.9195$ by considering the available initial value (\ref{y0}) as $y_0$. Therefore, by assuming absolute errors in the initial temperatures of the three floors not larger than $0.5 \ ^{\circ}\mathrm{C}$ , we have $\delta(0)\leq 0.1353$ and then
$$
\lim\limits_{t\rightarrow +\infty}\delta(t)\leq 0.5670.
$$

In the case of the GDP-ND model, we have $K_\infty(y_0,\widehat{z}_0)=2.29$ for $B(0)=0.60$, the available initial value. Therefore, by assuming an absolute error in the initial ND not larger than $0.01$, we have $\delta(0)\leq 0.0086$ and then
$$
\lim\limits_{t\rightarrow +\infty}\delta(t)\leq 0.0196.
$$
\end{remark}

\subsection{Another building heating model}

The paper \cite{Tol2023} presents a building heating model of dimension $4$, where the first state component is the indoor temperature and the other three are the temperatures in three layers of the wall separating the indoor space from the outdoor space. The matrix of the system of differential equations is
\begin{equation}
	A=\left[ 
	\begin{array}{rrrr}
		-5.7215   &    5.7215    &        0      &       0 \\
		0.23076   & -0.39276     &      0.162    &        0  \\
		0      &    0.081     &     -0.162    &    0.081 \\
		0    &       0      &     0.162     & -0.91116 
	\end{array}
	\right] \label{matrixA}
\end{equation}
when the time is measured in hours. The matrix has the four real negative eigenvalues $-0.038938\ \mathrm{h}^{-1}$, $-0.26104\ \mathrm{h}^{-1}$, $-0.92864\ \mathrm{h}^{-1}$ and  $-5.9588\ \mathrm{h}^{-1}$. The characteristic time is $25.7\ \mathrm{h}$. The left eigenvector corresponding to the rightmost eigenvalue is
$$
w^{(1)}=\left[-0.029943\     -0.73735\      -1.1058\      -0.1027\right]
$$

\begin{figure}[tbp]
	\includegraphics[width=0.8\textwidth]{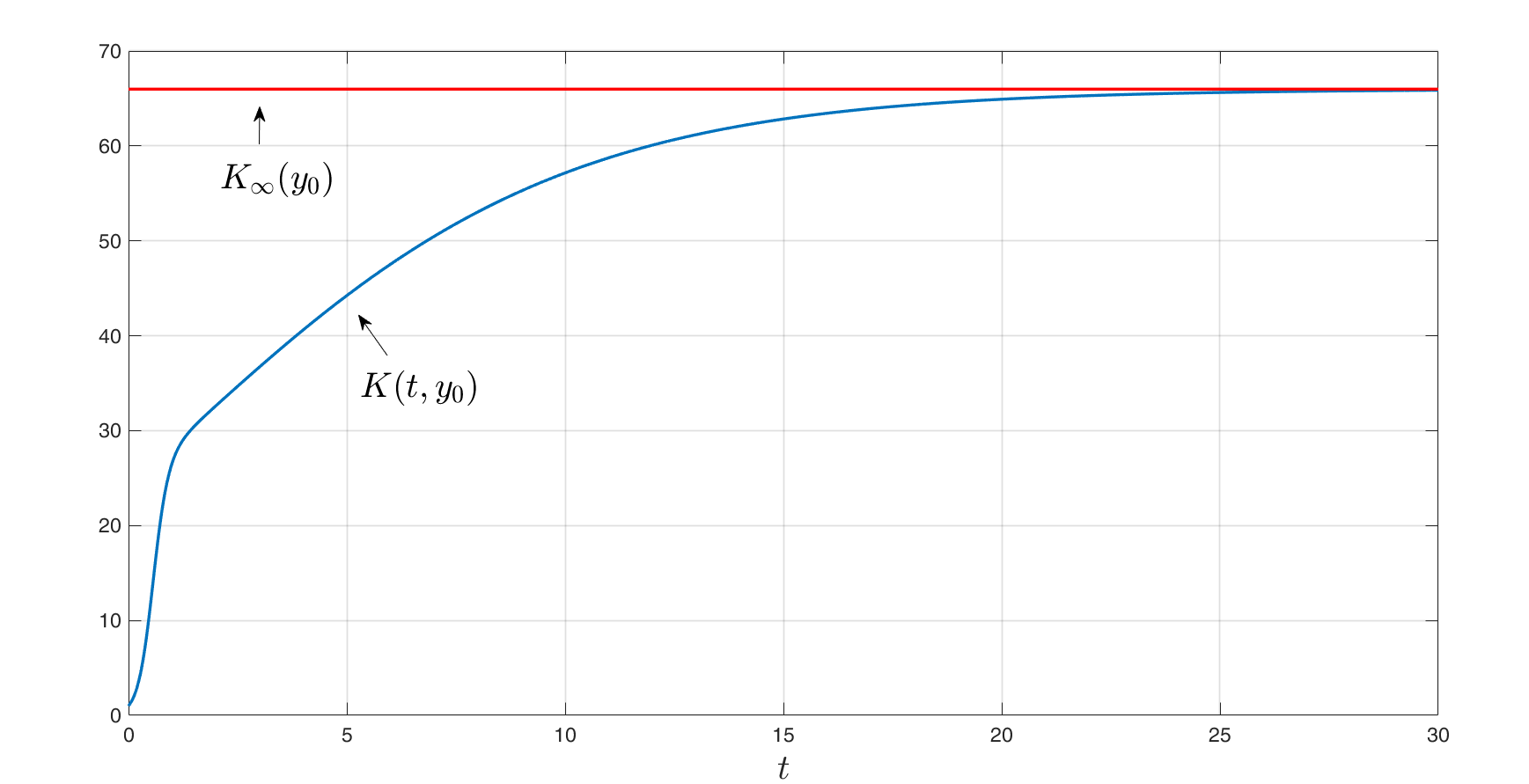}
	\caption{Condition number $K(t,y_0)$ and asymptotic condition number $K_\infty(y_0)$  for the matrix (\ref{matrixA}) and $y_0$ with $y_{02}=y_{03}=y_{04}=0$. The vector norm is the $\infty$-norm.}
	\label{Figure222}
\end{figure}

Consider initial temperatures $y_0$ (remind that these are temperatures with respect to the equilibrium temperatures) with $y_{02}=y_{03}=y_{04}=0$. This corresponds to the case where the uncertainty in the indoor temperature is much greater than that in the layer temperatures. We have, for a vector $p$-norm,  
$$
K_\infty(y_0)=\frac{1}{\vert \widehat{w}^{(1)}\widehat{y}_0\vert}=\frac{\Vert w^{(1)}\Vert}{\vert w^{(1)}_1\vert}.
$$
In Figure \ref{Figure222}, for the $\infty$-norm as the vector norm, we see the condition number $K(t,y_0)$, $t\in [0,30\ \mathrm{h}]$, and  the asymptotic condition number $K_\infty(y_0)=65.987$. The asymptotic behavior of the condition number is achieved within a characteristic time.

\subsection{Charged particle subject to a magnetic field with viscous frictional force}
The first three examples described situations in which the matrix $A$ in (\ref{ode}) had real eigenvalues. This fourth example examines a situation in which the rightmost eigenvalues of $A$ form a complex conjugate pair.

We consider a charged particle subject to a magnetic field and moving in a medium with an anisotropic viscous frictional force (see the paper \cite{Butanas2022}, which also considers Brownian motion and an electric field). Let $m$ and $q$ be the mass and the charge, respectively, of the particle, let $B=(B_x,B_y,B_z)$ be the constant magnetic field and let $\Gamma_x$, $\Gamma_y$ and $\Gamma_z$ be the constant coefficients of the viscous frictional force in the three spatial directions. The motion of the particle is described by the first order linear ODE for its velocity $v=(v_x,v_y,v_z)$:
\begin{equation*}
	v^{\prime}(t)=Av(t),
\end{equation*}
where
\begin{equation}
	A=\frac{1}{m}\left[ 
	\begin{array}{rrr}
		-\Gamma_x   &   qB_z    &       -qB_y    \\
		-qB_z   & -\Gamma_y    &     qB_x   \\
		qB_y      &    -qB_x   &     -\Gamma_z 
	\end{array}
	\right].  \label{Amf}
\end{equation}
We examine the instance (comparable with the examples given in \cite{Butanas2022}) where, in suitable units,
\begin{eqnarray*}
	&&m=1\\
	&&qB_x=-0.5,\ qB_y=-1.3,\ qB_z=-0.4,\\
	&&\Gamma_x=0.4,\ \Gamma_y=0.8,\ \Gamma_z=0.2.
	\label{instance}
\end{eqnarray*}
In such an instance, the matrix $A$ has the eigenvalues
$$
\lambda_1=-0.3433 +1.4326\mathrm{i},\ \overline{\lambda_1}=-0.3433 -1.4326\mathrm{i},\ \lambda_2=-0.7133,
$$
with a complex conjugate pair $\lambda_1$ and $\overline{\lambda_1}$ as rightmost eigenvalues.

In the Euclidean norm, the normalized left eigenvector corresponding to $\lambda_1$ is
$$
\widehat{w}^{(1)}=\left[\begin{array}{rrr} 0.6543 + 0.0111\mathrm{i} &-0.2821 + 0.1288 \mathrm{i}&-0.1031 -0.6819\mathrm{i}\end{array}\right].
$$ 
The oscillation scale factor of the asymptotic condition number $K_\infty(t,v_0)$ is 
$$
\mathrm{OSF}(v_0)=\frac{1}{\left\vert \widehat{w}^{(1)}\widehat{v}_0\right\vert}
=\frac{\Vert v_0\Vert_2}{\sqrt{\left(\mathrm{Re}(\widehat{w}^{(1)})v_0\right)^2+\left(\mathrm{Im}(\widehat{w}^{(1)})v_0\right)^2}},
$$
where $v_0$ is the initial velocity.

\begin{figure}[tbp]
	\includegraphics[width=0.8\textwidth]{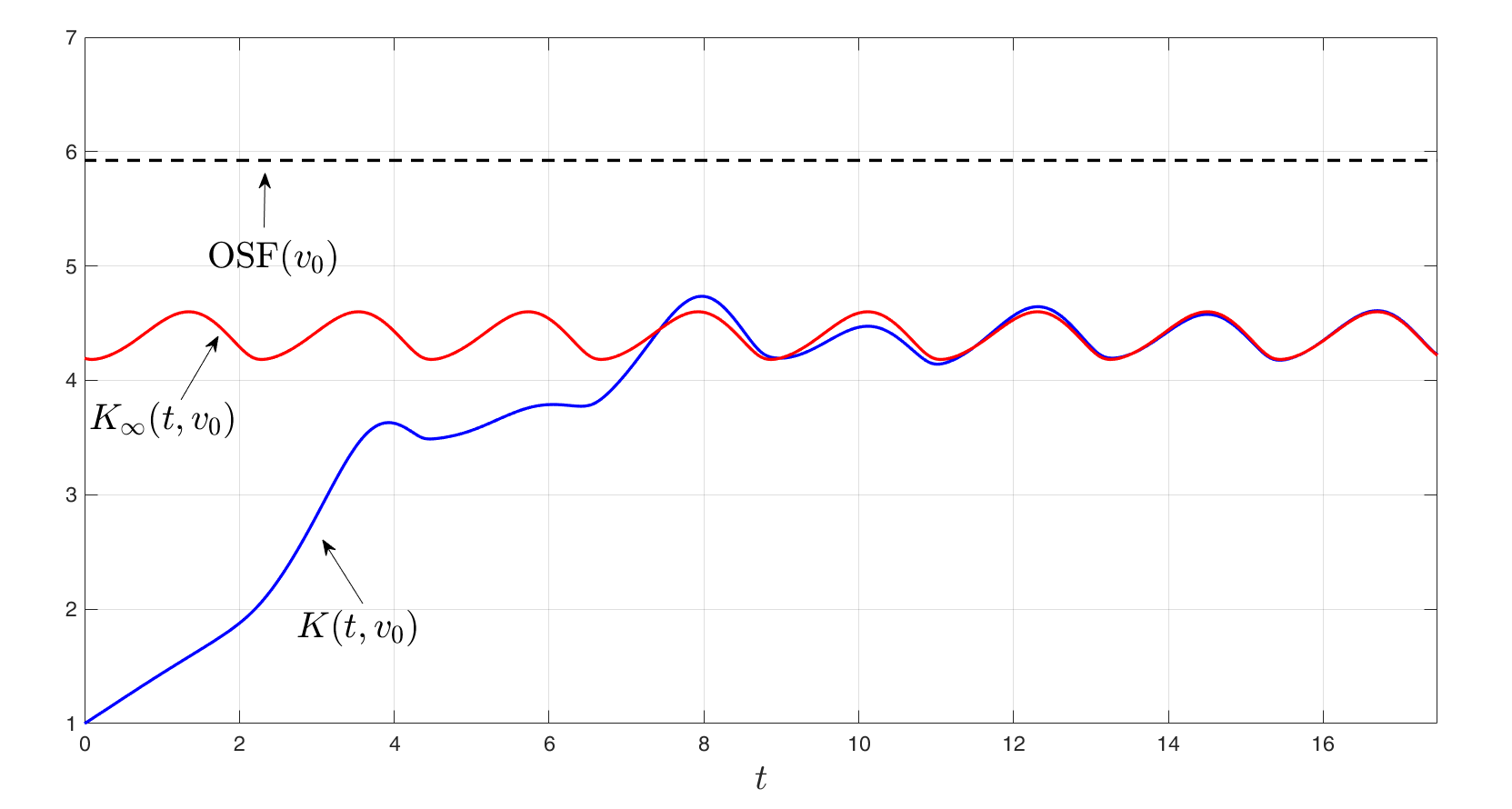}
	\caption{Condition number $K(t,v_0)$, asymptotic condition number $K_\infty(t,v_0)$ and oscillating scale factor $\mathrm{OSF}(v_0)$ for the matrix (\ref{Amf}) and $v_0=(0.5,1,0.5)$. The vector norm is the Euclidean norm.}
	\label{Figmagf}
\end{figure}

In Figure \ref{Figmagf}, we see, for $v_0=(0.5,1,0.5)$, the condition number $K(t,v_0)$ and the asymptotic condition number $K_\infty(t,v_0)$, $t$ over $6$ characteristic times $\frac{1}{-r_1}=2.9125$, and the oscillation scale factor $\mathrm{OSF}(v_0)=5.92$. The asymptotic behavior is achieved in approximately three characteristic times.

Since $V_1=0.0587$ and $W_1=0.0937$ in this example, the oscillations of $K_\infty(t,v_0)$ are tight for any $v_0$: for the oscillating term $\mathrm{OT}(t,y_0)$, we have (recall Theorem \ref{theorem9z})
$$
a^{\min}_{V_1W_1}=0.6973\text{\ \ and\ \ }a^{\max}_{V_1W_1}=0.7842.
$$

\subsection{A highly non-normal matrix $A$} \label{hnn}
In this final example, we consider in the ODE (\ref{ode}) the highly non-normal matrix $A=VDV^{-1}$, where 
$V$ is the Hilbert matrix of order $8$ and
$$
D=\mathrm{diag}(-0.1,-0.2,\ldots,-0.8).
$$
The Euclidean norm is used as the vector norm.

The high non-normality of the matrix $A$ appears in the large values $f_j$ defined in (\ref{fj}): we
have 
\begin{eqnarray*}
	&&f_1=5.2554\cdot 10^5,\ f_2=1.677\cdot 10^7,\ f_3=1.6347\cdot 10^8,\   f_4=7.1819\cdot 10^8,\\
	&&f_5=1.6407\cdot 10^9,\ f_6=2.0252\cdot 10^9,\ f_7=1.2815\cdot 10^9,\ f_8=3.2603\cdot 10^8.
\end{eqnarray*}

\begin{figure}[tbp]
	\includegraphics[width=0.8\textwidth]{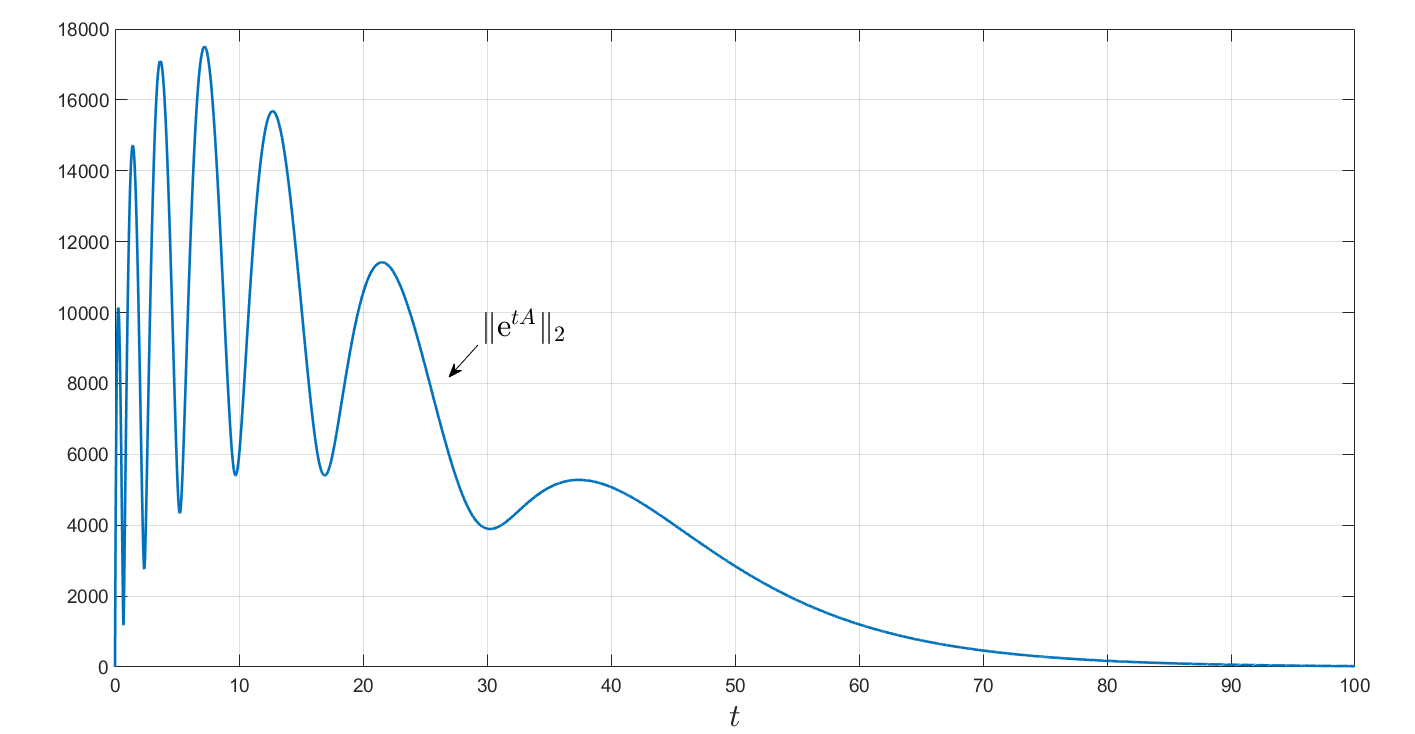}
	\caption{Norm $\Vert \mathrm{e}^{tA}\Vert_2$ for $A=VDV^{-1}$, where $V$ is the Hilbert matrix of order $8$.}
	\label{absolute}
\end{figure}

In Figure \ref{absolute}, we see the norm $\Vert \mathrm{e}^{tA}\Vert_2$ for $t \in [0,100]$, i.e., for $t$ over ten characteristic times $\frac{1}{-r_1} = 10$. Note that $\Vert \mathrm{e}^{tA}\Vert_2$ is the absolute condition number of the problem (\ref{due}), i.e., when we consider absolute errors instead of relative errors. Despite $\Vert \mathrm{e}^{tA} \Vert_2\rightarrow 0$, $t\rightarrow +\infty$, this fact holds little significance. This is because, owing to the high non-normality of the matrix $A$, the norm remains large for numerous characteristic times and it only drops below $1$ after 13 characteristic times.

Hence, it is of little interest to ascertain that the absolute error of a perturbation in the initial value decays to zero in the long-time, as this ``long-time'' may be very distant. What we observe within a possible time span of interest is a very large amplification of the absolute error of the perturbation.

\begin{figure}[tbp]
	\centering
	\begin{subfigure}[b]{1\textwidth}
		\centering
		\includegraphics[width=0.8\textwidth]{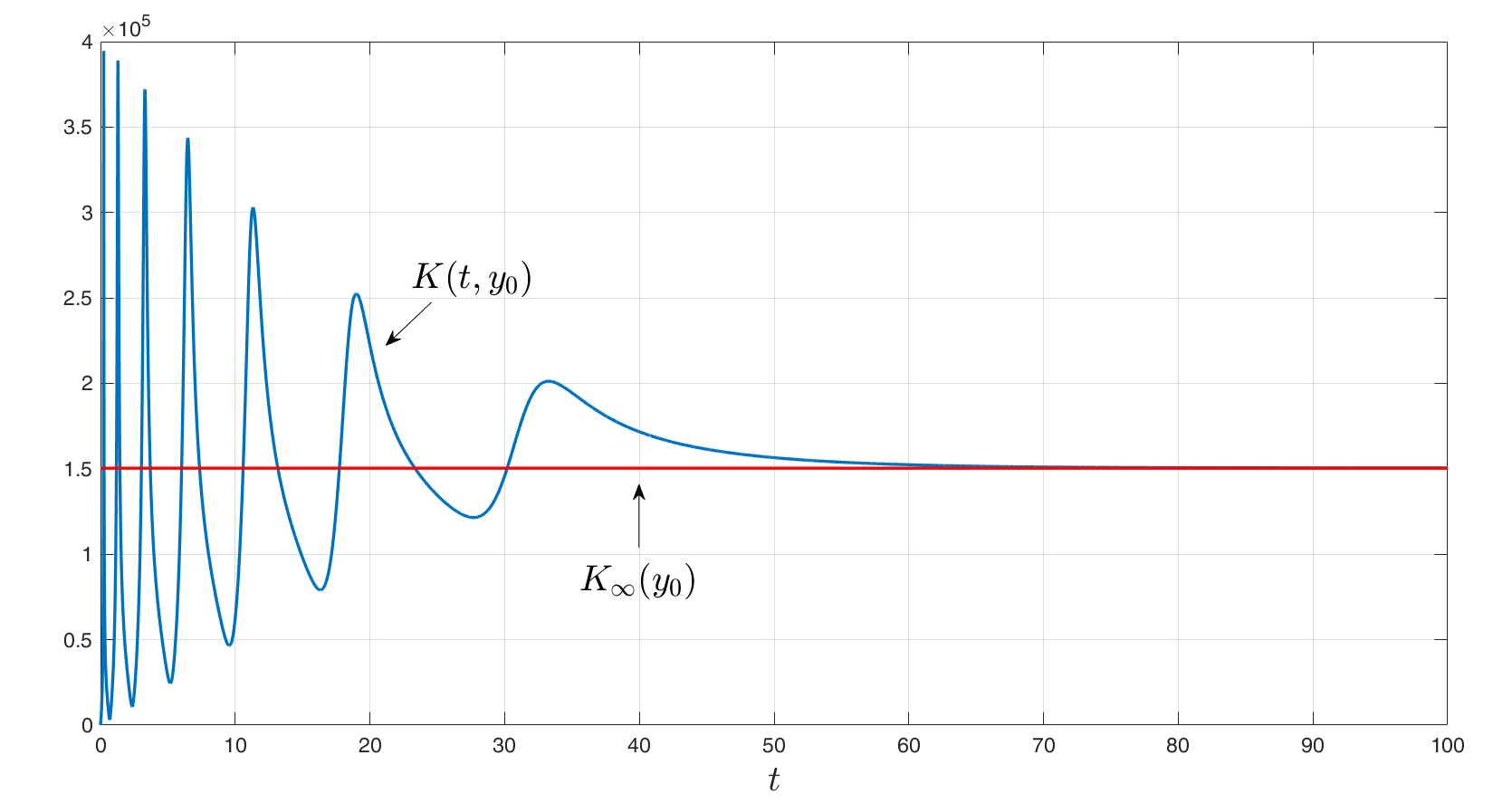}
		\caption{$y_0=(1,1,1,1,1,1,1,1)$}
	\end{subfigure}
	\hfill  
	\begin{subfigure}[b]{1\textwidth}
		\centering
		\includegraphics[width=0.8\textwidth]{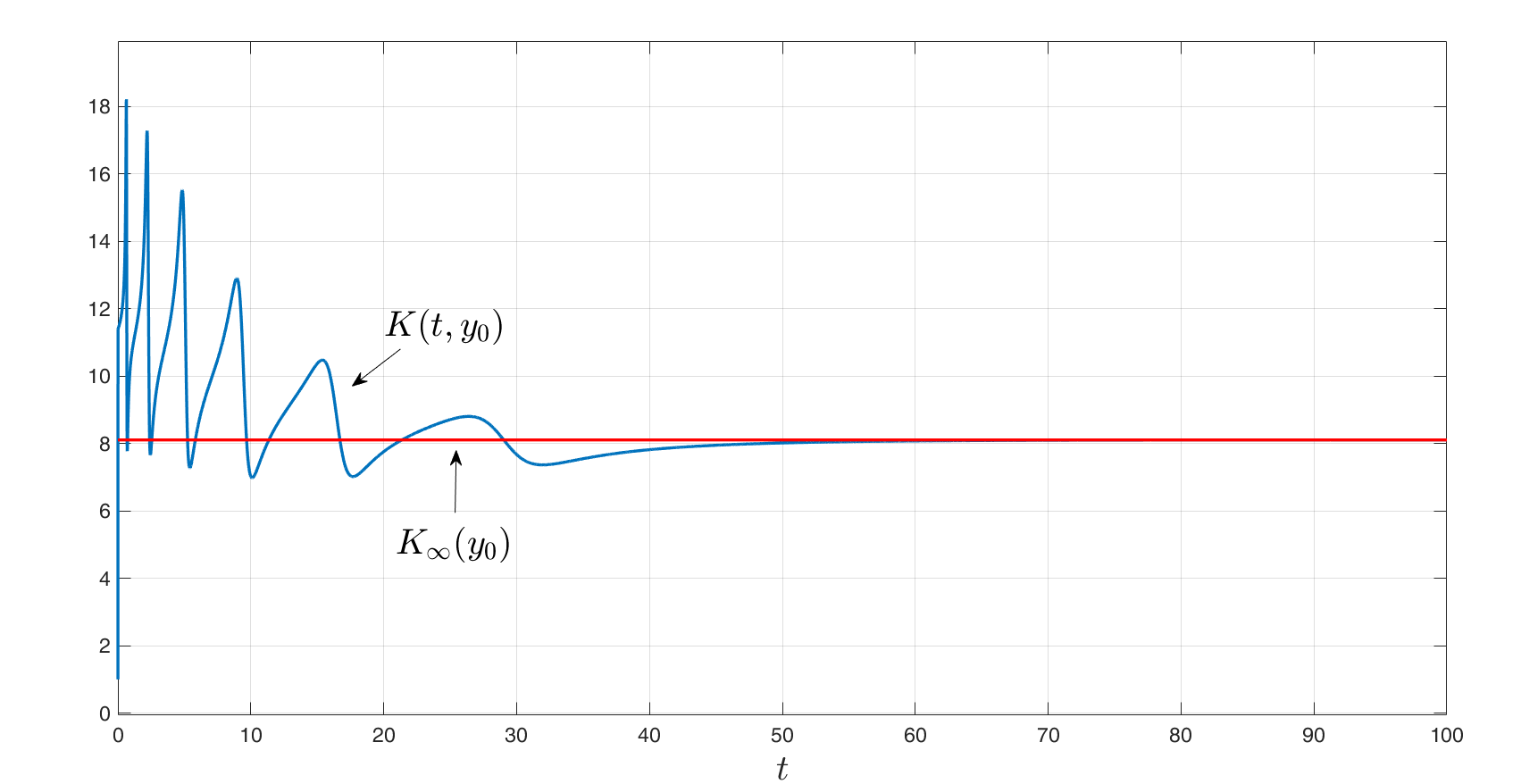}
		\caption{$y_0=(1,1,1,1,-1,-1,-1,-1)$}
	\end{subfigure}
	\hfill
	\caption{Condition numbers $K(t,y_0)$ and $
		K_\infty(y_0)$ for $A=VDV^{-1}$, where $V$ is the Hilbert matrix of order $8$. The vector norm is the Euclidean norm.}
	\label{Figure6}
\end{figure}

On the other hand, in Figure \ref{Figure6} we see, for two different initial values $y_0$ and for $t\in[0,100]$, 
the relative condition number $K(t,y_0)$ of the problem (\ref{due}), i.e., when we consider relative errors instead of absolute errors, and its asymptotic value
$
K_\infty(y_0)=\frac{1}{\vert \widehat{w}^{(1)}\widehat{y}_0\vert}.
$
Similar pictures are obtained
by randomly varying the initial value.

Note that the high non-normality does not have a critical impact, as is the case with the absolute condition number $\Vert \mathrm{e}^{tA}\Vert_2$. In fact, we can observe the following.
\begin{itemize}
	\item [1)] Although $K(t,y_0)$ oscillates similarly to $\Vert \mathrm{e}^{tA}\Vert_2$ as $t$ varies, the maximum of $K(t,y_0)$ as $t$ varies and $K_\infty(y_0)$ have the same order of magnitude. The high value of $K_\infty(y_0)$ in the top figure is due to a small value of $\vert \widehat{w}^{(1)}\widehat{y}_0\vert$, rather than high non-normality.
	\item [2)] The constant asymptotic behavior is achieved after a few (approximately five) characteristic times, similar to many examples where the matrix is not highly non-normal. For comparison, consider the values of $\Vert \mathrm{e}^{tA} \Vert_2$ in Figure \ref{absolute} after five characteristic times. Indeed, the fact that the high non-normality does not critically impact the attainment of asymptotic behavior can be explained by observing that, when the characteristic time $\frac{1}{-r_1}$ is used as the time unit $\widehat{t}$, we have
	$$  
	\max\limits_{j\in\{2,\ldots,8\}}\frac{1}{(r_1-r_j)\widehat{t}}\log\frac{f_j}{f_1}=\max\limits_{j\in\{2,\ldots,8\}}
	\frac{1}{\underset{=j}{\underbrace{\frac{-r_j}{-r_1}}}-1}\log\frac{f_j}{f_1}=3.4629
	$$
	in (\ref{realepsilon<}). Therefore, the high non-normality affects the achievement of the asymptotic behavior for only $3.4629$ characteristic times.
\end{itemize}

Whereas the fact that the absolute condition number $\Vert \mathrm{e}^{tA}\Vert_2$ decays to zero in the long-time is not significant, the fact that the relative condition number $K(t,y_0)$  asymptotically becomes $K_\infty(y_0)$, and its order of magnitude never exceeds that of $K_\infty(y_0)$, is highly significant.

\section{Asymptotic and non-asymptotic conditioning} \label{NAsy}
In this paper, we have studied the asymptotic conditioning of the problem (\ref{due}), i.e. the conditioning for large $t$. On the other hand, we are also interested in the conditioning of (\ref{due}) for any $t$, not only for large $t$.

Of course, if the problem (\ref{due}) is \emph{asymptotically ill-conditioned} at $y_0$, i.e., $K_\infty(t,y_0)$ assumes large values as $t$ varies, then the problem (\ref{due}) is also \emph{ill-conditioned} at $y_0$, i.e., $K(t,y_0)$ assumes large values as $t$ varies.

Vice versa, of course, we cannot say that if the problem (\ref{due}) is \emph{asymptotically well-conditioned} at $y_0$,  i.e., $K_\infty(t,y_0)$ does not assume large values as $t$ varies, then the problem (\ref{due}) is also \emph{well-conditioned} at $y_0$, i.e., $K(t,y_0)$ does not assume large values as $t$ varies. In fact, the asymptotic behavior  $K_\infty(t,y_0)$ of $K(t,y_0)$ cannot take into account an initial increase of $K(t,y_0)$ to values much larger than those of $K_\infty(t,y_0)$.

However, we now show experimentally that such an initial increase in $K(t,y_0)$ is very rare, and then the asymptotic well conditioning can indeed be considered equivalent to the well conditioning, for the problem (\ref{due}). This stresses the crucial and central importance of the asymptotic condition number.

We consider $10\ 000$ instances of an ODE (\ref{ode}) of dimension $n=5,25,100$, where the entries of $A$ and the components of the initial value $y_0$ are sampled from the standard normal distribution. Of these $10\ 000$ instances, $5000$ have a matrix $A$ with a rightmost real eigenvalue and $5000$ have a matrix $A$ with a rightmost pair of complex conjugate eigenvalues. Since these are random instances, the eigenvalues are always simple: therefore, we have $5000$ instances with $\Lambda_1$ simple real and  $5000$ instances with $\Lambda_1$ simple complex.
 
 We compute the ratios
 \begin{equation}
 R=\frac{\max\limits_{t\in[0,T]}K(t,y_0)}{ \max\limits_{t\in[0,T]}K_\infty(t,y_0)} \label{ratioR}
\end{equation}
for these $10\ 000$ instances, where $T$ is fifty times the characteristic time $\widehat{t}=\frac{1}{\vert r_1\vert}$. The maxima are determined using $1000$ equally spaced sampling points over the interval $[0,T]$, i.e. with a stepsize $\frac{1}{20}$ of the characteristic time. Observe that, for $\Lambda_1$ simple real, we have 
\begin{equation*}
	R=\frac{ \max\limits_{t\in[0,T]}K(t,y_0)}{K_\infty(y_0)},
\end{equation*}
since $K_\infty(t,y_0)=K_\infty(y_0)$ remains constant over time $t$.

We see in Tables \ref{tablereaL} and \ref{tablecomplex} percentiles of the ratio $R$ for the $5000$ instances with $\Lambda_1$ simple real and the $5000$ instances with $\Lambda_1$ simple complex, respectively. Euclidean norm, $\infty$-norm and $1$-norm are considered as vector norm.

\begin{table}
	
	Euclidean Norm\\
	\begin{tabular}{|l|l|l|l|l|l|}
		\hline
		& median & 90th & 99th & 99.9th & maximum value\\
		\hline
		$n=5$ & 1.0174 & 1.2905 & 2.1235 & 3.3036 & 6.0996\\
		\hline
		$n=25$ & 1.0394 & 1.3385 & 2.1400 & 4.8167 & 5.6612\\
		\hline
		$n=100$ & 1.0348 & 1.3238 & 1.9519 & 3.0598 & 4.6122\\
		\hline	  
	\end{tabular}
	\quad \\
	\quad\\
	$\infty$-norm\\
	\begin{tabular}{|l|l|l|l|l|l|}
		\hline
		& median & 90th & 99th & 99.9th & maximum value\\
		\hline
		$n=5$ & 1.0698 & 1.5358 & 2.5534 & 4.2668 & 6.1588\\
		\hline
		$n=25$ & 1.1202 & 1.6013 & 2.6554 & 4.6707 & 7.4487\\
		\hline
		$n=100$ & 1.1113 & 1.5331 & 2.3101 & 3.6398 & 4.7976\\
		\hline	  
	\end{tabular}
	\quad \\
	\quad \\
	$1$-norm\\
	\begin{tabular}{|l|l|l|l|l|l|}
		\hline
		& median & 90th & 99th & 99.9th & maximum value\\
		\hline
		$n=5$ & 1.0474 & 1.4590 & 2.4822 & 4.2311 & 6.8434\\
		\hline
		$n=25$ & 1.0687 & 1.4478 & 2.2913 & 4.6761 & 6.3929\\
		\hline
		$n=100$ & 1.0536 & 1.4274 & 2.1011 & 3.0116 & 4.5574\\
		\hline	  
	\end{tabular}
	\caption{Percentiles of the ratio $R$ in (\ref{ratioR}) for $5000$ instances of $A$ with $\Lambda_1$ simple real.}
	\label{tablereaL}
\end{table}
 
\begin{table}
	
	Euclidean Norm\\
	\begin{tabular}{|l|l|l|l|l|l|}
		\hline
		& median & 90th & 99th & 99.9th & maximum value\\
		\hline
		$n=5$ & 1.0033 & 1.1752 & 1.7349 & 3.2651 & 7.7723\\
		\hline
		$n=25$ & 1.0181 & 1.2917 & 2.0735 & 3.0539 & 5.5074\\
		\hline
		$n=100$ & 1.0164 & 1.3112 & 2.0263 & 3.2667 & 4.0842\\
		\hline	  
	\end{tabular}
	\quad \\
	\quad\\
	$\infty$-norm\\
	\begin{tabular}{|l|l|l|l|l|l|}
		\hline
		& median & 90th & 99th & 99.9th & maximum value\\
		\hline
		$n=5$ & 1.0000 & 1.2002 & 1.7929 & 3.1110 & 6.7306\\
		\hline
		$n=25$ & 1.0282 & 1.3542 & 2.1083 & 3.2852 & 5.6688\\
		\hline
		$n=100$ & 1.0248 & 1.3526 & 2.0401 & 3.4070 & 4.0496\\
		\hline	  
	\end{tabular}
	\quad \\
	\quad \\
	$1$-norm\\
	\begin{tabular}{|l|l|l|l|l|l|}
		\hline
		& median & 90th & 99th & 99.9th & maximum value\\
		\hline
		$n=5$ & 1.0016 & 1.2190 & 1.8220 & 2.7602 & 7.6305\\
		\hline
		$n=25$ & 1.0244 & 1.3308 & 2.0630 & 3.0337 & 5.1067\\
		\hline
		$n=100$ & 1.0224 & 1.3484 & 2.0166 & 3.3294 & 4.8190\\
		\hline	  
	\end{tabular}
	\caption{Percentiles of the ratio $R$ in (\ref{ratioR}) for $5000$ instances of $A$ with $\Lambda_1$ simple complex.}
	\label{tablecomplex}
\end{table}

From these experiments, we can conclude that, at least for the $p$-norm with  $p\in\{1,2,\infty\}$ as vector norm, the ratio $R$ is very rarely much larger than $1$. Indeed, it is never much larger than $1$ in these $10\ 000$ instances.

Therefore, we have the strong evidence of the following fact.

\begin{fact}\label{Fact}
Suppose that the $p$-norm with  $p\in\{1,2,\infty\}$ is used as vector norm. In the vast majority of cases for the ODE (\ref{ode}),
\begin{equation*}
\max\limits_{t\geq 0}K(t,y_0) \ \text{and} \ \max\limits_{t\geq 0}K_\infty(t,y_0)
\end{equation*}
are of the same order of magnitude. Consequently, in the vast majority of cases, the problem (\ref{due}) is well-conditioned at $y_0$ if and only if it is asymptotically well-conditioned at $y_0$.
\end{fact}

In the following two facts, we separate the cases $\Lambda_1$ simple real and $\Lambda_1$ simple complex. For the latter, recall Subsection \ref{V1c1}.

\begin{fact}
Suppose that the $p$-norm with $p\in\{1,2,\infty\}$ is used as vector norm. In the vast majority of cases of the ODE (\ref{ode}) with $\Lambda_1$ simple real,
\begin{equation}
	\max\limits_{t\geq 0}K(t,y_0) \ \text{and} \ K_\infty(y_0)=\frac{1}{\vert \widehat{w}^{(1)}\widehat{y}_0\vert} \label{KandK+}
\end{equation}
are of the same order of magnitude. Consequently, the problem (\ref{due}) is well-conditioned at $y_0$ if and only if 
\begin{equation}
	\frac{1}{\vert \widehat{w}^{(1)}\widehat{y}_0\vert}=\frac{\left\Vert w^{(1)} \right\Vert \left\Vert y_0 \right\Vert}{\vert  w^{(1)} y_0\vert }
\end{equation}
is not large.
\end{fact}

\begin{fact}
Suppose that the Euclidean norm is used as vector norm. In the vast majority of cases of the ODE (\ref{ode}) with $\Lambda_1$ simple complex,
\begin{equation*}
\max\limits_{t\geq 0}K(t,y_0),\ \ \max\limits_{t\geq 0}K_\infty(t,y_0)\text{\ \ and\ \ }\mathrm{OSF}(y_0)=\frac{1}{\vert \widehat{w}^{(1)}\widehat{y}_0\vert}
\end{equation*}
are of the same order of magnitude. Consequently, the problem (\ref{due}) is well-conditioned at $y_0$ if and only if
\begin{equation*}
	\frac{1}{\vert \widehat{w}^{(1)}\widehat{y}_0\vert}=\frac{\left\Vert w^{(1)} \right\Vert_2 \left\Vert y_0 \right\Vert_2}{\vert  w^{(1)} y_0\vert } 
\end{equation*}
is not large.
\end{fact}

\subsection{Highly non-normal matrices} \label{Hnn}

	One might get the impression that the previous conclusion regarding what happens in the vast majority of cases for the ODE (\ref{ode}) does not hold for highly non-normal matrices, since sampling the entries of the matrix $A$ from the standard normal distribution rarely produces highly non-normal matrices.
	
	However, this conclusion seems to be valid also for highly non-normal matrices, as demonstrated by the example in Subsection \ref{hnn}, where the quantities in (\ref{KandK+}) are of the same order of magnitude.  
	
	To show experimentally this, we consider a total of $5000$ instances of an ODE (\ref{ode}) of dimension $n=5,25,50$, where the initial value $y_0$ has components drawn from the standard normal distribution, and the matrix $A$ is given by $A = Q U Q^T$ with  $Q$ being the orthogonal matrix obtained from the QR factorization of a $n\times n$ random matrix with entries from the standard normal distribution; and $U$ is a $n\times n$ random upper triangular matrix with entries also drawn from the standard normal distribution. Generating matrices $A$ in this manner yields highly non-normal matrices with real eigenvalues.
	
Table \ref{tablehnn} is  analogous to Tables \ref{tablereaL} and \ref{tablecomplex}. We now observe larger values of the ratio \( R \), but they are rare. Overall, the global picture remains unchanged compared to the previous experiment.

\bigskip

	\begin{table}
		
		Euclidean Norm\\
		\begin{tabular}{|l|l|l|l|l|l|}
			\hline
			& median & 90th & 99th & 99.9th & maximum value\\
			\hline
			$n=5$ & 1.0259 & 1.7523 & 5.4707 & 15.9212 & 24.5415\\
			\hline
			$n=25$ & 1.1548 & 3.0735 & 10.3924 & 23.8093 & 64.9178\\
			\hline
			$n=100$ & 1.5150 & 4.9395 & 17.4353 & 49.1532 & 157.2478\\
			\hline	  
		\end{tabular}
		\quad \\
		\quad\\
		$\infty$-norm\\
		\begin{tabular}{|l|l|l|l|l|l|}
			\hline
			& median & 90th & 99th & 99.9th & maximum value\\
			\hline
			$n=5$ & 1.0889 & 2.0923 & 5.7941 & 17.6035 & 22.5019\\
			\hline
			$n=25$ & 1.2998 & 3.4813 & 10.8553 & 24.1466 & 53.6882\\
			\hline
			$n=100$ & 1.6693 & 5.3319 & 19.5200 & 54.3432 & 147.7216\\
			\hline	  
		\end{tabular}
		\quad \\
		\quad \\
		$1$-norm\\
		\begin{tabular}{|l|l|l|l|l|l|}
			\hline
			& median & 90th & 99th & 99.9th & maximum value\\
			\hline
			$n=5$ & 1.0666 & 1.9550 & 5.6517 & 20.7590 & 24.3172\\
			\hline
			$n=25$ & 1.2235 & 3.2105 & 10.5126 & 24.9469 & 68.3449\\
			\hline
			$n=100$ & 1.5765 & 5.0773 & 17.1928 & 60.0687 & 174.2289\\
			\hline	  
		\end{tabular}
		\caption{Percentiles of the ratio $R$ in (\ref{ratioR}) for $5000$ instances of $A=QUQ^T$.}
		\label{tablehnn}
	\end{table}

\section{Non-normal dynamics} \label{NND}

Suppose we want to simulate the transient phase of a real-world system, where this transient $y$ satisfies an ODE (\ref{ode}) with a non-normal stable matrix $A$, where \emph{stable} means that all the eigenvalues of $A$ have negative real part. For a non-linear real-world system, this ODE is the linearization around  an asymptotically stable equilibrium of a non-linear ODE. 

Due to the non-normality of the matrix $A$, it is expected that the transient exhibits an initial (possibly large) growth before decaying to zero. This phenomenon is what characterizes non-normal dynamics and can destabilize a non-linear system subject to small perturbations from the asymptotically stable equilibrium (see \cite{Trefethen2005} and \cite{Asilani2018}).

Suppose that, through simulation, we need to determine the maximum growth of the transient, i.e. $\max\limits_{t\geq 0} \Vert y(t) \Vert$. Additionally, suppose there is uncertainty in the initial value of our simulation. Here, we consider $y_0$ as the initial value available for the simulation and $\widetilde{y}_0$ as the actual initial value.

Since there are uncertainties in the initial value, it is important to assess how close the simulation result $\max\limits_{t\geq 0} \Vert y(t) \Vert$ is to the actual maximum growth of the transient, i.e. $\max\limits_{t\geq 0} \Vert \widetilde{y}(t) \Vert$.

Information on relative closeness can be obtained by the next result.
\begin{theorem}
Suppose that the initial value of an ODE (\ref{ode}) with $A$ stable is perturbed by a normwise relative error $\varepsilon$. We have
\begin{equation}
\left\vert \frac{\max\limits_{t\geq 0} \Vert \widetilde{y}(t) \Vert}{\max\limits_{t\geq 0} \Vert y(t) \Vert}-1\right\vert\leq E:=\max\limits_{t\geq 0}K(t,y_0)\cdot \varepsilon, \label{E}
\end{equation}
whenever $E<1$.
\end{theorem}
\begin{proof}
Since
$$
\Vert \widetilde{y}(t)\Vert =\Vert y(t)\Vert(1+\chi(t)),
$$
where
$$
\vert \chi(t)\vert \leq \delta(t)=K(t,y_0,\widehat{z}_0)\varepsilon\leq K(t,y_0)\varepsilon\leq E
$$
(see (\ref{magnification})), we have
\begin{equation}
\max\limits_{t\geq 0} \Vert \widetilde{y}(t) \Vert\leq \max\limits_{t\geq 0} \Vert y(t) \Vert(1+E). \label{E1}
\end{equation}
Moreover, since 
$$
\Vert y(t)\Vert =\frac{\Vert \widetilde{y}(t)\Vert}{1+\chi(t)},
$$
we have
\begin{equation}
\max\limits_{t\geq 0} \Vert y(t) \Vert\leq \frac{\max\limits_{t\geq 0} \Vert \widetilde{y}(t)\Vert}{1-E}, \label{E2}
\end{equation}
whenever $E<1$. By (\ref{E1}) and (\ref{E2}), (\ref{E}) follows.
\end{proof}
The quantity $E$ is an upper bound for the relative error of the actual maximum transient growth with respect to the simulated maximum transient growth, demonstrating how the condition number $K(t,y_0)$ can be useful in non-normal dynamics simulations.
 
By recalling Fact \ref{Fact}, we can conclude what follows.

 \begin{fact}\label{FactII}
 	Suppose that the $p$-norm with $p\in\{1,2,\infty\}$ is used as vector norm. In the vast majority of cases of the ODE (\ref{ode}) with $A$ stable, $E$ in (\ref{E}) and
\begin{equation*}
E_\infty:=\max\limits_{t\geq 0}K_\infty(t,y_0)\cdot \varepsilon 
\end{equation*}
are of the same order of magnitude.
\end{fact}
Therefore, also the asymptotic condition number $K_\infty(t,y_0)$ can be useful in non-normal dynamics simulations.

\begin{example}
	Consider an ODE (\ref{ode}) with
	\begin{equation}
A=\left[ 
\begin{array}{rr}
	-1 & a \\ 
	0 & -2 
\end{array}
\right], \label{Annormal}
\end{equation}
where $a>0$ is large. In this example,  the $\infty$-norm is used as the vector norm.

The matrix $A$ represents a classical example of a highly non-normal stable matrix (see the ``quiz'' at the beginning of \cite{Trefethen2005}), exhibiting a large initial  growth of the solution of (\ref{ode}). In Figure \ref{three-expm}, for $a=50,500,5000$, we see that $\Vert \mathrm{e}^{tA}\Vert_\infty$ initially grows significantly over time before eventually decaying to zero.

Since  
$$  
\left\vert \max_{t\geq 0} \Vert \widetilde{y}(t) \Vert_\infty - \max_{t\geq 0} \Vert y(t) \Vert_\infty \right\vert \leq \max_{t\geq 0} \Vert \mathrm{e}^{tA} \Vert_\infty \Vert \widetilde{y}_0 - y_0 \Vert_\infty,  
$$  
it is expected that $\max\limits_{t \geq 0} \Vert \widetilde{y}(t) \Vert_\infty$ will not be as close to $\max\limits_{t \geq 0} \Vert y(t) \Vert_\infty$ as $\widetilde{y}_0$ is close to $y_0$, when the closeness is measured in terms of absolute error.

On the other hand, as we will see below, this is not the case when the closeness is measured in terms of relative error. Therefore, if controlling the relative error of $\max\limits_{t \geq 0} \Vert \widetilde{y}(t) \Vert_\infty$ is important, considering the growth of $\Vert \mathrm{e}^{tA}\Vert_\infty$ can be misleading, as it pertains to the absolute error rather than the relative error. We need to consider $K(t,y_0)$, not $\Vert \mathrm{e}^{tA}\Vert_\infty$.

\begin{figure}[tbp]
	\centering
	\par
	\includegraphics[width=1\textwidth]{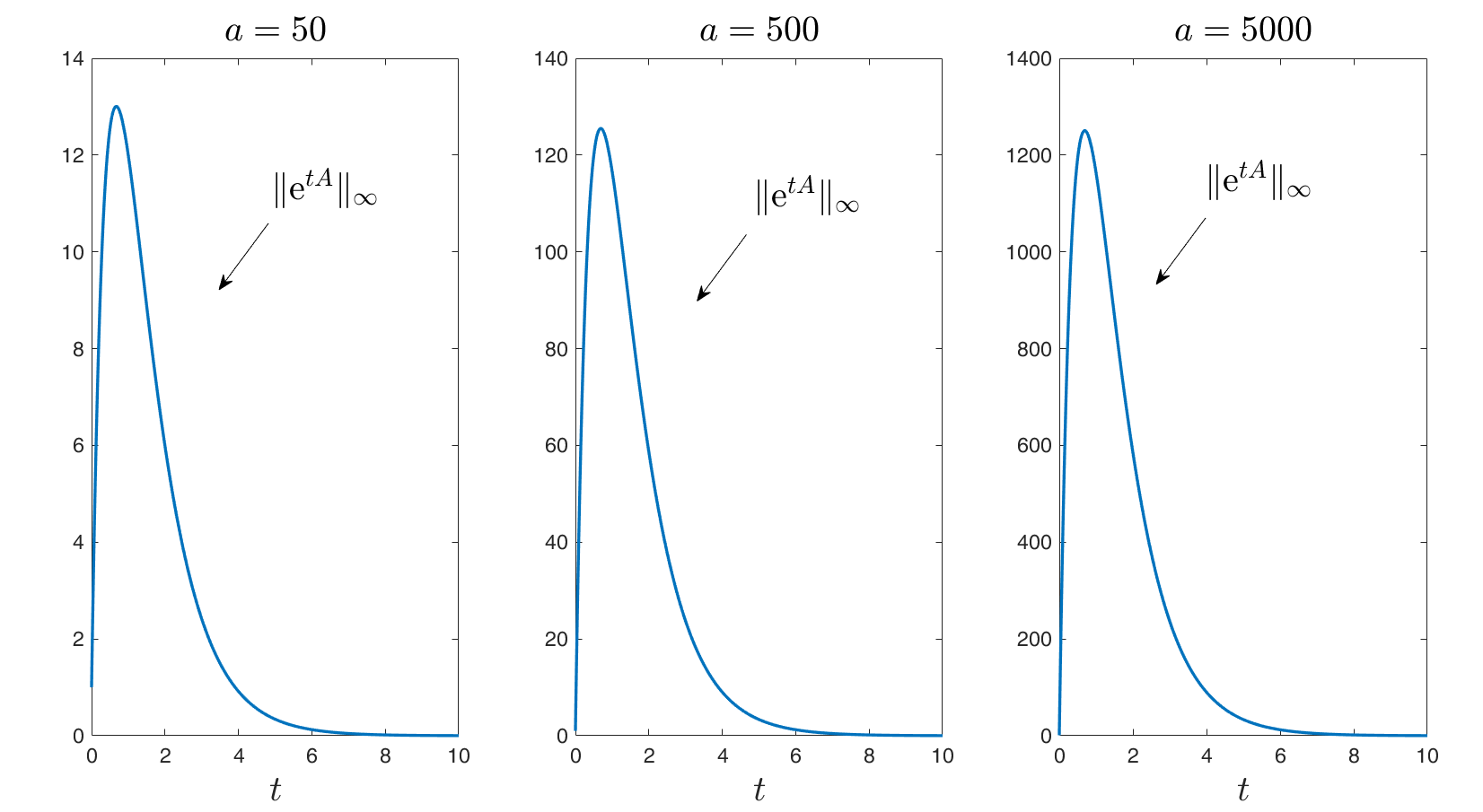}
	\caption{$\Vert \mathrm{e}^{tA}\Vert_\infty$, $t\in[0,10]$, for $a=50,500,5000$ in (\ref{Annormal}).}
	\label{three-expm}
\end{figure}

For $a=50,500,5000$, we present in the table below percentiles of the ratio
$$
R=\frac{\max\limits_{t\in[0,T]}K(t,y_0)}{K_\infty(y_0)}
$$
computed for $10\ 000$ instances of $y_0$ sampled from the standard normal distribution. As in Section \ref{NAsy}, we set $T$ to fifty times the characteristic time, i.e. $T=50$, and the maximum is determined using $1000$ equally spaced sampling points.
$$
\begin{tabular}{|l|l|l|l|l|l|}
	\hline
	& \emph{median} & \emph{90th} & \emph{99th} & \emph{99.9th} &\emph{maximum value}\\
	\hline
	$a=50$ & $1.9125$ &  $4.6042$ &  $20.4262$ &  $45.3000$ & $49.0576$   \\
	\hline
	$a=500$ & $1.3983$ & $3.4958$ & $27.7753$  & $179.9569$  & $364.8947$  \\
	\hline
	$a=5000$  & $1.0400$    & $1.1898$  & $8.8109$ & $79.8017$ & $1206.5313$ \\
	\hline	  
\end{tabular} 
$$
This table confirms what is expected (see Fact \ref{FactII}): in the vast majority of cases for $y_0$, the maxima in (\ref{KandK+}) are of the same order of magnitude.

We explain this example in more detail. The left eigenvector corresponding to the rightmost eigenvalue $-1$ is $$w^{(1)}=\left[\begin{array}{rr}
1 & a 
\end{array}
\right].$$ 
Therefore,
\begin{equation}
E_\infty=K_\infty(y_0)\cdot \varepsilon =\frac{1+a}{\vert \widehat{y}_{01}+a\widehat{y}_{02}\vert}\cdot\varepsilon. \label{Einftynw}
\end{equation}
In the vast majority of cases for $y_0$, the relative error of $\max\limits_{t \geq 0} \Vert \widetilde{y}(t) \Vert_\infty$ with respect to $\max\limits_{t \geq 0} \Vert y(t) \Vert_\infty$ is not much larger than $E_\infty$.

The expression in (\ref{Einftynw}) involves the normwise relative error $\varepsilon$ of $\widetilde{y}_0$. However, by recalling point 1) in Subsection \ref{nwcw}, we have
$
\varepsilon\leq \max\{\varepsilon_1,\varepsilon_2\},
$
where $\varepsilon_1$ and $\varepsilon_2$ are the componentwise relative errors of $\widetilde{y}_0$. Therefore, we have that, in the vast majority of cases for $y_0$, the relative error of $\max\limits_{t \geq 0} \Vert \widetilde{y}(t) \Vert_\infty$ with respect to $\max\limits_{t \geq 0} \Vert y(t) \Vert_\infty$ is not much larger than
$$
\frac{1+a}{\vert \widehat{y}_{01}+a\widehat{y}_{02}\vert}\cdot \max\{\varepsilon_1,\varepsilon_2\},
$$

To confirm this last conclusion, consider $y_0 =\widehat{y}_0= (1,1)$. It is expected  that the relative error of $\max\limits_{t \geq 0} \Vert \widetilde{y}(t) \Vert_\infty$ with respect to $\max\limits_{t \geq 0} \Vert y(t) \Vert_\infty$ is not much larger than $\max\{\varepsilon_1,\varepsilon_2\}$. Indeed, for $\widetilde{y}_0 = (1.01,0.99)$, the following table presents the signed absolute and relative errors of  
\(\max\limits_{t \geq 0} \Vert \widetilde{y}(t) \Vert_\infty\) with respect to $\max\limits_{t \geq 0} \Vert y(t) \Vert_\infty$. 
$$  
\begin{tabular}{|l|l|l|}  
	\hline  
	& \emph{Absolute error} & \emph{Relative error} \\  
	\hline  
	$a=50$ & $-0.1198$ &  $-0.0092$ \\  
	\hline  
	$a=500$ & $-1.2450$ & $-0.0099$ \\  
	\hline  
	$a=5000$  & $-12.4949$    & $-0.0100$ \\  
	\hline	  
\end{tabular}  
$$  
While the absolute error grows significantly with respect to the absolute errors $\pm 0.01$ of the components of \(\widetilde{y}_0\), due to the large values attained by \(\Vert \mathrm{e}^{tA} \Vert_\infty\), the relative error remains close to the relative errors $\pm 0.01$ of the components of \(\widetilde{y}_0\).
\end{example}
\section{Conclusion}

The present paper explores how a perturbation of the initial value of the ODE (\ref{ode}) is propagated to the solution by examining the relative error, i.e., by comparing the perturbation of the solution to the solution itself. Considering the relative error, rather than the absolute error, provides a clearer perspective on the perturbation as the solution $y(t)$ evolves with time $t$: at any given time $t$, the perturbation of the solution $y(t)$ can be large (small) when compared to the initial value $y_0$, and simultaneously small (large) when compared to $y(t)$.

Understanding the behavior of the relative error of the perturbation can offer new insights and impact our approach to understanding the propagation of uncertainties. For example, the non-normality of the matrix $A$ appears to have a lesser effect on the relative error compared to the absolute error.

When we examine the \emph{absolute conditioning} of the problem (\ref{due}), the pointwise condition number is
$$
K_{\mathrm{abs}}(t)=\Vert \mathrm{e}^{tA}\Vert.
$$
Therefore, in the context of absolute error, a perturbation of the initial value grows in the worst-case scenario as a solution of $y'(t)=Ay(t)$ grows in the worst-case scenario. There is nothing new: by knowing how the solutions grow, we also know how the perturbations grow, since the perturbations are themselves solutions. Moreover, the initial value $y_0$ does not play any role, since the pointwise condition number is independent of $y_0$.

On the other hand, when we consider the \emph{relative conditioning} of the problem (\ref{due}), the pointwise condition number is
$$
K_{\mathrm{rel}}(t,y_0)=K(t,y_0)=\frac{\Vert \mathrm{e}^{tA}\Vert \Vert y_0\Vert }{\Vert \mathrm{e}^{tA}y_0\Vert}.
$$
In the context of relative error, the perturbations do not grow in the same manner as the solutions, but in a completely different way: compare fact A) at page 1  with fact B) at page 17. Unlike the absolute conditioning, the initial value $y_0$ plays a role in the relative conditioning. Indeed, the time growth of $\Vert \mathrm{e}^{tA}\Vert$ can be mitigated by  the time growth of $\Vert \mathrm{e}^{tA}y_0\Vert$. 

These considerations represent the novelties in focusing on relative error rather than absolute error, offering new perspectives in linear dynamics. In particular, when simulating the transient behavior and there is uncertainty in the initial value of the transient: see Subsection \ref{buildingheating} in the context of determining the end of the transient, and Section \ref{NND}, concerning non-normal dynamics, in the context of determining the maximum transient growth.

The present paper analyzes the asymptotic (long-time) behavior of the relative error of the perturbed solution. This asymptotic behavior also provides information on the non-asymptotic behavior. In fact, the strong experimental evidence included in the paper suggests that, in the vast majority of cases, the maximum values over time $t$   of the condition number $K(t,y_0)$ and its asymptotic form $K_\infty(t,y_0)$ share the same order of magnitude. 

The main result of the paper can be summarized as follows. Consider the Euclidean norm as vector norm. In the vast majority of cases for the ODE (\ref{ode}), the relative error of the perturbed initial value is magnified in the perturbed solution, in the worst-case scenario for the initial perturbation, by a factor of the order of magnitude of
\begin{equation}
\frac{\left\Vert w^{(1)} \right\Vert_2 \left\Vert y_0 \right\Vert_2}{\vert  w^{(1)} y_0\vert }, \label{final1}
\end{equation}
where the row vector $w^{(1)}$ is a left eigenvector of the matrix $A$ corresponding to the rightmost real or complex eigenvalue. Indeed, the asymptotic magnification factor is of the same order of magnitude as the quantity in (\ref{final1}).

\bigskip

\noindent {\bf Acknowledgements:} the research was supported by the INdAM
Research group GNCS (Gruppo Nazionale di Calcolo Scientifico).

\end{document}